\documentclass[11pt,reqno]{amsart}
\usepackage[margin=1in]{geometry}
\usepackage{amsmath,amssymb,amsthm,graphicx,amsxtra, setspace}
\usepackage[utf8]{inputenc}
\usepackage{mathrsfs}
\usepackage{hyperref}
\usepackage{upgreek}
\usepackage{mathtools}
\usepackage{xcolor}
\usepackage{mathrsfs}
\usepackage{calrsfs} 
\usepackage[T1]{fontenc}  
\usepackage[utf8]{inputenc}
\usepackage[mathcal]{euscript}
\allowdisplaybreaks

\usepackage[pagewise]{lineno}

\DeclareMathAlphabet{\mathpzc}{OT1}{pzc}{m}{it}

\usepackage[cyr]{aeguill}

\colorlet{darkblue}{blue!50!black}

\hypersetup{
	colorlinks,%
	citecolor=blue,%
	filecolor=red,%
	linkcolor=red,%
	urlcolor=blue,%
	pdfnewwindow=true,%
	pdfstartview={FitH}
}

\newtheorem{theorem}{Theorem}[section]
\newtheorem{lemma}[theorem]{Lemma}

\newtheorem{corollary}[theorem]{Corollary}
\newtheorem{definition}[theorem]{Definition}

\newtheorem{ex}[theorem]{Example}
\newtheorem{remark}[theorem]{Remark}

\allowdisplaybreaks

\let\originalleft\left
\let\originalright\right
\renewcommand{\left}{\mathopen{}\mathclose\bgroup\originalleft}
\renewcommand{\right}{\aftergroup\egroup\originalright}

\newcommand{\norm}[1]{\left\|#1\right\|}

\newcommand{\Addresses}{{
		\footnote{
			
			\noindent \textsuperscript{1,3}Department of Mathematics, Indian Institute of Technology Roorkee-IIT Roorkee,
			Haridwar Highway, Roorkee, Uttarakhand 247667, INDIA.\par\nopagebreak
			\noindent \textsuperscript{2}Universidad de Talca, Instituto de Matemáticas, Casilla 747, Talca, Chile.\par\nopagebreak
			\noindent  \textit{e-mail:} \texttt{Manil T. Mohan: maniltmohan@ma.iitr.ac.in, maniltmohan@gmail.com.}
			
			\textit{e-mail:} \texttt{Jyoti Jindal: jyoti@ma.iitr.ac.in.}
			
			\textit{e-mail:} \texttt{Sumit Arora: arorasumit10623@gmail.com, sumit.arora@utalca.cl.}
			
			\noindent \textsuperscript{*}Corresponding author.
			
			\textit{Key words:} Approximate controllability, Hemivariational inequalities, Differential inclusion, Resolvent operator.
			
}}}

\begin{document}
	
	\title[Approximate Controllability of Neutral Integro-Differential Inclusions]{Approximate Controllability of Neutral Integro-Differential Inclusions Associated with  Hemivariational Inequalities in Banach Spaces with Uniformly Convex Dual
		\Addresses}
	\author[J. Jindal, S. Arora, and M. T. Mohan]
	{Jyoti Jindal\textsuperscript{1}, Sumit Arora\textsuperscript{2}, and Manil T. Mohan\textsuperscript{3*}}
	
	\maketitle
	\begin{abstract}
		In this paper, we investigate the approximate controllability of neutral integro-differential inclusions driven by hemivariational inequalities with nonlocal inclusion conditions in a Banach space $\mathbb{V}$ whose dual $\mathbb{V}^*$ is uniformly convex, using resolvent operator theory. We first establish an equivalence between the original problem and an associated differential inclusion involving the Clarke subdifferential. Under suitable assumptions, we then prove the existence of mild solutions corresponding to admissible controls in $\mathrm{L}^2(I,\mathbb{U})$. A key contribution of this work lies in the analysis of the control construction in Banach spaces, where the duality mapping $\mathcal{J} = \partial\big(\tfrac{1}{2}\|\cdot\|_{\mathbb{V}}^2\big)$ is generally nonlinear, in contrast to the Hilbert space setting, while such nonlinearity has been viewed as a technical obstruction in earlier studies. Based on this framework, we prove approximate controllability of the original system by relating it to the controllability of the associated linearized inclusion and constructing appropriate controls involving the duality mapping. Finally, we present applications to illustrate the effectiveness of the proposed approach.
	\end{abstract}
	
	\section{Introduction}
	
	\subsection{Literature survey}
	Variational and hemivariational inequalities constitute two fundamental classes of nonlinear problems that naturally emerge in the mathematical modeling of equilibrium and optimization phenomena. Variational inequalities are well-suited for systems governed by convex energy functionals and find extensive application in fields such as contact mechanics \cite{VB}, economic equilibrium theory \cite{ZNMN}, and fluid flow \cite{JTONK}. Hemivariational inequalities extend this paradigm to non-convex and non-smooth regimes, offering a framework capable of handling multi-valued and non-differentiable energy functionals. As a result, they are well-suited for describing phenomena such as frictional effects and material nonlinearities. These multi-valued relationships are derived from non-smooth and non-convex superpotentials employing Clarke’s generalized gradient.  When the superpotential is convex, the hemivariational inequalities reduce to standard variational inequalities.
	
	Taken together, these two frameworks provide powerful analytical tools for addressing a wide range of problems in mechanics, physics, engineering, and economics. Among the foundational milestones in the theory of variational inequalities was Fichera’s groundbreaking work \cite{GFi} on frictionless contact between a linearly elastic body and a rigid foundation, a problem originally formulated by Signorini \cite{ASi}. Interest in hemivariational inequalities emerged in a similar manner, driven primarily by mechanical applications. The theory of hemivariational inequalities was first formulated and analyzed by Panagiotopoulos during the early 1980s (see \cite{PDPH}). By rigorously treating contact conditions through inequalities rather than classical equalities, this contribution laid a crucial foundation for modern analysis. Readers interested in the historical development and underlying mathematical theory of hemivariational inequalities are encouraged to consult \cite{ZNPDP, PDPH}.
	
	As is well known, controllability plays a vital role in system analysis and design of control systems. In infinite-dimensional settings, two fundamental notions of controllability, namely exact controllability and approximate controllability, are naturally distinguished, and each gives rise to a broad range of applications. However, exact controllability is difficult to achieve in the case of infinite-dimensional control systems (see \cite{NIM, RT, EZ1}). This motivates the study of the weaker and more practical notion of controllability, known as approximate controllability, which means that the system can steer its state arbitrarily close to a desired final state. The authors in \cite{ZLLX, ZLXLD} analyzed the existence of mild solutions and approximate controllability for control problems driven by nonlinear evolution hemivariational inequalities in Hilbert spaces. 
	Moreover, boundary value problems for first- and second-order (functional) differential equations have been widely studied. Notably, Byszewski and Lakshmikantham \cite{LBVL} showed that many practical applications require conditions imposed throughout the whole time interval, thereby including boundary value conditions. Such conditions are known in the literature as non-local conditions \cite{LB2, LB4}. Consequently, the study of differential equations with non-local initial conditions has become an active area of research (see \cite{ACDNP, KEXFKH}). Additionally, numerous studies have investigated the controllability of systems governed by non-local conditions (see, e.g., \cite{FGHCZ, KJDB} and references therein).
	
	In the context of optimal control of evolution equations involving hemivariational inequalities, the recent papers \cite{SV-2025-2,SV-2026} are formulated in the framework of feedback control systems; a careful analysis reveals that they do not provide a constructive synthesis of a feedback control law in the classical sense. In particular, the control is defined through a multivalued map 
	$\mathcal{U}(\cdot,q(\cdot))$, where $q$ is the state variable and the main results establish the existence of a measurable selection $u(t)$ satisfying $u(t)\in \mathcal{U}(t,q(t)) $ a.e. $t \in [0,a]$. The use of compactness arguments, the Filippov theorem, and the Cesari property ensures the existence of such admissible controls, which is an open-loop construction, as part of feasible pairs $(q,u)$. However, this approach neither produces an explicit feedback operator of the form $u(t)=\mathrm{K}(t,q(t))$ nor offers a constructive method for determining the control in terms of the current state. Therefore, the results should be interpreted as establishing the existence of admissible state-dependent controls rather than the synthesis of a true feedback control law. This distinction is important, particularly in stochastic and infinite-dimensional settings, where adaptedness, measurability, and implementability of the control play a crucial role.
	
	The term neutral in the theory of functional differential equations has a precise and well-established meaning, originating from the classical works of Hale \cite{JKL} and Hale-Verduyn Lunel \cite{JKLSMV}. A differential equation is classified as neutral if the highest-order derivative of the unknown function appears implicitly through a functional that depends on the history of the state, rather than solely on its present value. The canonical form is given by
	\begin{align*}
		\frac{d}{dt}[q(t)+\mathrm{g}(t,q_t)]=\mathrm{F}(t,q_t),\ t\in [0,a],
	\end{align*}
	where $q_t(s)=q(t+s),s\in[-r,0],$ denotes the history segment and $\mathrm{g}, \mathrm{F} :[0,a]\times \mathrm{C}([-r,0];\mathbb{V})\to\mathbb{V}$ are functions acting on the history. A crucial requirement is that the function $\mathrm{g}$ depends on the past trajectory $q_t$, not merely on the current state $q(t)$, so that the derivative cannot be separated to act only on $q(t)$, but is intrinsically coupled with the system’s history. Neutral fractional differential systems have attracted considerable attention due to their ability to model hereditary effects arising from such history-dependent structures. However, a closer examination of several recent works, including \cite{SV-2024, SV-2025} reveals that the neutral character of the equations considered therein requires further clarification. In these works, the function $\mathrm{g}(t,q(t))$ appearing inside the derivative depends only on the current state $q(t)$, as reflected in the imposed growth conditions (e.g., assumption (H4) in \cite{SV-2025} and analogous hypotheses in related works), with no dependence on the history segment $q_t$ or any delayed argument. Consequently, it reduces to standard differential evolution inclusions with nonlinear perturbations, and the defining structural features of neutral systems, such as dependence on past trajectories and the requirement to prescribe initial data over a history interval, do not arise in the analysis. Therefore, it would enhance the mathematical precision of these works to either justify the neutral terminology through assumptions that genuinely incorporate history dependence or to reclassify the models as non-neutral differential stochastic inclusions with nonlinear perturbations, which more accurately reflect the underlying structure of the problems studied.
	
	In \cite{vijayakumar2018approximate}, the authors utilized the resolvent operator method to demonstrate the approximate controllability results of abstract neutral integro-differential inclusions (NIDI) with infinite delay in Hilbert spaces. For further studies, 
	we refer the reader to \cite{arora2025controllability, mokkedem2014approximate} and the references therein. Recently, Pradeesh et al. \cite{JPSVTA} studied Sobolev-type neutral integro-differential evolution hemivariational inequalities in Hilbert spaces and developed results related to approximate controllability. Furthermore, Kattan and Hammad \cite{DAKHA} investigated optimal and approximate control strategies for a class of neutral integro-differential systems governed by hemivariational inequalities in Hilbert spaces. It is worth noting that establishing controllability results in Banach spaces is more challenging than in Hilbert spaces. Subsequently, the authors in \cite{arora2022existence, kumbhakar2023approximate} investigated the approximate controllability of semilinear differential inclusions in reflexive Banach spaces using resolvent operator condition and fixed point techniques. In \cite{JJBKDNP, BKDNP}, the authors recently established the approximate controllability for evolution hemivariational inequalities in super-reflexive Banach spaces and overcame the convexity issues arising from the non-linear nature of the duality mapping in the control representation. To the best of the authors’ knowledge, existing literature does not address the approximate controllability results for abstract NIDI associated with hemivariational inequalities and non-local inclusion conditions in Banach spaces. Thus, this manuscript aims to fill this research gap. 
	
	\subsection{The model problem}
	Motivated by the above-mentioned works, the aim of this manuscript is to develop the existence and approximate controllability of the following hemivariational inequality problem:
	\begin{equation}\label{HV1}
		\begin{cases}
			\begin{aligned}
				&\left\langle -\frac{d}{dt}\left(q(t)+\int_0^t N(t-s)q(s)ds\right)+\mathscr{A}q(t)+\int_0^t G(t-s)q(s)ds +Bu(t), v^*\right\rangle_{\mathbb{V},\mathbb{V}^*} \\&\quad+ \mathbb{F}^0(t, q(t); H^*v^*)\ge 0\ \text{ for all }\  v^*\in\mathbb{V}^*,~ t\in(0,a],\\
				&q(0)\in q_0+h(q),
			\end{aligned}
		\end{cases}
	\end{equation}
	where $\mathbb{V}$ is a Banach space with uniformly convex dual, the notation $\langle \cdot, \cdot \rangle_{\mathbb{V}, \mathbb{V}^*}$ denotes the duality pairing between $\mathbb{V}$ and its dual space $\mathbb{V}^*$, $ \mathscr{A} : \mathscr{D}(\mathscr{A}) \subset \mathbb{V}\to \mathbb{V}$ and $G(t) \colon \mathscr{D}(G(t)) \subset \mathbb{V} \to \mathbb{V} $ for $t \geq 0$ are closed linear operators and $I=[0,a]$. The linear operators $N(t) \colon \mathbb{V}\to \mathbb{V}, t \geq 0$ and $B \colon \mathbb{U}\to \mathbb{V}$ are bounded. In addition, the control function $u(\cdot) \in \mathrm{L}^2(I, \mathbb{U})$, where $\mathbb{U}$ is a separable Hilbert space identified with its dual and the notation $\mathbb{F}^0(t, \cdot;\cdot)$ stands for the generalized Clarke directional derivative of a locally Lipschitz function $\mathbb{F}(t,\cdot): \mathbb{V}\to \mathbb{R}$, which is defined later. The linear operator $H: \mathbb{V}^*\to \mathbb{V}$ is bounded and $h : \mathrm{C}(I, \mathbb{V}) \to \mathcal{P}(\mathbb{V})$ is a non-local condition. The main results of this work are listed below: 
	\begin{enumerate}
		\item Under the assumptions \textit{(F1)}-\textit{(F2)} on $\mathbb{F}$ and \textit{(H1)}–\textit{(H3)} on $h$ (see Section \ref{sec3}), we first establish the existence of mild solutions to problem \eqref{PC1} by assuming the compactness of the resolvent family $(\mathscr{R}(t))_{t>0}$ introduced in Section \ref{sec-resolvent} (cf. Theorem \ref{Theorem-2}). The main idea of the proof relies on the application of a multivalued fixed point theorem (see Theorem \ref{Fixed-4}) combined with compactness arguments. This approach effectively handles the multivalued nature arising from the Clarke subdifferential. 
		
		Further, under the same compactness assumption on the resolvent family, we establish an additional existence result for mild solutions to problem \eqref{PC1} in Theorem \ref{Theorem-1}. This is achieved by applying Theorem \ref{Fixed-1} under the conditions \textit{(F1)}-\textit{(F2)} and \textit{(F)} on the Clarke subdifferential of $\mathbb{F}$, together with assumption \textit{(H4)} and the Lipchitz condition \eqref{haus-h} on the set valued function $h$.
		
		In addition, we obtain the approximate controllability results under more relaxed conditions, that is, \textit{(F1)}, \textit{(F2)}$^\prime$ supplemented by condition \textit{(H1)}-\textit{(H3)} and \textit{(H6)} with \textit{(R)} condition as presented in Theorem \ref{MR}.
		\item  We also establish existence results without assuming the compactness of the resolvent family $(\mathscr{R}(t))_{t \geq 0}$. In this case, under conditions \textit{(F1)}-\textit{(F2)} and \textit{(F)}$^\prime$ on $\mathbb{F}$ and assumptions \textit{(H1)}, \textit{(H5)} on $h$, the result is obtained by applying fixed point theorems for condensing multivalued maps (see Theorem \ref{Theorem-3}).
	\end{enumerate}
	
	\subsection{Novelties, difficulties and approaches}
	Upon conducting an extensive literature review, 
	our results extend the earlier works in the following directions:
	\begin{itemize}
		\item[(1)] Unlike previous studies that rely on a separable Hilbert space setting, our approach extends the current framework to a general Banach space $\mathbb{V}$ with a uniformly convex dual $\mathbb{V}^*$. Consequently, the geometric structure of $\mathbb{V}$ plays a central role, and our results provide a natural generalization of those presented in \cite{DAKHA, vijayakumar2018approximate}. Furthermore, the techniques developed in this paper extend the results in \cite{JPSVTA} to Banach spaces with a uniformly convex dual.
		\item[(2)] The present work overcomes the difficulties that arise when $\mathbb{V}$ is a Banach space, in which case the duality mapping appearing in the control representation \eqref{Y2} is generally nonlinear. In particular, the image of a convex selection set under the corresponding controlled trajectory map $\Gamma_z^\varepsilon(\cdot)$ defined in \eqref{CM} need not be convex, which creates a difficulty in applying multivalued fixed point theorems requiring convex values. This particular difficulty does not arise in the Hilbert space setting, where the duality mapping reduces to the identity under the standard identification of the space with its dual. Related difficulties caused by the nonlinear duality mapping have also been discussed in previous studies (see, for instance, \cite{kumbhakar2023approximate, BKDNP}). To overcome this issue, we formulate the fixed point problem in the space of selection sets rather than in the space of trajectories and introduce the multimap $\Upsilon: Q_{\mathbb{F}}\multimap \mathrm{L}^1(I,\mathbb{V}^*)$ defined by
		$$
		\Upsilon(f)=S_{\partial\mathbb{F}}\bigl(\Gamma_z^\varepsilon(f)\bigr),
		\ \  f\in Q_{\mathbb{F}},
		$$
		where $\Gamma_z^\varepsilon$ and $S_{\partial\mathbb{F}}$ are defined in \eqref{CM} and \eqref{def-S}, respectively. For each fixed $f$, $\Upsilon(f)$ consists of the selections of the Clarke subdifferential along the fixed trajectory $\Gamma_z^\varepsilon(f)$ for $\varepsilon>0$. Therefore, \cite[Lemma 5.3]{SMAOMS} ensures that $\Upsilon(f)$ is nonempty, convex, and weakly compact. Moreover, a fixed point $f\in\Upsilon(f)$ yields $f\in S_{\partial\mathbb{F}}(q^{\varepsilon})$ with $q^\varepsilon=\Gamma_z^\varepsilon(f)$ for $\varepsilon>0$, as required by the mild solution formulation. Thus, this construction supplies the convex-valued multimap needed for the fixed point argument without requiring $\Gamma_z^\varepsilon(\cdot)$ to preserve convexity.
		
		
		\item[(3)] The present work further distinguishes itself by establishing existence results without relying on the compactness of the resolvent family. More precisely, we derive results both under the compactness assumption and without compactness of the resolvent family. This unified treatment broadens the class of admissible systems and relaxes a restrictive condition, that is, compactness of the resolvent family, which is commonly imposed in the existing literature (see, for instance, \cite{arora2025controllability, JPSVTA}).
		
		\item[(4)] The applicability of the developed results provides several new insights into the underlying framework. In particular, we present a class of transport equations generating a non-compact resolvent family (see Example \ref{ex 3.11}). We also study a controlled heat equation with Dirichlet boundary conditions and a nonlocal initial condition in the Banach space $\mathrm{L}^p([0,\pi],\mathbb{R})$ (see Example \ref{5.1111}), where the approximate controllability of the associated linear system is established through the semigroup representation and controllability operator. Furthermore, a class of heat equations with weakly singular kernels is considered in Example \ref{5.111}, illustrating the broad applicability of the obtained results.
		
		In addition, we present an example involving the Walsh--Laplacian (dyadic Laplacian) in the Banach space $\mathrm{L}^p([0,1],\mathbb{R})$, formulated through the Haar basis and dyadic harmonic analysis (see Example \ref{Example 5.6}). Using our theoretical findings, we investigate the approximate controllability of the considered system.
	\end{itemize}
	\subsection{Organization of the paper}
	This manuscript is arranged in the following order: In Section \ref{sec2}, preliminary results and basic definitions are provided, which will
	be helpful in proving the main results and also provide the connection between the original problem \eqref{HV1} and an associated differential inclusion problem \eqref{PC1} which involves the Clarke subdifferential. 
	
	Section \ref{sec3} is devoted to establish existence results for the proposed inclusion system \eqref{PC1} under the assumption that the resolvent family $(\mathscr{R}(t))_{t>0}$ is compact, considering various conditions on $\mathbb{F}$ and the non-local term $h$ (see Theorems \ref{Theorem-2} and \ref{Theorem-1}). Subsequently, we address the existence of mild solutions in the case of a non-compact resolvent family $(\mathscr{R}(t))_{t \ge 0}$ (cf. Theorem \ref{Theorem-3}).
	
	In Section \ref{sec4}, we investigate the approximate controllability of the inclusion system \eqref{PC1} under the assumption that the resolvent family $(\mathscr{R}(t))_{t>0}$ is compact.
	
	Finally, Section \ref{sec5} is devoted to illustrative examples. In particular, Example \ref{ex 3.11} demonstrates the existence of a non-compact resolvent family $(\mathscr{R}(t))_{t \geq 0}$ and verifies all the hypotheses of Theorem \ref{Theorem-3} to ensure the existence of a mild solution. Moreover, we present an example of a heat equation involving a weakly singular kernel (see Example \ref{5.111}), which is used to illustrate the obtained approximate controllability results and highlight the applicability of the developed theory.
	
	\section{Preliminaries}\label{sec2}\setcounter{equation}{0}
	In this section, we introduce fundamental notations and definitions used throughout the manuscript. Let $\mathbb{V}$ and $\mathbb{W}$ be Banach spaces. We denote by $\mathcal{L}(\mathbb{V},\mathbb{W})$ the space of bounded linear operators from $\mathbb{V}$ into $\mathbb{W}$, writing $\mathcal{L}(\mathbb{V})$ when $\mathbb{V} = \mathbb{W}$. For a closed linear operator $P \colon \mathscr{D}(P) \subseteq \mathbb{V}\to \mathbb{W}$, $\mathscr{D}_1(P)$ represents its domain equipped with the graph norm $\Vert{}z\Vert{}_1 = \Vert{}z\Vert{}_{\mathbb{V}} + \Vert{}Pz\Vert{}_{\mathbb{W}}$.
	Additionally, $\mathbb{V}_w$ denotes the Banach space $\mathbb{V}$ with the weak topology and
	let $\mathrm{C}(I,\mathbb{V})$ be the space of continuous functions $u \colon I \to \mathbb{V}$, endowed with supremum norm
	\begin{align*}
		\|u\|_{\infty}:=	\|u\|_{\mathrm{C}(I,\mathbb{V})}=\sup_{t\in I}\|u(t)\|_{\mathbb{V}}.
	\end{align*}
	For $1 \leq p < \infty$, let $\mathrm{L}^p(I,\mathbb{V})$ denote the Bochner space of strongly measurable functions $u \colon I \to \mathbb{V}$ such that
	\begin{align*}
		\|u\|_{\mathrm{L}^p(I,\mathbb{V})}
		:=\left(\int_I\|u(t)\|_{\mathbb{V}}^p dt\right)^{1/p}<\infty,
	\end{align*}
	where functions equal almost everywhere are identified. 
		These spaces are Banach spaces with their respective norms.
		
		\subsection{Resolvent Family}\label{sec-resolvent}
		Here, we introduce the resolvent family associated with the following abstract Cauchy problem for neutral integro-differential equations:
		\begin{equation}\label{NIDE}
			\begin{cases}
				\begin{aligned}
					\frac{d}{dt}\left(q(t)+\int_0^t N(t-s)q(s)ds\right)&= \mathscr{A}q(t)+\int_0^t G(t-s)q(s)ds,~t\in I,\\
					q(0)&=q_0.
				\end{aligned}
			\end{cases}
		\end{equation}
		
		\begin{definition}[{\cite[Defintion 2.1]{dos2011existence}}]\label{def-res}
			A one-parameter family $(\mathscr{R}(t))_{t\geq 0} \subset \mathcal{L}(\mathbb{V})$ is a resolvent family for \eqref{NIDE} if the following conditions hold:
			\begin{itemize}
				\item[(a)]  The function $ \mathscr{R}(\cdot) \colon [0, \infty) \to \mathcal{L}(\mathbb{V})$ is strongly continuous, exponentially bounded, with $\mathscr{R}(0)x = x$ for all $x \in \mathbb{V}$.
				\item[(b)] $\mathscr{R}(\cdot)x \in \mathrm{C}([0, \infty); \mathscr{D}_1(\mathscr{A})) \cap \mathrm{C}^1((0, \infty); \mathbb{V}) \text{ for } x \in \mathscr{D}(\mathscr{A})$, and 
				\begin{align}\label{2.3}
					\frac{d}{dt}\left(\mathscr{R}(t)+\int_0^t N(t-s)\mathscr{R}(s)xds\right)&= \mathscr{A}\mathscr{R}(t)x+\int_0^t G(t-s)\mathscr{R}(s)x ds,  \\
					\frac{d}{dt}\left(\mathscr{R}(t)+\int_0^t \mathscr{R}(t-s)N(s)x ds)\right)&= \mathscr{R}(t)\mathscr{A}x +\int_0^t \mathscr{R}(t-s)G(s)x ds, \nonumber
				\end{align}
			\end{itemize}
			for every $t \geq 0$. 
		\end{definition}
		Throughout this paper, we considered the following conditions:
		\begin{itemize}
			\item[\textit{(A1)}] The linear operator $\mathscr{A} \colon \mathscr{D}(\mathscr{A})\subseteq \mathbb{V} \to \mathbb{V}$ generates an analytic semigroup $(\mathcal{S}(t))_{t\geq 0}$. Thus, for constants $M_0 > 0$ and $\vartheta \in (\frac{\pi}{2}, \pi)$, the resolvent set 
			$\rho(\mathscr{A})$ contains $\widetilde{\Lambda}_\vartheta = \{\lambda \in \mathbb{C} \backslash \{0\} \colon |\mathrm{arg}(\lambda)| < \vartheta\}$, and $\|R(\lambda, \mathscr{A})\|_{\mathcal{L}(\mathbb{V})} \leq \frac{M_0}{|\lambda|}$ for all $\lambda \in \widetilde{\Lambda}_\vartheta$.
			\item[\textit{(A2)}] The function $ N : [0,\infty) \to \mathcal{L}(\mathbb{V})$ is strongly continuous and its Laplace transform $\widehat{N}(\lambda)$ extends analytically to $\widetilde{\Lambda}_\vartheta$ and
			$\widehat{N}(\lambda)x$ is absolutely convergent for $x \in \mathbb{V}$ and $Re(\lambda) > 0.$ Furthermore, there exist constants $N_0, N_1, \widetilde{\alpha} > 0$ such that for all $\lambda \in \widetilde{\Lambda}_\vartheta$:
			$$\|\widehat{N}(\lambda)\|_{\mathcal{L}(\mathbb{V})}  \leq \frac{N_0}{|\lambda|^{\widetilde{\alpha}}} \text{  and } \|\widehat{N}(\lambda)x\|_{\mathbb{V}} \leq \frac{N_1}{|\lambda|}\|x\|_1\ \text{ for all }\ x \in \mathscr{D}(\mathscr{A}).$$
			\item[\textit{(A3)}] For each $t \geq 0$, $G(t) \colon \mathscr{D}(G(t)) \subseteq \mathbb{V} \to \mathbb{V}$ is a linear and closed  operator with $\mathscr{D}(\mathscr{A}) \subseteq \mathscr{D}(G(t))$. For every $x \in \mathscr{D}(\mathscr{A})$, $G(\cdot)x$ is strongly measurable on $(0,\infty)$, and there exists $b \in \mathrm{L}^1_{\mathrm{loc}}(\mathbb{R}^+)$ such that
			$$\Vert{}G(t)x\Vert{}_{\mathbb{V}} \le b(t)\Vert{}x\Vert{}_1 \ \text{ for all } \ t > 0.$$
			Furthermore, the Laplace transform $\widehat{G} \colon \widetilde{\Lambda}_{\pi/2} \to \mathcal{L}(\mathscr{D}_1(\mathscr{A}), \mathbb{V})$ extends analytically to $\widetilde{\Lambda}_\vartheta$, with $\|\widehat{G}(\lambda)x\|_{\mathbb{V}} \leq \|\widehat{G}(\lambda)\|_{\mathcal{L}(\mathscr{D}_1(\mathscr{A}),\mathbb{V})}\|x\|_1 $ for all $x \in \mathscr{D}(\mathscr{A})$, and $\|\widehat{G}(\lambda)\|_{\mathcal{L}(\mathscr{D}_1(\mathscr{A}),\mathbb{V})}\to 0$ as $|\lambda| \to \infty$.
			\item[\textit{(A4)}] There exists a dense subspace $\widetilde{D} \subseteq \mathscr{D}(\mathscr{A})$ in $\mathscr{D}_1(\mathscr{A})$ and $C_i>0,~ i = 1, 2$ such that $\mathscr{A}(\widetilde{D}) ,\widehat{G}(\lambda)(\widetilde{D}),\widehat{N}(\widetilde{D}) \subseteq \mathscr{D}(\mathscr{A})$ with $\|\mathscr{A}\widehat{G}(\lambda)x\|_{\mathbb{V}} \leq C_1\|x\|_1$ and $\|\widehat{N}(\lambda)x\|_1 \leq \frac{C_2}{|\lambda|^{\widetilde{\alpha}}}\|x\|_1$ for every $x\in \widetilde{D}$ and all $\lambda \in \widetilde{\Lambda}_\vartheta$.
		\end{itemize}
		In the sequel, for $ r>0 \text{ and } \theta \in (\pi/2,\nu)$, we define
		$$ \widetilde{\Lambda}_{r,\theta} = \{\lambda \in \mathbb{C}\setminus\{0\} : |\lambda| > r,\ |\mathrm{arg}(\lambda)| < \theta\},$$
		and 
		$$\widetilde{\Gamma}_{r,\theta} = \bigcup_{k=1}^{3} \widetilde{\Gamma}^{k}_{r,\theta}, \text{ oriented counterclockwise},$$
		where
		\[
		\widetilde{\Gamma}^{1}_{r,\theta} = \{ t e^{i\theta} : t \ge r\},~
		\widetilde{\Gamma}^{2}_{r,\theta} = \{ r e^{i\phi} : -\theta \le \phi \le \theta\} \text{ and } \widetilde{\Gamma}^{3}_{r,\theta} = \{ t e^{-i\theta} : t \ge r\}.
		\]
		\begin{theorem}[{\cite[Theorem 2.1]{dos2011existence}\label{Thm-2.2}}]
			Assume conditions (A1)–(A4) are fulfilled. For problem \eqref{NIDE}, there exists a unique resolvent family $\mathscr{R}(\cdot)$ given by 
			\begin{equation}\label{2.31}
				\mathscr{R}(t)= 
				\begin{cases}
					\frac{1}{2 \pi \iota} \int_{\widetilde{\Gamma}_{r,\theta}} e^{\lambda t}\mathcal{F}(\lambda)d \lambda, \ t>0,\\
					\mathtt{I}, ~t=0,
				\end{cases}
			\end{equation}
			where, $\mathtt{I}$ is the identity operator and $\Omega(\mathcal{F})=\{\lambda \in \mathbb{C}: \mathcal{F}(\lambda):=(\lambda \mathtt{I}+\lambda \widehat{N}(\lambda)-\mathscr{A}-\widehat{G}(\lambda))^{-1} \in \mathcal{L}(\mathbb{V}) \}.$ 
		\end{theorem}
		
		\begin{lemma}[{\cite[Lemma 2.9]{arora2025controllability}}]\label{lemma 3.1}
			If  $(\mathscr{R}(t))_{t>0}$ is compact operator, then the following holds:
			\begin{itemize}
				\item[(a)] $\lim\limits_{h \to 0^+} \|\mathscr{R}(t+ h)-\mathscr{R}(h)\mathscr{R}(t)\|_{\mathbb{V}} = 0$ for all $t > 0$.
				\item[(b)] $\lim\limits_{h \to 0^+} \|\mathscr{R}(t)-\mathscr{R}(h)\mathscr{R}(t-h)\|_{\mathbb{V}} = 0$ for all $t > 0$.
			\end{itemize}
		\end{lemma}
		
		\subsection{Multivalued maps}
		In this subsection, we recall some basic concepts of multivalued analysis needed later. Here, $\multimap$ is used to denote a set-valued mapping. Consider $(\Omega,d)$ is a metric space. Let $\mathcal{P}(\Omega)$ represent the collection of all non-empty subsets of $\Omega$, and
		$$\mathcal{P}_b(\Omega)=\{K \subset \Omega~\colon K \text{ is bounded}\}  \ \mbox{and}\  \mathcal{P}_{cb}(\Omega)=\{K \subset \Omega~\colon K \text{ is bounded and closed}\} .$$
		The Hausdorff metric $d_{H} \text{ on } \mathcal{P}_{cb}(\Omega)$ is given by
		$$   d_{H}(X,Y) = \max\left\{\sup\limits_{x \in X} d(x,Y), \sup\limits_{y \in Y} d(y,X)\right\},$$
		where $X$ and $Y$ are nonempty, closed, and bounded subsets of the metric space $(\Omega, d)$
		and $d(x,Y) = \inf\limits_{y \in Y} d(x,y).$ 
		
		\begin{definition}
			Let $F \colon \Omega \to \mathcal{P}(\Omega)$ be a multivalued map. A point $x \in \Omega$ such that $x \in F(x)$ is called a fixed point of $F$, and $\mathrm{Fix}(F)$ denotes the set of all such fixed points.
		\end{definition}
		\begin{definition}
			Given a metric space $(\Omega^1, d^1)$, a multivalued map $F \colon \Omega \to \mathcal{P}_{cb}(\Omega^1)$ is called a $k$-contraction if there exists $k \in [0, 1)$ such that$$d^1_H (Fx,Fy) \leq k d(x,y) \ \text{ for all }\ x,y \in \Omega.$$
		\end{definition}
		\begin{theorem}[{\cite[Theorem 5]{SBNJ}}]\label{Fixed-1}
			Let $(\Omega, d)$ be a complete metric space. If $F \colon \Omega \to \mathcal{P}_{cb}(\Omega)$ is a $k$-contraction, then $F$ has a fixed point.
		\end{theorem} 
		
		\begin{definition}[\cite{granas2003fixed}]
			For a bounded subset $K$ of a metric space $\Omega$, its Hausdorff measure of noncompactness (Hausdorff MNC) is given as
			$$ \chi(K) = \inf\{\varepsilon > 0~\colon K \text{ can be covered by finitely many balls of radius } \varepsilon\}.$$
		\end{definition}
		For the sake of clarity, let $\chi$ and $\beta$ denote the measures of noncompactness in a general space and a function space, respectively.  For a nonempty bounded subset $K \subseteq \mathrm{C}(I, \mathbb{V})$, we define
		$$\beta(K) = \sup\limits_{t \in I} \chi(K(t)),$$
		where $\chi$ is the Hausdorff MNC in $\mathbb{V} \text{ and } K(t) = \{y(t) \colon y \in K\}.$
		
		\begin{definition}[{\cite[Definition 2.2.6]{kamenskii}}]
			Let $\mathfrak{C}$ be a closed subset of a Banach space $\mathbb{W}$. A map $F \colon \mathfrak{C} \to \mathcal{P}(\mathbb{W})$ is $\chi$-condensing if for any non-relatively compact set $K \subset \mathfrak{C},$
			$$\chi(F(K)) \ngeq \chi(K).$$
		\end{definition}
		Next, we state a key lemma on the measure of noncompactness in $\mathbb{V}$-valued function spaces.
		
		\begin{lemma}[{\cite[Lemma 2.2]{MFPH}}]\label{lemma-2}
			Let $\mathrm{G}\colon I \to  \mathcal{L}(\mathbb{V})$ be a strongly continuous and $K \subset \mathbb{V}$ be a bounded set. Then 
			$$ \beta(\{\mathrm{G}(\cdot)x \colon x \in K\}) \leq \sup_{ t \in I }\norm{\mathrm{G}(t)}_{\mathcal{L}(\mathbb{V})}\chi(K).$$  
		\end{lemma} 
		We recall some basic results on multivalued maps in a Banach space $\mathbb{V}$ equipped with the norm $\norm{\cdot}_{\mathbb{V}}$ (see \cite{kamenskii} for details).
		\begin{definition}\label{def-2.10}
			A multivalued map $F: \mathbb{V} \to \mathcal{P}(\mathbb{W})$ is said to be
			\begin{itemize}
				\item  convex and closed valued, if $F(x)$ is convex and $F(x)$ is closed for all $x\in \mathbb{V}$, respectively.
				\item  bounded, 
				if for any bounded set $B \subset \mathbb{V}$, the image  $F(B)$ is bounded in $\mathbb{W}$.
				
				\item  upper semicontinuous (abbreviated as u.s.c.), if the preimage  $F^{-1}(V)$ defined as $F^{-1}(V )= \{w \in \mathbb{V}~\colon F(w) \subseteq V\}$ is open in $\mathbb{V}$ for every open set $V \subseteq \mathbb{W}.$
			\end{itemize}
		\end{definition}
		From now, the notations $\mathcal{C}(\mathbb{W}) \text{ and } \mathcal{C}v(\mathbb{W})$ denotes the following sets:
		\begin{align*}
			\mathcal{C}(\mathbb{W})&= \{K \in \mathcal{P}(\mathbb{W})~\colon K\text{ is closed}\},\\ 
			\mathcal{C}v(\mathbb{W}) &= \{K \in \mathcal{P}(\mathbb{W})~\colon K\text{ is closed and convex}\},\\ 
			\mathcal{K}v(\mathbb{W}) &= \{K \in \mathcal{P}(\mathbb{W})~\colon K \text{ is compact and convex}\}.
		\end{align*}
		The following result plays a crucial role in the analysis of our work, providing the fundamental tool to establish a fixed point result.
		
		\begin{theorem}[\cite{ILG}]\label{Fixed-4}
			Consider a Hausdorff locally convex topological vector space $\mathbb{V}$ with a compact convex subset $ \mathcal{M}$, and  $\Sigma: \mathcal{M} \multimap  \mathcal{M} $ is a u.s.c. multimap with convex and closed values. Then  $\Sigma$ admits a fixed point in $ \mathcal{M}$.
		\end{theorem}
		
		
	Let us introduce the concept of uniformly integrability and integrably bounded families, which plays a crucial role in our analysis.
		\begin{definition}
			A subset $\Xi\subset \mathrm{L}^1(I, \mathbb{V})$ is uniformly integrable, if for every $\varepsilon>0$, there exists $\delta>0$ such that 
			$\displaystyle \int_{E}\norm{\tilde{g}(s)}_{\mathbb{V}}ds\leq \varepsilon \text{ for all measurable } E \subset I \text{ with } \mathrm{m}(E) \le \delta \text{ and all } \tilde{g} \in \Xi,$
			where $\mathrm{m}(E)$ is the Lebesgue measure.
		\end{definition}
		\begin{definition}
			A multimap $\Phi: I \to \mathcal{P}(\mathbb{V})$ is integrably bounded, if there exists $\eta \in \mathrm{L}^1(I, [0, +\infty))$ such that
			$$\norm{\Phi(t)}_{\mathbb{V}}:= \sup\limits_{z \in \Phi(t)}\norm{z}_{\mathbb{V}}\leq \eta(t) \ \text{ for a.e. } \ t \in I.  $$  
		\end{definition}
		\begin{remark}
			If $\{f_n\}_{n\in \mathbb{N}}\subset \mathrm{L}^1(I,\mathbb{V})$ is integrably bounded, then it is uniformly integrable.
		\end{remark}
		
		\begin{definition}[{\cite[Definition 4.2.1]{kamenskii}}]
			The sequence $\{f_n\}_{n\in \mathbb{N}} \in  \mathrm{L}^1(I,\mathbb{V}) $ is semicompact if it is integrably bounded and $\{f_n(t)\}_{n\in \mathbb{N}}$ is relatively compact for \mbox{a.e.} $t\in I$.
		\end{definition}
		Now, let us define the Clarke generalized directional derivative (see \cite[page 27]{clarke1990optimization}) of a locally Lipschitz functional $\mathbb{F} \colon \mathbb{V}\to \mathbb{R}$ at $y$ in the direction $\tilde{x}$, denoted by $\mathbb{F}^0(y; \tilde{x})$ and defined as follows:
		\begin{equation*}
			{\mathbb{F}}^0(y;\tilde{x})=\lim\limits_{\epsilon\to 0^+}\sup_{\xi\to y}\frac{\mathbb{F}(\xi+\epsilon \tilde{x})-\mathbb{F}(\xi)}{\epsilon}.
		\end{equation*} 
		Additionally, the Clarke subdifferential of $\mathbb{F}$ at $y$, denoted by $\partial \mathbb{F}(y) \subset \mathbb{V}^*$, is given by
		\begin{equation}\label{subdiff-1}
			\partial \mathbb{F}(y)=\{y^*\in \mathbb{V}^*: {\mathbb{F}}^0(y;\tilde{x})\ge \langle y^*,\tilde{x}\rangle_{\mathbb{V}^*,\mathbb{V}} \ \text{ for all }\ \tilde{x}\in \mathbb{V}\}.
		\end{equation}

		\subsection{Differential inclusion problem}
		This subsection shows  that the hemivariational control problem \eqref{HV1} is connected to the following differential inclusion problem 
		:
		\begin{equation}\label{PC1}
			\begin{cases}
			\begin{aligned}
				\frac{d}{dt}\left(q(t)+\int_0^t N(t-s)q(s)ds\right)&\in \mathscr{A}q(t)+\int_0^t G(t-s)q(s)ds+H\partial \mathbb{F}(t, q(t))+Bu(t)\ \text{ a.e.} \ t\in I,\\
				q(0)&\in q_0+h(q).
		\end{aligned}
		\end{cases}
	\end{equation}
		For this purpose,  we define the selection map $S_{\partial \mathbb{F}}: \mathrm{C}(I, \mathbb{V})\multimap \mathrm{L}^1(I,\mathbb{V}^*)$ as follows:
		\begin{equation*}
			S_{\partial \mathbb{F}}(q)=\{f\in \mathrm{L}^1(I,\mathbb{V}^*): f(t)\in \partial \mathbb{F}(t, q(t)) \ \text{ a.e. }\ t\in I\}.
		\end{equation*}
		By preassuming that the multimap $S_{\partial \mathbb{F}}$ is well-defined, we introduce the notion of the solution for the problem \eqref{PC1}.
		\begin{definition}\label{solution of DI}
			Let $u(\cdot)\in \mathrm{L}^2(I,\mathbb{U})$. A function $q(\cdot)$ is a classical (respectively mild) solution of the problem \eqref{PC1} if $q(\cdot)$ is a classical (respectively mild) solution of the problem
			\begin{equation*}
				\begin{cases}
				\begin{aligned}
\frac{d}{dt}\left(q(t)+\int_0^t N(t-s)q(s)ds\right)&=\mathscr{A}q(t)+\int_0^t G(t-s)q(s)ds+Hf(t)+Bu(t), ~t\in I,\\
					q(0)&=q_0+z,
			\end{aligned}
		\end{cases}
			\end{equation*}
			for some $f\in S_{\partial \mathbb{F}}(q) \text{ and } z \in h(q)$.
		\end{definition}
		Suppose $q$ is a solution of the inclusion problem \eqref{PC1}. By means of Definition \ref{solution of DI}, it follows that 
		\begin{equation}\label{3.4}
			\frac{d}{dt}\left(q(t)+\int_0^t N(t-s)q(s)ds\right)=\mathscr{A}q(t)+\int_0^t G(t-s)q(s)ds+Hf(t)+Bu(t), ~t\in I,
		\end{equation}
		for some $f\in S_{\partial F}(q)$. Taking the duality pairing on both sides in \eqref{3.4}, we obtain for all $v^*\in \mathbb{V}^*$ that 
		\begin{equation}\label{Class}
			\left\langle  \frac{d}{dt}\left(q(t)+\int_0^t N(t-s)q(s)ds\right),v^*\right\rangle=\langle \mathscr{A}q(t)+\int_0^t G(t-s)q(s)ds+Bu(t)+Hf(t), v^*\rangle, 
		\end{equation}
		for a.e. $t\in I$. Since $f(t) \in \partial \mathbb{F}(t, q(t))$ and using the definition of the Clarke subdifferential, we get
		\begin{equation*}
			\langle f(t), v\rangle_{\mathbb{V}^*,\mathbb{V}}\le \mathbb{F}^0(t, q(t);v)\ \text{ for all } \ v\in \mathbb{V},~t\in I. 
		\end{equation*}
		Choosing $v = H^*v^*$ for $v^* \in \mathbb{V}^*$ and using the definition of the adjoint operator yields
		\begin{equation}\label{CHV4}
			\langle Hf(t), v^* \rangle_{\mathbb{V},\mathbb{V}^*}\le \mathbb{F}^0(t, q(t); H^*v^*),~t\in I.
		\end{equation}
		Combining \eqref{Class} and \eqref{CHV4}, we deduce that $q$ satisfies
		\begin{align*}
			&\left\langle -\frac{d}{dt}\left(q(t)+\int_0^t N(t-s)q(s)ds\right)+\mathscr{A}q(t)+\int_0^t G(t-s)q(s)ds+Bu(t), v^*\right\rangle_{\mathbb{V},\mathbb{V}^*} \nonumber\\&\qquad+\mathbb{F}^0(t, q(t); H^*v^*)\ge 0,
		\end{align*}
		for all $ v^*\in \mathbb{V}^*.$ Hence, we prove that the solution of the inclusion problem \eqref{PC1} is also the solution of the hemivariational control problem \eqref{HV1}.
		
		\section{Existence of Mild Solution}\label{sec3}\setcounter{equation}{0}
		In this section, we study the existence of mild solution for the NIDI given in problem \eqref{PC1} under various assumptions of non-local inclusion conditions $h.$ We impose the hypotheses on the function $\mathbb{F}: I\times \mathbb{V} \to \mathbb{R}$ as stated below:
		\begin{itemize}
			\item[\textit{(F1)}] For each $x \in \mathbb{V}$, $\mathbb{F}(\cdot, x)$ is measurable, and for a.e. $t \in I$, $\mathbb{F}(t, \cdot)$ is locally Lipschitz continuous.
			
			\item[\textit{(F2)}] There exists $\gamma\in \mathrm{L}^1(I, \mathbb{R}^+)$ such that
			\begin{equation*}
				\norm{\partial \mathbb{F}(t,x)}_{\mathbb{V}^*}=\sup\{\norm{y}_{\mathbb{V}^*}: y\in \partial \mathbb{F}(t,x)\}\le \gamma(t)\ \mbox{ a.e.} \ t\in I, x\in \mathbb{V}.
			\end{equation*}
			\item[\textit{(F2)$^\prime$}] There exists $\gamma_{R}\in \mathrm{L}^1(I, \mathbb{R}^+)$ such that
			\begin{equation*}
				\norm{\partial \mathbb{F}(t,x)}_{\mathbb{V}^*}=\sup\{\norm{y}_{\mathbb{V}^*}: y\in \partial \mathbb{F}(t,x)\}\le \gamma_{R}(t)\ \mbox{ a.e.}\ t\in I, x\in \mathbb{V} \text{ with } \|x\|_{\mathbb{V}} \leq R, 
			\end{equation*}
			where 
			\begin{equation*}
				\liminf_{R \to\infty}\frac{1}{R}\int_{0}^{a}\gamma_{R}(s)ds=0.
			\end{equation*}
		\end{itemize}
		Under hypotheses \textit{(F1)}–\textit{(F2)}, and by \cite[Proposition 2.1.2]{clarke1990optimization} together with \cite[Proposition 3.44]{SMAOMS}, the multimap $\mathbb{G} \colon I \times \mathbb{V} \to \mathcal{P}(\mathbb{V}^*)$ defined by $\mathbb{G}(t,x) = \partial \mathbb{F}(t,x)$ satisfies the following properties:
		\begin{itemize}
			\item[\textit{(G1)}] The multimap $\mathbb{G}$ has nonempty, convex and weakly compact values. Also, the multimap $\mathbb{G}(\cdot,x): I\to \mathcal{P}(\mathbb{V}^*)$ is measurable for all $x\in \mathbb{V}$.
			\item[\textit{(G2)}] The multimap $\mathbb{G}(t,\cdot): \mathbb{V}\to \mathcal{P}(\mathbb{V}^*)$ has strongly weakly closed graph, that is, if $x_n\to x$ in $\mathbb{V}$, $y_n \xrightarrow{w} y$ in $\mathbb{V}^*$ with $y_n\in \mathbb{G}(t,x_n)$, then $y\in \mathbb{G}(t,x)$.
			\item[\textit{(G3)}] There exists $\gamma\in \mathrm{L}^1(I,\mathbb{R}^+)$ such that
			\begin{equation*}
				\norm{\mathbb{G}(t,x)}_{\mathbb{V}^*}=\sup\{\norm{y}_{\mathbb{V}^*}: y\in \mathbb{G}(t,x)\}\le \gamma(t) \ \mbox{ a.e.} \ t\in I, x\in \mathbb{V}.
			\end{equation*} 
		\end{itemize}
		Applying \cite[Lemma 5.3]{SMAOMS}, we conclude that $S_\mathbb{G}(q)$ is nonempty for all $q \in \mathrm{C}(I, \mathbb{V})$.
		Let us now introduce the notion of a mild solution for problem \eqref{PC1} as follows:
		\begin{definition}
			Fix $u(\cdot)\in \mathrm{L}^2(I,\mathbb{U})$. A function $q\in \mathrm{C}(I, \mathbb{V}) $ is a mild solution of \eqref{PC1} if it satisfies
			\begin{equation*}
				q(t)=\mathscr{R}(t)(q_0+z)+\int_{0}^{t}\mathscr{R}(t-s)Hf(s)ds+\int_{0}^{t}\mathscr{R}(t-s)Bu(s)ds, ~t\in I,
			\end{equation*}
			for some $f\in S_{\partial \mathbb{F}}(q) \text{ and } z \in h(q)$.
		\end{definition}
		We now investigate the existence of mild solutions to problem \eqref{PC1} under 
		both compact and non-compact cases of $\mathscr{R}(t)$ together with the additional non-local inclusion condition which is written as follows: 
		\begin{itemize}
			\item[\textit{(H1)}] The function $h : \mathrm{C}(I, \mathbb{V}) \to \mathcal{P}(\mathbb{V})$ maps bounded sets into bounded sets. In particluar, for each $\mathfrak{R}>0,$ there exists $\Phi_{\mathfrak{R}} > 0$ such that
			$$\sup\{\norm{z}_{\mathbb{V}} : z \in h(q)\} \leq \Phi_{\mathfrak{R}} \text{ for all } q \in \mathrm{C}(I, \mathbb{V}) \text{ with } \norm{q}_{\infty} \leq \mathfrak{R}.$$ 
			\item[\textit{(H2)}] The multimap  $h : \mathrm{C}(I, \mathbb{V}) \to \mathcal{C}v(\mathbb{V})$ is compact.
			\item[\textit{(H3)}] The multimap $h$ has a closed graph, that is, if
			$z_n \to z $ in $\mathbb{V}$ and $q_n \to q$ in $\mathrm{C}(I, \mathbb{V})$ with $z_n \in h(q_n)$, then $z \in h(q).$
		\end{itemize}
		\subsection{The operator $\mathscr{R}(\cdot)$ is compact}
		In this subsection, we proceed to investigate the existence of mild solutions to problem \eqref{PC1} by assuming the compactness of the family $(\mathscr{R}(t))_{t > 0}$. Throughout this work, we further take 
		$$\sup\limits_{t \in I} \|\mathscr{R}(t)\|_{\mathcal{L}(\mathbb{V})} \leq M.$$
		\begin{theorem}\label{Theorem-2}
			Fix $u(\cdot)\in \mathrm{L}^2(I, \mathbb{U})$, and consider the multivalued map $h : \mathrm{C}(I, \mathbb{V}) \to \mathcal{P}(\mathbb{V})$ satisfying conditions (H1)-(H3). Additionally, suppose that conditions (F1)-(F2) are satisfied and there there exists $\mathfrak{R}> 0$ such that 
			\begin{align}\label{assum-1}
				M \bigg( \norm{q_0}_{\mathbb{V}}+ \Phi_{\mathfrak{R}} + \norm{H}_{\mathcal{L}(\mathbb{V}^*,\mathbb{V})} \int_{0}^{a} \gamma(s)ds+ \norm{B}_{\mathcal{L}(\mathbb{U},\mathbb{V})}a^{1/2}\norm{u}_{\mathrm{L}^2(I, \mathbb{U})}\bigg) \leq \mathfrak{R}.
			\end{align}
			Then the problem \eqref{PC1} admits a mild solution.
		\end{theorem}
		\begin{proof}
			Let us define the sets 
			\begin{equation*}
				B_{\mathfrak{R}}=\{q\in \mathrm{C}(I, \mathbb{V}): \norm{q}_{\infty} \le \mathfrak{R}\}, ~\mathfrak{R}>0,
			\end{equation*}
			and
			\begin{equation*}
				Q_{\mathbb{F}}=\{f\in \mathrm{L}^1(I, \mathbb{V}^*): \norm{f(t)}_{\mathbb{V}^*}\le \gamma(t)~\text{a.e.}~t\in I\},
			\end{equation*}
			where $\gamma$ is as in Hypothesis \textit{(F2)}. For fixed $z \in h(B_{\mathfrak{R}})$, we define the map $$\Gamma_z:Q_{\mathbb{F}} \to \mathrm{C}(I, \mathbb{V})$$ as: 
			\begin{equation}\label{equation-1}
				\Gamma_z(f)(t)=\mathscr{R}(t)(q_0+z)+\int_{0}^{t}\mathscr{R}(t-s)Hf(s)ds+\int_{0}^{t}\mathscr{R}(t-s)Bu(s)ds \ \text{ for all } \ t \in I. 
			\end{equation}   
			The proof is divided into several steps written as follows:
			\vskip 0.2cm
			\noindent		
			\textbf{STEP-I:} \emph{The map $\Gamma_z$ maps the set $Q_{\mathbb{F}}$ into $B_{\mathfrak{R}}$ for some $\mathfrak{R}>0$.} Let
			$y = \Gamma_z(f) \text{ for } f \in Q_{\mathbb{F}}$. Then, from \eqref{equation-1} and H{\"o}lder’s inequality, we compute
			\begin{align*}
				\norm{y(t)}_{\mathbb{V}}&\leq M\norm{q_0+z}_{\mathbb{V}}+M\int_{0}^{t}\norm{H}_{\mathcal{L}(\mathbb{V}^*,\mathbb{V})}\norm{f(s)}_{\mathbb{V}^*}ds+\norm{\int_{0}^{t}\mathscr{R}(t-s)Bu(s)ds}_{\mathbb{V}}  \\
				& \leq  M  (\norm{q_0}_{\mathbb{V}} +\|z\|_{\mathbb{V}})+  M \norm{H}_{\mathcal{L}(\mathbb{V}^*,\mathbb{V})} \int_{0}^{a} \gamma(s)ds+  M \norm{B}_{\mathcal{L}(\mathbb{U},\mathbb{V})}a^{1/2}\norm{u}_{\mathrm{L}^2(I, \mathbb{U})}  \\
				& \leq  M  \left(\norm{q_0}_{\mathbb{V}} +\Phi_{\mathfrak{R}} + \norm{H}_{\mathcal{L}(\mathbb{V}^*,\mathbb{V})} \int_{0}^{a} \gamma(s)ds+ \norm{B}_{\mathcal{L}(\mathbb{U},\mathbb{V})}a^{1/2}\norm{u}_{\mathrm{L}^2(I, \mathbb{U})}\right),
			\end{align*}
			where we used condition \textit{(F2)} and \textit{(H1)}. Thus, by using \eqref{assum-1}, we conclude that $y \in B_{\mathfrak{R}}$ for some $\mathfrak{R}>0 $.
			\vskip 0.2cm
			\noindent
			\textbf{STEP-II:} \emph{The map $\Gamma_z$ is continuous from ${Q_{\mathbb{F}}}_w$ into $\mathrm{C}(I, \mathbb{V})$.} Observe that the set $Q_{\mathbb{F}}$ is a compact convex metrizable set in $\mathrm{L}^1(I, \mathbb{V}^*)_w$. Therefore, it suffices to prove the sequential continuity of the operator $\Gamma_z$. For proving this, consider $\{f_n\}_{n\in \mathbb{N}}\subset Q_{\mathbb{F}}$ such that 
			$$ f_n \xrightarrow{w} f \text{ in } \mathrm{L}^1(I, \mathbb{V}^*).$$
			Let $y_n=\Gamma_z(f_n)$, so that the following equation holds:
			\begin{equation}\label{4.711}
				y_n(t)=\mathscr{R}(t)(q_0+z)+\int_{0}^{t}\mathscr{R}(t-s)Hf_n(s)ds+\int_{0}^{t}\mathscr{R}(t-s)Bu(s)ds \ \text{ for all } \ t \in I.
			\end{equation}
			Furthermore, consider $\varphi \colon \mathrm{L}^1(I,\mathbb{V}^*) \to \mathrm{L}^1(I, \mathbb{V})$ defined by $$  \varphi(f)(t)= H f(t), $$
			where $H \colon \mathbb{V}^* \to \mathbb{V}$ is a bounded linear operator.
			For $\zeta \in (\mathrm{L}^1(I,\mathbb{V}))^*$, define $\Psi: \mathrm{L}^1(I,\mathbb{V}^*)\to \mathbb{R}$ as $\Psi=\zeta \circ \varphi$. Since $f_n \xrightarrow{w} f \text{ in } \mathrm{L}^1(I,\mathbb{V}^*)$, 
			weak convergence implies that
			\begin{equation*}
				\zeta(\varphi(f_n))=\Psi(f_n)\to \Psi(f)=\zeta(\varphi(f)) \ \text{ as } \ n\to\infty.
			\end{equation*}
			Subsequently
			\begin{align*}
				\varphi(f_n) \xrightarrow{w} \varphi(f) \ \text{ in } \ \mathrm{L}^1(I,\mathbb{V}).
			\end{align*}
			By the compactness of $(\mathscr{R}(t))_{t > 0}$  and \cite[Lemma 2.11]{arora2025controllability}, we obtain
			$$   \int_{0}^{t}\mathscr{R}(t-s)\varphi(f_n)(s)ds \to \int_{0}^{t}\mathscr{R}(t-s)\varphi(f)(s)ds, $$
			and consequently, we get
			\begin{align}\label{compact-con}
				\int_{0}^{t}\mathscr{R}(t-s)H f_n(s)ds \to \int_{0}^{t}\mathscr{R}(t-s)H f(s)ds \ \text{ in } \ \mathrm{C}(I, \mathbb{V}).
			\end{align}
			Therefore, passing limit in \eqref{4.711} as $n \to \infty$, we obtain
			$$ y_n(t) \to y(t)=\mathscr{R}(t)(q_0+z)+\int_{0}^{t}\mathscr{R}(t-s)Hf(s)ds+\int_{0}^{t}\mathscr{R}(t-s)Bu(s)ds \ \text{ for all } \ t \in I.$$
			Hence, we conclude that the map $\Gamma_z$ is continuous from $Q_{\mathbb{F}}$ with weak topology into $\mathrm{C}(I, \mathbb{V})$. Since, $Q_{\mathbb{F}}$ is weakly compact in $\mathrm{L}^1(I, \mathbb{V}^*)$, it follows that $\mathpzc{K}= \Gamma_z(Q_{\mathbb{F}})$ is compact in $\mathrm{C}(I, \mathbb{V})$.
			\vskip 0.2cm
			\noindent
			\textbf{STEP-III:} \emph{Construction of multimap $\Upsilon: Q_{\mathbb{F}}\multimap \mathrm{L}^1(I,\mathbb{V}^*)$.} Let us define  $\Upsilon: Q_{\mathbb{F}} \multimap \mathrm{L}^1(I,\mathbb{V}^*)$ as follows:
			\begin{equation*}
				\Upsilon(f)=S_{\partial \mathbb{F}}(\Gamma_z(f)),~ f\in Q_{\mathbb{F}}.
			\end{equation*}
			In the above,  $S_{\partial \mathbb{F}}: \Gamma_z(Q_{\mathbb{F}})\multimap \mathrm{L}^1(I,\mathbb{V}^*)$ is given by
			\begin{equation}\label{def-S}
				S_{\partial \mathbb{F}}(q)=\{f \in \mathrm{L}^1(I,\mathbb{V}^*): f(t)\in \partial \mathbb{F}(t, q(t)) \ \text{a.e.} \ t\in I\}, \ q\in \Gamma_z(Q_{\mathbb{F}}).
			\end{equation}
			Clearly, $\Upsilon$ is well-defined and using the characteristics of $S_{\partial \mathbb{F}}$ (see \cite[Lemma 5.3]{SMAOMS}), $\Upsilon(f)$ is convex, weakly compact subset of $\mathrm{L}^1(I,\mathbb{V}^*)$. Furthermore, Hypothesis \textit{(G3)} implies $\Upsilon(Q_{\mathbb{F}}) \subset Q_{\mathbb{F}}$.
			\vskip 0.2cm
			\noindent
			\textbf{STEP-IV:} \emph{The multimap $\Upsilon$ is u.s.c. in $\mathrm{L}^1(I,\mathbb{V}^*)$ endowed with weak topology.} By virtue of Hypothesis \textit{(G3)}, the multimap $\Upsilon$ maps bounded subsets of $Q_{\mathbb{F}}$ into weakly compact subsets of $Q_{\mathbb{F}}$, and is therefore weakly compact. Consequently, by \cite[Theorem 1.1.12]{kamenskii}, it remains only to prove that $\Upsilon$  has a weakly closed graph. To prove this, let $\{g_n\}_{n\in \mathbb{N}},\{f_n\}_{n\in \mathbb{N}}\subset Q_{\mathbb{F}}$ be such that 
			$$f_n \xrightarrow{w} f, ~ g_n \xrightarrow{w} g \ \text{ in } \ \mathrm{L}^1(I,\mathbb{V}^*)$$
			with $g_n\in \Upsilon(f_n)$ for all $n\in \mathbb{N}$. We need to demonstrate that $g \in \Upsilon(f)$. By using the fact that $g_n\in \Upsilon(f_n),$ we obtain
			\begin{equation*}
				g_n\in \Upsilon(f_n)=S_{\partial \mathbb{F}}(\Gamma_z(f_n)),~ f_n\in Q_{\mathbb{F}}.
			\end{equation*}
			Since $f_n \xrightarrow{w} f$ in $\mathrm{L}^1(I,\mathbb{V}^*)$ and the map $\Gamma_z$ is continuous from ${Q_{\mathbb{F}}}_w$ into $\mathrm{C}(I, \mathbb{V})$ (see STEP-II), we obtain
			\begin{equation*}
				\Gamma_z(f_n)\to \Gamma_z(f)  \ \text{ in } \ \mathrm{C}(I, \mathbb{V}).
			\end{equation*}
			Furthermore, by \cite[Lemma 5.3]{SMAOMS}, the multifunction $S_{\partial \mathbb{F}}$ has a strongly-weakly closed graph. Therefore, passing limit $n \to \infty$ in the relation
			$$ g_n\in S_{\partial \mathbb{F}}(\Gamma_z(f_n)),$$
			we conclude that
			$g \in S_{\partial \mathbb{F}}(\Gamma_z(f))=\Upsilon(f)$. 
			\vskip 0.2cm
			\noindent
			\textbf{STEP-V:} \emph{Existence of a fixed Point of the multimap $\Upsilon$.} By arguments in Steps III and IV, we verify that the map $\Upsilon: {Q_{\mathbb{F}}}_w \multimap {Q_{\mathbb{F}}}_w$ fulfills the hypotheses of Theorem \ref{Fixed-4}. Hence, by Theorem \ref{Fixed-4}, $\Upsilon$ admits a fixed point. Therefore, we have 
			\begin{equation*}
				f \in S_{\partial \mathbb{F}}(\Gamma_z(f)), ~ f\in Q_{\mathbb{F}}.
			\end{equation*}
			Let $q=\Gamma_z(f)$, then $f \in S_{\partial \mathbb{F}}(q)$. The definition of the map $\Gamma_z$ gives 
			\begin{align}\label{mild solution-1}
				q(t)=\mathscr{R}(t)(q_0+z)+\int_{0}^{t}\mathscr{R}(t-s)Hf(s)ds+\int_{0}^{t}\mathscr{R}(t-s)Bu(s)ds \ \text{ for all } \ t \in I.
			\end{align}
			\vskip 0.2cm
			\noindent
			\textbf{STEP-VI:} \emph{Initail data $z \in h(q)$.} To conclude that $q$ given by \eqref{mild solution-1} is a mild solution of problem \eqref{PC1}, it remains to verify that $z \in h(q)$. For this, we can define a map
			$\mathcal{E}: h(\mathpzc{K})\rightarrow \mathrm{C}(I, \mathbb{V})$ as follows:
			\begin{equation*}
				\mathcal{E}(z)=q,
			\end{equation*}
			where $q \in \mathrm{C}(I, \mathbb{V})$ is
			defined as below:
			$$ q(t)=\mathscr{R}(t)(q_0+z)+\int_{0}^{t}\mathscr{R}(t-s)Hf(s)ds+\int_{0}^{t}\mathscr{R}(t-s)Bu(s)ds \ \text{ for all } \ t \in I.$$ Clearly, the map $\mathcal{E}$ is well defined and $\|\mathcal{E}(h(\mathpzc{K}))\|_{\infty} \leq \mathfrak{R}$. Using this mapping, we define the multivalued map $\aleph: h(\mathpzc{K}) \multimap \mathbb{V}$ as follows:
			\begin{equation*}
				\aleph(z)=h(q)
			\end{equation*}
			where $q=\mathcal{E}(z)$ for each $ z\in h(\mathpzc{K})$. It follows immediately that the multimap $\aleph$ is well defined. Moreover, since the values of $h$ are closed and convex, the multimap $\aleph$ also has closed and convex values. Also, observe that, by means of Hypotheses \textit{(H2)}  and \textit{(H3)}, the set $h(\mathpzc{K})$ is compact
			in $\mathbb{V}$.
			In order to apply Theorem \ref{Fixed-4}, it remains to prove the upper semicontinuity of the multimap $$\aleph:h(\mathpzc{K})  \multimap h(\mathpzc{K}).$$ 
			Firstly, we verify that \(\aleph\) maps \(h(\mathpzc{K})\) into itself. Let
			$$
			y\in\aleph(z)=h(q), \quad z\in h(\mathpzc{K}).$$
			It is enough to show that $ q\in\mathpzc{K}$. Since $ f\in S_{\partial \mathbb{F}}(q)$, Hypothesis \textit{(F2)} yields implies that $  f\in Q_{\mathbb{F}}.$ Thus, we have
			$$ q=\mathcal{E}(z)=\Gamma_z(f)\in\Gamma_z(Q_{\mathbb{F}})=\mathpzc{K}.$$
			By means of
			\cite[Theorem 1.1.12]{kamenskii}, it suffices to prove the multimap 
			$\aleph$ has a closed graph. For this, consider
			$y_n\in \aleph(z_n)$ with 
			$$ y_n \to y,~ z_n\to z\ \text{ in } \ \mathbb{V}.$$ 
			We have to prove $y\in \aleph(z)$. Since $y_n\in \aleph(z_n)$, it further implies $y_n\in h(q_n)$, where
			\begin{equation*}
				q_n(t)=\mathscr{R}(t)(q_0+z_n)+\int_{0}^{t}\mathscr{R}(t-s)Hf_n(s)ds+\int_{0}^{t}\mathscr{R}(t-s)Bu(s)ds \ \text{ for all } \ t \in I,
			\end{equation*}
			for some $f_n \in S_{\partial \mathbb{F}}(q_n).$ Then, by virtue of assumption \textit{(F2)}, the sequence $\{f_n\}_{n \in \mathbb{N}}$ is uniformly integrable in $\mathrm{L}^1(I,\mathbb{V}^*).$ Thus, by using Dunford-Pettis Theorem \cite[Theorem 4.75]{gasinski2016exercises}, $\{f_n\}_{n \in \mathbb{N}}$ is relatively weakly compact in $ \mathrm{L}^1(I,\mathbb{V}^*)$. Therefore, we can find a subsequence of $\{f_n\}_{n \in \mathbb{N}}$ (labeled with the same symbol) such that
			$$ f_n \xrightarrow{w} f \ \text{ in } \ \mathrm{L}^1(I,\mathbb{V}^*).$$
			By following the same argument as in STEP-II, we obtain \eqref{compact-con}. Moreover, since $z_n\to z$ in $\mathbb{V}$ and  the compactness of the family $(\mathscr{R}(t))_{t > 0}$ together with \cite[Lemma 2.11]{arora2025controllability}, yields
			\begin{equation}\label{exp-q}
				q_n(t)\to q(t)=\mathscr{R}(t)(q_0+z)+\int_{0}^{t}\mathscr{R}(t-s)Hf(s)ds+\int_{0}^{t}\mathscr{R}(t-s)Bu(s)ds \ \text{ for all } \ t \in I,
			\end{equation}
			uniformly in $\mathrm{C}(I, \mathbb{V})$. By means of Hypotheses \textit{(H3)}, we obtain $y \in h(q)$, where $q$ is given
			in \eqref{exp-q}. Thus, the definition of the multimap $\aleph$ implies $y \in \aleph(z)$.	In conclusion, $\aleph$ fulfills all hypotheses of Theorem \ref{Fixed-4}. Therefore, we obtain a fixed point of the multimap $\aleph$, say $z$. Thus, the definition of $\aleph$ implies
			$$ z  \in \aleph(z) =h(q),$$
			which completes the proof.
		\end{proof}
		Let us now prove the existence of a mild solution of the problem \eqref{PC1} under the following condition on the Clarke subdifferential of $\mathbb{F}$ and $h : \mathrm{C}(I, \mathbb{V})\to \mathcal{P}(\mathbb{V})$ written as follows:
		\begin{itemize}
			\item[\textit{(F)}] There exists $\eta \in \mathrm{L}^1(I,\mathbb{R}^+)$ such as
			\begin{align*}
				d_{H}(\partial \mathbb{F}(t,z_2(t)),\partial \mathbb{F}(t,z_1(t))) \leq \eta(t) \norm{z_2(t)-z_1(t)}_{\mathbb{V}} \ \text{ for a.e. }\ t \in I.
			\end{align*}
			\item[\textit{(H4)}] The function $h : \mathrm{C}(I, \mathbb{V}) \to \mathcal{P}_{cb}(\mathbb{V}).$
		\end{itemize}	
		We first establish a useful lemma based on the properties of u.s.c. mappings.
		\begin{lemma}\label{lemma-1}
			Under the assumption that $h$ is u.s.c. and fulfills condition \textit{(H4)}, the multimap $\mathcal{G} \colon \mathrm{C}(I, \mathbb{V}) \to \mathcal{P}(\mathrm{C}(I, \mathbb{V}))$ defined by 
			\begin{align}\label{Def-G}
				\mathcal{G}(q) = \{\mathscr{R}(\cdot)z : z \in h(q)\}
			\end{align}
			is closed-valued and u.s.c.
			\begin{proof}	Since $h(q)$ is bounded by condition \textit{(H4)}, it implies that $\mathcal{G}(q)$ is bounded for all $q \in \mathrm{C}(I, \mathbb{V})$. To prove that $\mathcal{G}(q)$ is closed, consider a sequence $(y_n)_{n \in \mathbb{N}} \in \mathcal{G}(q)$ such that 
				$$y_n \to y \in  \mathrm{C}(I, \mathbb{V}) \ \text{ as } \ n \to \infty.$$ 
				Since $y_n(t)=\mathscr{R}(t)z_n \text{ for some } z_n \in h(q)$ and $t\in I$, then
				\begin{equation}\label{Conv-1}
					y_n(0)=  \mathscr{R}(0)z_n=z_n \to y(0) \ \text{ as }\ n \to \infty.
				\end{equation}
				By using closedness of $h(q),$ we obtain $y(0) \in h(q).$ Consider $\Tilde{y}(t)=\mathscr{R}(t)y(0) \text{ for } t \in I.$ Then $\Tilde{y} \in \mathcal{G}(q)$ and
				\begin{align*}
					\norm{y_n-\Tilde{y}}_{\infty} \leq M \norm{z_n-y(0)}_{\mathbb{V}} \to 0 \ \text{ as }\ n \to \infty,
				\end{align*}
				where we have  used \eqref{Conv-1}. Thus, by uniqueness of limits, we obtain $\Tilde{y}=y \in  \mathcal{G}(q),$ which implies that $\mathcal{G}(q)$ is closed. Hence, $\mathcal{G} : \mathrm{C}(I, \mathbb{V}) \to \mathcal{P}_{cb}(\mathrm{C}(I, \mathbb{V})).$ 
				Let us now consider the map $P: \mathbb{V} \to \mathrm{C}(I, \mathbb{V}) $ defined by 
				\begin{align*}
					P(z)= \mathscr{R}(\cdot)z, ~z \in \mathbb{V}.
				\end{align*}
				Since, $\mathscr{R}(t)$ is linear and bounded operator, thus $P$ is linear and bounded. In addition, since $\mathcal{G}(q) = P(h(q))$, we deduce from Definition \ref{def-2.10} that $\mathcal{G}$ is an u.s.c. map.
			\end{proof}   
		\end{lemma}	
		\begin{theorem}\label{Theorem-1}
			Let $u(\cdot)\in \mathrm{L}^2(I, \mathbb{U})$ be fixed. Suppose that the conditions (F2), (F), and \eqref{assum-1} are satisfied. Furthermore, assume that the multivalued map $h : \mathrm{C}(I, \mathbb{V}) \to \mathcal{P}(\mathbb{V})$ satisfying (H4) and there exists a constant $L_h > 0$ such that
			\begin{align}\label{haus-h}
				d_{H}(h(q),h(p)) \leq L_h \norm{q-p}_{\infty} \ \text{ for all }\ q, p \in \mathrm{C}(I, \mathbb{V}).  
			\end{align}
			If $$M\bigg(L_{h}+ \norm{H}_{\mathcal{L}(\mathbb{V}^*,\mathbb{V})} \int_{0}^{a} \eta(s)ds\bigg) <1,$$
			then the problem \eqref{PC1} admits a mild solution.
			\begin{proof}  Let us consider the map $\Lambda:  \mathrm{C}(I, \mathbb{V}) \to \mathcal{P}(\mathrm{C}(I, \mathbb{V})) $ defined as:
				for $q \in \mathrm{C}(I, \mathbb{V}), \text{ if } y \in \Lambda(q),$ then
				\begin{equation}\label{def-y}
					y(t)=\mathscr{R}(t)q_0+\int_{0}^{t}\mathscr{R}(t-s)Hf(s)ds+\int_{0}^{t}\mathscr{R}(t-s)Bu(s)ds, ~t\in I,
				\end{equation}
				where $f \in S_{\partial \mathbb{F}}(q)$. We divide the proof into several steps written as follows:
				\vskip 0.1in 
				\noindent 
				\textbf{STEP-I:-} \emph{The multivalued map $\Lambda$ has compact values.} By making suitable modifications to the proof of \cite[Theorem~3.1]{JDNP}, it follows
				$$
				\Lambda:\mathrm{C}(I,\mathbb{V})\to \mathcal{P}(\mathrm{C}(I,\mathbb{V})) $$
				is completely continuous. To establish that $\Lambda $ has compact values, it is sufficient to show that its graph is closed. To this end, let $\{q_n\}_{n\in\mathbb{N}}\subset B_{\mathfrak{R}}$ and $\{y_n\}_{n\in\mathbb{N}}\subset \mathrm{C}(I,\mathbb{V})$ be sequences such that
				$$ q_n\to q, \quad y_n\to y  \ \text{ with } \ y_n\in\Lambda(q_n).$$
				Since $ y_n\in\Lambda(q_n)$, there exists $f_n\in S_{\partial\mathbb{F}}(q_n)$ such that
				\begin{align}\label{Conver-1}
					y_n(t)=\mathscr{R}(t)q_0+\int_{0}^{t}\mathscr{R}(t-s)Hf_n(s)ds+\int_{0}^{t}\mathscr{R}(t-s)Bu(s)ds, ~t\in I. 
				\end{align}
				By assumption \textit{(F2)}, $\{f_n\}_{n \in \mathbb{N}}$ is uniformly integrable in $\mathrm{L}^1(I,\mathbb{V}^*)$. Thus, by the Dunford–Pettis Theorem \cite[Theorem 4.75]{gasinski2016exercises}, $\{f_n\}_{n \in \mathbb{N}}$ is relatively weakly compact in $\mathrm{L}^1(I,\mathbb{V}^*)$. Hence, passing to a subsequence (still denoted by $\{f_n\}_{n \in \mathbb{N}}$), we have
				$$ f_n \xrightarrow{w} f \ \text{ in } \ \mathrm{L}^1(I,\mathbb{V}^*).$$
				By following the same argument as in STEP-II of Theorem \ref{Theorem-2} and passing limit in \eqref{Conver-1} as $n \to \infty$, we obtain
				$$ y_n(t) \to y(t)=\mathscr{R}(t)q_0+\int_{0}^{t}\mathscr{R}(t-s)Hf(s)ds+\int_{0}^{t}\mathscr{R}(t-s)Bu(s)ds, ~t \in I.$$
				Since $f_n \in S_{\partial \mathbb{F}}(q_n) \text{ for all } n \in \mathbb{N}$ with $f_n \xrightarrow{w} f \text{ in } \mathrm{L}^1(I, \mathbb{V}^*)$ and $q_n \to q \text{ in } \mathrm{C}(I, \mathbb{V}),$ therefore, by \cite[Lemma 5.3]{SMAOMS}, we conclude that $f \in S_{\partial \mathbb{F}}(q).$ Hence, by the definition of the multimap $\Lambda$, we confirm that $y \in \Lambda(q)$. Thus, the graph of $\Lambda $ is closed. Since $\Lambda $ is completely continuous, it follows that $\Lambda$ has compact values.
				\vskip 0.1in 
				\noindent 
				\textbf{STEP-II:-} \emph{Construction of closed and bounded valued operator $\Gamma$.} Let us now define the multivalued operator $\Gamma: \mathrm{C}(I, \mathbb{V}) \to \mathcal{P}_{cb}(\mathrm{C}(I, \mathbb{V}))$ by 
				\begin{equation}\label{Def.-Gamma}
					\Gamma(q)= \mathcal{G}(q)+\Lambda (q),
				\end{equation}
				where $\mathcal{G}(\cdot)$ is given in Lemma \ref{lemma-1}. To show that  the multivalued operator $\Gamma$ has closed values, let $\{x_n\}_{n \in \mathbb{N}} \subset \Gamma(q)$ be a sequence such that $x_n \to x$. Then, each $x_n$ can be written as
				$$x_n = g_n + \upsilon_n,~ g_n \in \mathcal{G}(q)\  \text{ and } \ \upsilon_n \in \Lambda (q).$$
				Since $\Lambda$ is compact-valued, there exists a subsequence (still denoted by the same symbol) such that
				$$ \upsilon_{n} \to \upsilon \  \text{ in } \ \Lambda (q).$$
				Furthermore, along a subsequence, we also have
				$$g_{n}=x_n-\upsilon_n \to  x-\upsilon.$$
				Since, $\mathcal{G}$ has closed values (see Lemma \ref{lemma-1}), it follows that
				$$x-\upsilon \in \mathcal{G}(q).$$
				Hence, 
				$$x = \upsilon + (x - \upsilon) \in \Lambda (q) + \mathcal{G}(q)=\Gamma(q).$$
				Thus, $\Gamma(q)$ is closed for each $q \in \mathrm{C}(I, \mathbb{V})$. 
				\vskip 0.1in 
				\noindent 
				\textbf{STEP-III:-} \emph{The multivalued operator $\Gamma$ is a set contraction.} Now, consider $p, q, \Tilde{p},\Tilde{q} \in \mathrm{C}(I, \mathbb{V})$ such that $\Tilde{p} \in  \Gamma(p) \text{ and } \Tilde{q} \in  \Gamma(q)$ such that 
				\begin{align*}
					\Tilde{p}(t) &=\mathscr{R}(t)z_p+ \mathscr{R}(t)q_0+\int_{0}^{t}\mathscr{R}(t-s)Hf_1(s)ds+\int_{0}^{t}\mathscr{R}(t-s)Bu(s)ds,\\
					\Tilde{q}(t) &=\mathscr{R}(t)z_q+ \mathscr{R}(t)q_0+\int_{0}^{t}\mathscr{R}(t-s)Hf_2(s)ds+\int_{0}^{t}\mathscr{R}(t-s)Bu(s)ds,
				\end{align*}
				for some $z_p \in h(p) \text{ and } z_q \in h(q)$.
				
				For sufficiently small $\epsilon > 0$, choose $f_1 \in \partial \mathbb{F}(t, \Lambda(p))$ such that
				\begin{equation}\label{estimate-20}
					\norm{f_2(s)-f_1(s)}_{\mathbb{V}^*} \leq  d(f_2(s),\partial \mathbb{F}(t, \Lambda(p)) + \epsilon \leq \eta(s)\norm{q-p}_{\infty} + \epsilon,
				\end{equation}
				and for a given $\epsilon > 0$, choose $z_p \in h(p)$ such that
				\begin{equation}\label{estimate-21}
					\norm{z_q-z_p}_{\mathbb{V}} \leq  d(z_q ,h(p)) + \epsilon \leq L_{h}\norm{q-p}_{\infty} + \epsilon.
				\end{equation}
				Hence, we can estimate
				\begin{align}\label{estimate-22}
					\norm{  \Tilde{q}(t)- \Tilde{p}(t)}_{\mathbb{V}} &\leq M\norm{z_q-z_p}_{\mathbb{V}}+\norm{\int_{0}^{t}\mathscr{R}(t-s)H(f_2(s)-f_1(s))ds}_{\mathbb{V}} \nonumber \\
					& \leq M\norm{z_q-z_p}_{\mathbb{V}}+M \norm{H}_{\mathcal{L}(\mathbb{V}^*,\mathbb{V})} \bigg(\int_{0}^{a} \eta(s)ds+ a \epsilon \bigg) \norm{q-p}_{\infty}.
				\end{align}
				Combining \eqref{estimate-21} and \eqref{estimate-22}, we obtain 
				$$ \norm{  \Tilde{q}- \Tilde{p}}_{\infty} \leq M\bigg(L_{h}+ \norm{H}_{\mathcal{L}(\mathbb{V}^*,\mathbb{V})} \int_{0}^{a} \eta(s)ds\bigg) \norm{q-p}_{\infty}+M(1+ \norm{H}_{\mathcal{L}(\mathbb{V}^*,\mathbb{V})}a)\epsilon.$$
				This implies that 
				\begin{align*}
					\sup_{\Tilde{q} \in \Gamma(q)}d(\Tilde{q},\Gamma(p))\leq M\bigg(L_{h}+ \norm{H}_{\mathcal{L}(\mathbb{V}^*,\mathbb{V})} \int_{0}^{a} \eta(s)ds\bigg) \norm{q-p}_{\infty}.
				\end{align*}
				Similarly, we get 
				\begin{align*}
					\sup_{\Tilde{p} \in \Gamma(p)}d(\Tilde{p},\Gamma(q))\leq M\bigg(L_{h}+ \norm{H}_{\mathcal{L}(\mathbb{V}^*,\mathbb{V})} \int_{0}^{a} \eta(s)ds\bigg) \norm{q-p}_{\infty}.
				\end{align*}
				Therefore, one can obtain
				\begin{align*}
					d_H{(\Gamma(q),\Gamma(p))}\leq M\bigg(L_{h}+ \norm{H}_{\mathcal{L}(\mathbb{V}^*,\mathbb{V})} \int_{0}^{a} \eta(s)ds\bigg) \norm{q-p}_{\infty}.
				\end{align*}
				From Theorem \ref{Fixed-1}, we conclude that the multivalued map $\Gamma$ has a fixed point which is a mild solution of problem \eqref{PC1}.
			\end{proof}    
		\end{theorem}
		
		\subsection{The operator $\mathscr{R}(\cdot)$ is not compact}
		In this subsection, we establish the existence of mild solutions to problem \eqref{PC1} using the theory of measures of noncompactness under the assumption that $(\mathscr{R}(t))_{t \geq  0}$ is not compact.  Let $\mathcal{T}  \colon  \mathrm{L}^1(I,\mathbb{V}^*) \to \mathrm{C}(I, \mathbb{V})$ be an operator satisfying:
		\begin{itemize}
			\item[\textit{(S1)}] there exists $\mathscr{M}\geq 0$ such that  
			$$ \|\mathcal{T} f(t) - \mathcal{T} g(t)\|_{\mathbb{V}} \leq \mathscr{M} \left( \int_{0}^{t} \|f(s) - g(s)\|_{\mathbb{V}^*} ds \right) \ \text{ for all }\ f,g \in \mathrm{L}^1(I,\mathbb{V}^*),~ t \in I,$$ 
			\item[\textit{(S2)}] for any compact set $\mathcal{K} \subset \mathbb{V}^* $ and sequence $\{g_n\}_{n\in \mathbb{N}} \subset \mathrm{L}^1(I,\mathbb{V}^*)$ with 
			$\{  g_n(t) \}_{n\in \mathbb{N}} \subset \mathcal{K}$ for a.e. $t \in I$, if $g_n \xrightarrow{w} g$, then
			$\mathcal{T}  g_n \to \mathcal{T}   g $.
		\end{itemize}
		Following the approach of \cite[Theorem 5.1.1]{kamenskii}, we state the following lemma:
		\begin{lemma}
			Let $\mathcal{T}   \colon  \mathrm{L}^1(I,\mathbb{V}^*) \to \mathrm{C}(I, \mathbb{V})$ satisfy conditions (S1)-(S2). For every semicompact sequence $\{g_n\}_{n\in \mathbb{N}} \subset \mathrm{L}^1(I,\mathbb{V}^*)$, the sequence $\{\mathcal{T} g_n\}_{n\in \mathbb{N}}$  is relatively compact in $\mathrm{C}(I, \mathbb{V})$. Furthermore, if $g_n \xrightarrow{w} g$, then $\mathcal{T}   g_n \to \mathcal{T}   g $.
		\end{lemma}
Before proceeding further, we introduce an operator $\Theta: \mathrm{L}^1(I,\mathbb{V}^*)\to \mathrm{C}(I, \mathbb{V})$ such that
		\begin{equation}\label{def-Theta}
			\Theta(g)(t)=\int_{0}^{t}\mathscr{R}(t-s)Hg(s)ds,~~ t\in I.
		\end{equation}
		In order to establish the properties required for our existence result, we next show that the operator $\Theta$ constructed in \eqref{def-Theta} fulfills assumptions \textit{(S1)}–\textit{(S2)}.
		\begin{lemma}
			Let $\mathscr{R} \colon I \to  \mathcal{L}(\mathbb{V})$ be a strongly continuous operator family. Then the operator $\Theta $ defined by \eqref{def-Theta}, satisfies the properties (S1)-(S2).
		\end{lemma}
		\begin{proof}
			Firstly, we estimate
			\begin{align*}
				\|\Theta(f)(t) - \Theta(g)(t)\|_{\mathbb{V}} &\leq  \int_{0}^{t} \|\mathscr{R}(t-s)(Hf(s) - Hg(s))\|_{\mathbb{V}} ds \\
				& \leq M \norm{H}_{\mathcal{L}(\mathbb{V}^*,\mathbb{V})} \int_{0}^{t} \|f(s) - g(s)\|_{\mathbb{V}^*} ds,
			\end{align*}
			which proves the property \textit{(S1)}.  
			
			To verify the property \textit{(S2)}, let $\mathcal{K} \subset \mathbb{V}^*$ be a compact set and consider a sequence
			$\{g_n\}_{n \in \mathbb{N}} \subset \mathrm{L}^1(I,\mathbb{V}^*)$ such that  
			$\{g_n(s)\}_{n \in \mathbb{N}} \subset \mathcal{K}$ for a.e. $s \in I$. Following the argument used in the proof of \cite[Lemma 4.2.1]{kamenskii}, one can prove that the
			sequence $\{\Theta(g_n)(t)\}_{n \in \mathbb{N}} \subset \mathbb{V}$ is relatively compact for every $t \in I$ and the sequence of functions  
			$\{\Theta(g_n)\}_{n \in \mathbb{N}} \subset \mathrm{C}(I, \mathbb{V})$ is equicontinuous. In addition, it is easy to prove that $\{\Theta(g_n)\}_{n \in \mathbb{N}}$ is uniformly bounded. Hence, by the Arzelá-Ascoli theorem, the sequence $\{\Theta(g_n)\}_{n \in \mathbb{N}}$ is relatively compact.
			
			Moreover, since $\Theta \in \mathcal{L}(\mathrm{L}^1(I,\mathbb{V}^*), \mathrm{C}(I, \mathbb{V})),$
			the weak convergence $g_n \xrightarrow{w} g$ implies
			\begin{align*}
				\Theta(g_n) \xrightarrow{w} \Theta(g).
			\end{align*}
			Combining the relative compactness of $\{\Theta(g_n)\}_{n \in \mathbb{N}}$ with the uniqueness of weak limits yields $$\Theta(g_n) \to \Theta(g) \ \text{ in } \ \mathrm{C}(I, \mathbb{V}),$$ which completes the proof.
		\end{proof}
		
		\begin{corollary}\label{corollary 4.7}
			For every semicompact sequence $\{g_n\}_{n \in \mathbb{N}} \subset \mathrm{L}^1(I,\mathbb{V}^*)$, its image $\{\Theta(g_n)\}_{n \in \mathbb{N}} \subset \mathrm{C}(I, \mathbb{V})$ is relatively compact. Moreover, if $g_n \xrightarrow{w} g$, then $\Theta(g_n) \to \Theta(g)$.
		\end{corollary}
		To establish an existence result for functions that do not satisfy condition \textit{(F)}, we instead assume that  $\partial \mathbb{F}$ and the non-local multifunction $h$ satisfy the following condition, respectively.
		\begin{itemize}
			\item[\textit{(F)$^{\prime}$}] There exists $k \in  \mathrm{L}^1(I,\mathbb{R}^+)$ such that
			\begin{align*}
				\chi(\partial \mathbb{F}(t, K)) \leq k (t)\chi(K)\ \text{ a.e. } t \geq 0 \ \text{ for all bounded set }\  K \subseteq \mathbb{V}.
			\end{align*}
			\item[\textit{(H5)}] There exists a constant $k_h \geq 0$ such that
			$$ \chi(h(W)) \leq k_h \beta(W) \ \text{ for all bounded set }\ W \subseteq \mathrm{C}(I, \mathbb{V}).$$
			\item[\textit{(H6)}]  There is $N_h> 0$ such that
			$$ \|h(\mathrm{C}(I, \mathbb{V}))\| = \sup\{\|h(q)\|_{\mathbb{V}}: q \in \mathrm{C}(I, \mathbb{V})\} \leq N_h.$$
		\end{itemize}
		\begin{theorem}\label{Theorem-3}
			Fix $u(\cdot)\in \mathrm{L}^2(I, \mathbb{U})$, and consider the map $h : \mathrm{C}(I, \mathbb{V}) \to \mathcal{K}v(\mathbb{V})$, which is u.s.c. and  fulfills the conditions (H1) and (H5). Moreover, assume that the conditions (F1)-(F2) with (F)$^{\prime}$, \eqref{assum-1} are satisfied and 
			\begin{align}\label{3.38}
				M \bigg(k_h+2 \norm{H}_{\mathcal{L}(\mathbb{V}^*,\mathbb{V})}  \int_{0}^{a} k(s)ds\bigg) <1.
			\end{align}
			Then the problem \eqref{PC1} has a mild solution.
		\end{theorem}
		\begin{proof}
			Define the map $\Lambda \colon \mathrm{C}(I, \mathbb{V}) \to \mathcal{P}(\mathrm{C}(I, \mathbb{V})) $ as follows: for $q \in \mathrm{C}(I, \mathbb{V}),~\Lambda(q)$ consists functions $y \in \mathrm{C}(I, \mathbb{V})$ such that
			\begin{equation*}
				y(t)=\mathscr{R}(t)q_0+\int_{0}^{t}\mathscr{R}(t-s)Hf(s)ds+\int_{0}^{t}\mathscr{R}(t-s)Bu(s)ds, ~t\in I, 
			\end{equation*}
			where $f \in S_{\partial \mathbb{F}}(q)$. Consider two sequences 
			$$y_n \to y,~ q_n \to q \text{ in }  \mathrm{C}(I, \mathbb{V}) \ \text{ with }\  y_n \in \Lambda(q_n).$$
			Since $y_n \in \Lambda(q_n)$, there exist $f_n \in S_{\partial \mathbb{F}}(q_n)$ such that, we can write 
			\begin{align}\label{Converge1}
				y_n(t)=\mathscr{R}(t)q_0+\int_{0}^{t}\mathscr{R}(t-s)Hf_n(s)ds+\int_{0}^{t}\mathscr{R}(t-s)Bu(s)ds, ~t\in I. 
			\end{align}
			Using hypothesis \textit{(F)$^{\prime}$}, we have 
			\begin{align*}
				\chi(\{f_n(t)\}_{n\in \mathbb{N}}) &\leq \chi\big(\partial \mathbb{F}(t,\{q_n(t)\}_{n\in \mathbb{N}})\big) \nonumber\\
				&\leq k(t)\chi(\{q_n(t)\}_{n\in \mathbb{N}})= 0 \ \text{ for a.e. } t \in I.
			\end{align*}
			The above argument ensures that $\{f_n(t)\}_{n\in \mathbb{N}}$ is relatively compact for a.e. $t \in I$. Using  assumption \textit{(F2)}, the sequence $\{f_n\}_{n\in \mathbb{N}}$ is uniformly integrable in $\mathrm{L}^1(I,\mathbb{V}^*).$ Therefore, by the Dunford–Pettis Theorem \cite[Theorem 4.75]{gasinski2016exercises}, $\{f_n\}_{n \in \mathbb{N}}$ is relatively weakly compact in $\mathrm{L}^1(I,\mathbb{V}^*)$. Hence, passing to a subsequence (still denoted by $\{f_n\}_{n \in \mathbb{N}}$), we have
			$$ f_n \xrightarrow{w} f \text{ in } \mathrm{L}^1(I,\mathbb{V}^*)\ \text{ as } \ n \to \infty.$$
			By using Corollary \ref{corollary 4.7}, it implies that 
			\begin{align}\label{conver-4}
				\Theta(f_n) \to \Theta(f)\ \text{ uniformly in } \ \mathbb{V}.
			\end{align}
			Therefore, using \eqref{conver-4} and passing limit in \eqref{Converge1} as $n \to \infty$, we obtain
			$$ y_n(t) \to y(t)=\mathscr{R}(t)q_0+\int_{0}^{t}\mathscr{R}(t-s)Hf(s)ds+\int_{0}^{t}\mathscr{R}(t-s)Bu(s)ds, ~t \in I.$$
			Note that $f_n \in S_{\partial \mathbb{F}}(q_n) \text{ for all } n \in \mathbb{N}$ with $f_n \xrightarrow{w} f \text{ in } \mathrm{L}^1(I, \mathbb{V}^*)$ and $q_n \to q \text{ in } \mathrm{C}(I, \mathbb{V}).$ It follows from \cite[Lemma 5.3]{SMAOMS}, we conclude that $f \in S_{\partial \mathbb{F}}(q).$ Hence, by the definition of the multimap $\Lambda$, we confirm that $y \in \Lambda(q),$  showing that $\Lambda$ is closed. Moreover, by Corollary \ref{corollary 4.7}, the sequence $\{\Lambda(q_n)\}_{n\in \mathbb{N}}$ is
			relatively compact. Since $\Lambda$ has closed values, it follows that the set $\Lambda(q)$ is compact. Consequently, by using \cite[Proposition 1.2]{DK} with \eqref{assum-1}, the map $\Lambda: B_{\mathfrak{R}} \to \mathcal{P}( B_{\mathfrak{R}}) $ is u.s.c. 
			Therefore, by \cite[Theorem 1.2.14]{kamenskii} together with Lemma \ref{lemma-1}, the multivalued map  
			$$\Gamma: B_{\mathfrak{R}} \to \mathcal{K}v(B_{\mathfrak{R}})$$ 
			defined by \eqref{Def.-Gamma} is u.s.c.     
			Finally, it remains to show that $\Gamma$ is $\beta$-condensing.
			To this end, let $W \subset \mathrm{C}(I, \mathbb{V})$ be a bounded set. Then, it follows from Lemma \ref{lemma-2} and condition \textit{(H5)} that
			\begin{equation}\label{ineq-0}
				\beta(\mathcal{G}(W)) \leq M\chi(h(W)) \leq k_h M\beta(W). 
			\end{equation}
			
			Consider the sequences 
			$\{w_n\}_{n\in \mathbb{N}} \subset W $ and $\{f_n\}_{n\in \mathbb{N}} \subset \mathrm{L}^1(I,\mathbb{V}^*)$ with $f_n \in  S_{\partial \mathbb{F}}(w_n)$ for all $n\in\mathbb{N}.$ Since, for $y_n \in \Lambda(w_n),$ we have 
			\begin{align*}
				y_n(t) = \mathscr{R}(t)q_0+\int_{0}^{t}\mathscr{R}(t-s)Hf_n(s)ds+\int_{0}^{t}\mathscr{R}(t-s)Bu(s)ds \ \text{ for all } \ t \in I.
			\end{align*}
			Then we obtain
			\begin{align*}
				\chi(y_n(t)) &= \chi \left( \left\{ \mathscr{R}(t)q_0+\int_{0}^{t}\mathscr{R}(t-s)Hf_n(s)ds+\int_{0}^{t}\mathscr{R}(t-s)Bu(s)ds \colon n \geq 1 \right\}\right) \\
				& \leq \chi\left(\left\{\int_{0}^{t}\mathscr{R}(t-s)Hf_n(s)ds \colon n \geq 1 \right\}\right).
			\end{align*}
			Using \textit{(F)$^{\prime}$}, we have
			\begin{align}\label{bound}
				\chi(\{  f_n(s)\}_{n\in \mathbb{N}})&\leq  \chi\big( \partial \mathbb{F}(s,\{w_n(s)\}_{n\in \mathbb{N}})\big) \nonumber\\
				&\leq   k(s)\chi(\{w_n(s)\}_{n\in \mathbb{N}}) \leq  
				k(s)\chi(W(s)) \leq  k(s) \beta(W).
			\end{align}
			It follows from \eqref{bound} and \cite[Theorem 4.2.2]{kamenskii} that
			\begin{align}\label{ineq-1}
				\beta(\Lambda(W)) \leq 2 M \norm{H}_{\mathcal{L}(\mathbb{V}^*,\mathbb{V})}  \int_0^t k(s) ds \beta(W).
			\end{align}
			Combining \eqref{Def.-Gamma}, \eqref{ineq-0} and \eqref{ineq-1}, we get
			\begin{align*}
				\beta(\Gamma(W)) \leq  M \bigg(k_h+2 \norm{H}_{\mathcal{L}(\mathbb{V}^*,\mathbb{V})} \int_{0}^{a} k(s)ds\bigg) \beta(W),
			\end{align*}
			that implies the map $\Gamma$ is $\beta$-condensing on $B_{\mathfrak{R}}$. By applying \cite[Corollary 3.3.1]{kamenskii}, we obtain that $\Gamma$ has a fixed point.
		\end{proof}
		
		\section{Approximate Controllability of systems}\label{sec4}\setcounter{equation}{0}
		In this section, we concern the approximate controllability of the inclusion system with a general non-local condition given by \eqref{PC1}. First, we recall some preliminary results concerning the approximate controllability for the linear system:
		\begin{equation}\label{ineq-14}
			\begin{cases}
				\begin{aligned}
					\frac{d}{dt}\left(q(t)+\int_0^t N(t-s)q(s)ds\right)&=\mathscr{A}q(t)+\int_0^t G(t-s)q(s)ds+Bu(t) \ t\geq 0,\\
					q(0)&=q_0 \in \mathbb{V}.
				\end{aligned}
			\end{cases}
		\end{equation}
		\begin{definition}
			The system \eqref{ineq-14} is said to be approximately controllable on $I$, if for any $q_0,\tilde{q} \in \mathbb{V}$ and  $\varepsilon>0$, there exist  $u(\cdot) \in \mathrm{L}^2(I, \mathbb{U})$ such that the mild solution $q(\cdot)$ of \eqref{ineq-14} satisfies:
			$$\|q(a)-\tilde{q}\|_{\mathbb{V}} \leq \varepsilon.$$
		\end{definition}
		
		\begin{lemma}[{\cite[Lemma 3.3]{arora2025controllability}}]\label{lemma 5.1}
			The system \eqref{ineq-14} is approximately controllable on $I$ if and only if,
			\begin{equation*}
				\text{for every}~y^*\in \mathbb{V}^*, \ B^*\mathscr{R}(a-t)^*y^*=0, ~\text{for all }~t\in I \implies y^*=0.
			\end{equation*}
		\end{lemma}
		
		\begin{theorem}
			Assume that $(\mathscr{R}(t))_{t > 0}$ be a compact operator and condition (F) holds and the linear system \eqref{ineq-14} is approximately controllable on $[0, t]$ for any $0 < t \leq a$. Then the
			system \eqref{PC1} is approximately controllable on $I$.\end{theorem}
		\begin{proof}
			Fix $\tilde{q} \in \mathbb{V}.$ By Theorem \ref{Theorem-1}, the system  has a unique mild solution on $I$. For $u = 0$, let $q(\cdot) = q(\cdot, q_0, 0)$  denote the corresponding solution, yielding:
			$$ q(t)= \mathscr{R}(t)(q_0+z)+\int_{0}^{t}\mathscr{R}(t-s)Hf(s)ds \ \text{ for all } \ t \in I, $$
			for some $z \in h(q) \text{ and } f \in S_{\partial \mathbb{F}}(q).$ Consider a sequence $0< b_n <a, $ such that $b_n \to a \text{ as } n \to \infty$. We denote $q_n = q(b_n, q_0, 0).$ Since the linear system \eqref{ineq-14} is approximately controllable on $[0, a-b_n]$. Consequently, there is a control function $w_n(\cdot) \in \mathrm{L}^2([0, a-b_n],\mathbb{U})$ such that 
			\begin{align}\label{ineq-15}
				\bigg\|\mathscr{R}(a-b_n)q_n+\int_{0}^{a-b_n}  \mathscr{R}(a-b_n- s) B w_n(s) ds -\tilde{q}\Big\|_{\mathbb{V}} \to 0\  \text{ as } \ n \to \infty.
			\end{align}
			Now, we define
			\begin{align*}
				u_n(s) = 
				\begin{cases}
					0, & 0 \leq s \leq b_n, \\
					w_n(s - b_n), & b_n < s \leq a.
				\end{cases}
			\end{align*}
			Using Theorem \ref{Theorem-1}, we obtain that there exists a solution $y_n(\cdot)$ on $I$ of the integral equation 
			\begin{equation}\label{ineq-16}
				y_n(t) = \mathscr{R}(t)(q_0 + z) + \int_0^t \mathscr{R}(t - s) Hf_n(s) \, ds + \int_0^t \mathscr{R}(t - s) B u_n(s) ds,~ t \in I, \text{ for all }n \in \mathbb{N}.
			\end{equation}
			Since the solution of above equation is unique, we have that $q(s) = y_n(s)$ for all $0 \leq s \leq b_n$. Substituting in \eqref{ineq-16}, we have
			\begin{align*}
				y_n(a) &= \mathscr{R}(a)(q_0 + z) + \int_0^a  \mathscr{R}(a-s)Hf_n(s) ds + \int_0^a  \mathscr{R}(a - s) B u_n(s)ds \\
				&= \mathscr{R}(a)(q_0 + z) + \int_0^{b_n}  \mathscr{R}(a - s)H f_n(s) ds + \int_{b_n}^a  \mathscr{R}(a - s)Hf_n(s) ds + \int_0^a  \mathscr{R}(a - s) B u_n(s)ds.
			\end{align*}
			Finally, we evaluate 
			\begin{align*}
				&\norm{y_n(a)-\tilde{q}}_{\mathbb{V}} \\
				&\leq \Big\| \bigg(\mathscr{R}(a)-\mathscr{R}(a-b_n)\mathscr{R}(b_n)\bigg)(q_0 + z)+\mathscr{R}(a-b_n)\bigg(\mathscr{R}(b_n)(q_0 + z)+\int_0^{b_n}  \mathscr{R}(b_n - s)Hf_n(s) ds\bigg)\\
				&\quad+ \int_0^{b_n} \bigg(\mathscr{R}(a - s)- \mathscr{R}(a-b_n) \mathscr{R}(b_n - s)\bigg)Hf_n(s) ds + \int_{b_n}^a  \mathscr{R}(a - s)Hf_n(s)ds\\
				&\quad +\int_{b_n}^{a}  \mathscr{R}(a- s) B w_n(s-b_n) ds -\tilde{q}\Big\|_{\mathbb{V}} \\
				&\leq \Big\| \bigg(\mathscr{R}(a)-\mathscr{R}(a-b_n)\mathscr{R}(b_n)\bigg)(q_0 + z)\bigg\|_{\mathbb{V}} + \bigg\|\mathscr{R}(a-b_n)q_n+\int_{0}^{a-b_n}  \mathscr{R}(a- b_n-s) B w_n(s) ds -\tilde{q}\Big\|_{\mathbb{V}} \\
				&\quad + \int_0^{b_n} \bigg\|\bigg(\mathscr{R}(a - s)- \mathscr{R}(a-b_n) \mathscr{R}(b_n - s)\bigg)Hf_n(s)\bigg\|_{\mathbb{V}} ds\\
				&\quad + \int_{b_n}^a \bigg\|\mathscr{R}(a - s)Hf_n(s)\bigg\|_{\mathbb{V}}ds \to 0\  \text{ as } \ n \to \infty,
			\end{align*}
			where we have used Lemma \ref{lemma 3.1}, the convergence \eqref{ineq-15} and the Lebesgue dominated convergence theorem. Since $\tilde{q}$ was chosen arbitrarily, thus we infer that the system \eqref{PC1} is approximately controllable on $I$.
		\end{proof}
		
		Now, we initially introduce the notion of duality mapping and the special geometry of Banach spaces with its examples.
		\begin{definition}\label{defn-duality}
			The mapping $\mathcal{J}: \mathbb{V} \to \mathcal{P}(\mathbb{V}^*)$ denotes the duality mapping, which is defined by 
			\begin{equation*}
				\mathcal{J}(y)=\{y^*\in \mathbb{V}^*: \langle y, y^*\rangle=\norm{y}_{\mathbb{V}}^2=\norm{y^*}^2_{\mathbb{V}^*}\} \ \text{ for all } y\in \mathbb{V}.
			\end{equation*}
		\end{definition}
		\begin{definition}
			A normed space $\mathbb{V}$ is said to be uniformly convex if for any $\varepsilon\in (0,2]$ there exists a $\delta>0$ such that if $x,y\in \mathbb{V}$ with $\norm{x}_{\mathbb{V}}\le 1, \norm{y}_{\mathbb{V}}\le 1$ and $\norm{x-y}_{\mathbb{V}}\ge \varepsilon$, then
			\begin{equation*}
				\norm{\frac{x+y}{2}}_{\mathbb{V}}\le 1-\delta.
			\end{equation*}
		\end{definition}
		\begin{ex}[\cite{ZBXGF}] Some examples of uniformly convex spaces are Hilbert spaces, $\mathrm{L}^p([0,\pi],\mathbb{R})$
			for $2 \leq p < \infty$, and Sobolev space $W^{m,p}([0,\pi],\mathbb{R}) \text{ for } m \in \mathbb{N} \cup \{0\},~ 2\leq p < \infty$.
		\end{ex}
		In the remaining part of this article, we assume that $\mathbb{V}$ is a separable  Banach space having a uniformly convex dual, unless otherwise specified. Next, we define the operators
		\begin{equation}\label{2.16}
			\left\{
			\begin{aligned}
				&R(t):= \int_{0}^{t} \mathscr{R}(t-s) B B^* \mathscr{R}(a-s)^* ds,~t\in I, \\[4pt]
				&\mathcal{R}(\varepsilon, R(t)) := (\varepsilon \mathtt{I} + R(t)\mathcal{J})^{-1},~\varepsilon> 0.
			\end{aligned}
			\right.
		\end{equation}
		From the first expression in \eqref{2.16}, we see that the operator 
		$R(t) : \mathbb{V}^{*} \to \mathbb{V}$ is non-negative and symmetric.  Moreover, if $\mathbb{V}$ is a separable Hilbert space 
		identified with its dual, then the duality mapping $\mathcal{J}$ reduces to the identity operator $\mathtt{I}$. In this case, 
		$$\mathcal{R}(\varepsilon, R(t)) := (\varepsilon \mathtt{I} + R(t))^{-1}, ~\varepsilon > 0,$$
		is a linear operator.

		We now establish the controllability result under an additional assumption that the linear system \eqref{ineq-14} is approximately controllable with Hypothesis \textit{(R)}, stated below. This hypothesis introduced in \cite[pp. 4633]{MFPH}, will be employed throughout the proof.
		\begin{itemize} 
			\item[\textit{(R)}] $\norm{R(t)R(a)^{-1}}_{\mathcal{L}(\operatorname{Im}R(a),\mathbb{V})}\le C_0,~t\in I,$ for some constant $C_0\ge 1.$
		\end{itemize}
		\begin{theorem}\label{MR}
			Assume that $(\mathscr{R}(t))_{t > 0}$ be a compact operator and satisfies (R) condition. Let $h : \mathrm{C}(I, \mathbb{V}) \to \mathcal{P}(\mathbb{V})$ satisfies condition (H1)-(H3) and (H6). Additionally, suppose that the conditions (F1), (F2)$^\prime$ are satisfied. Then, for every $\varepsilon >0$, the problem \eqref{PC1} is approximately controllable in $I$.
		\end{theorem}	
		\begin{proof}
			Fix $\varepsilon>0$ and let $\tilde{q} \in \mathbb{V}$ be the target state we want to achieve in time $a$. Let us define the sets 
			\begin{equation*}
				B_{\mathfrak{R}}=\{q\in \mathrm{C}(I, \mathbb{V}): \norm{q}_{\infty} \leq \mathfrak{R}\}, ~\mathfrak{R}>0,
			\end{equation*}
			and
			\begin{equation*}
				Q_{\mathfrak{R}}=\{f\in \mathrm{L}^1(I, \mathbb{V}^*): \norm{f(t)}_{\mathbb{V}^*} \leq \gamma_{\mathfrak{R}}(t)~\text{a.e.}~t\in I\},
			\end{equation*}
			where $\gamma_{\mathfrak{R}}$ is as in Hypothesis \textit{(F2)}$^\prime$. For fixed $z \in h(B_{\mathfrak{R}})$, we define the map $$\Gamma_z^\varepsilon:Q_{\mathfrak{R}}\to \mathrm{C}(I, \mathbb{V})$$ as: 
			\begin{align}\label{CM}
				\Gamma_z^\varepsilon(f)(t)=\mathscr{R}(t)(q_0+z)+\int_{0}^{t}\mathscr{R}(t-s)Hf(s)ds+\int_{0}^{t}\mathscr{R}(t-s)Bu^{\varepsilon}(s)ds,~ t\in I, 
			\end{align}   
			for $f \in  S_{\partial \mathbb{F}}(q)$ and $u^{\varepsilon}\in \mathrm{L}^2(I,\mathbb{U})$ is expressed as
			\begin{equation}\label{Y2}
				u^{\varepsilon}(t)=B^*\mathscr{R}(a-t)^*\mathcal{J}\left((\varepsilon \mathtt{I}+R(a)\mathcal{J})^{-1}Z(z,f)\right),
			\end{equation}	
			where $$Z(z,f)=\tilde{q}-\mathscr{R}(a)(q_0+z)-\int_{0}^{a}\mathscr{R}(a-s)Hf(s)ds.$$ We divide the proof into several steps written as follows:
			\vskip 0.1in 
			\noindent 
			\textbf{STEP-I:}  \emph{The map $\Gamma_z^\varepsilon$ maps the set $Q_{\mathfrak{R}}$ into the set $B_{\mathfrak{R}}$ for some $\mathfrak{R}>0$.} 
			Consider $y^{\varepsilon}= \Gamma_z^\varepsilon(f)$ for some $f \in S_{\partial \mathbb{F}}(q) \text{ and } z \in h(B_{\mathfrak{R}})$. Then we estimate
			\begin{align}\label{Yec}
				\norm{y^{\varepsilon}(t)}_{\mathbb{V}}\leq & M\norm{q_0+z}_{\mathbb{V}}+M\norm{H}_{\mathcal{L}(\mathbb{V}^*,\mathbb{V})}\int_{0}^{t}\norm{f(s)}_{\mathbb{V}^*}ds+\norm{\int_{0}^{t}\mathscr{R}(t-s)Bu^{\varepsilon}(s)ds}_{\mathbb{V}}.
			\end{align}
			Using the expression \eqref{Y2}, we obtain
			\begin{align*}
				\norm{\int_{0}^{t}\mathscr{R}(t-s)Bu^{\varepsilon}(s)ds}_{\mathbb{V}}
				=&\norm{\int_{0}^{t}\mathscr{R}(t-s)BB^*\mathscr{R}(a-t)^*\mathcal{J}\left((\varepsilon \mathtt{I}+R(a)\mathcal{J})^{-1}Z(z,f)\right)ds}_{\mathbb{V}}\\
				\le & \norm{R(t)\mathcal{J}\left(((\varepsilon \mathtt{I}+R(a)\mathcal{J})^{-1}Z(z,f)\right)}_{\mathbb{V}}\\
				\le & \norm{R(t)R(a)^{-1}R(a)\mathcal{J}\left((\varepsilon\mathtt{I}+R(a)\mathcal{J})^{-1}Z(z,f)\right)}_{\mathbb{V}}.	\end{align*}
			Owing to Hypothesis \textit{(R)}, we have 
			\begin{align*}
				\norm{\int_{0}^{t}\mathscr{R}(t-s)Bu^{\varepsilon}(s)ds}_{\mathbb{V}}
				&	\leq  C_0 \norm{R(a)\mathcal{J}\left((\varepsilon\mathtt{I}+R(a)\mathcal{J})^{-1}Z(z,f)\right)}_{\mathbb{V}}\\
				&	\leq  C_0 \norm{(\varepsilon\mathtt{I}+R(a)\mathcal{J}-\varepsilon\mathtt{I})\left((\varepsilon\mathtt{I}+R(a)\mathcal{J})^{-1}Z(z,f)\right)}_{\mathbb{V}} \\
				& \leq  C_0 \norm{Z(z,f)-\varepsilon (\varepsilon\mathtt{I}+R(a)\mathcal{J})^{-1}Z(z,f)}_{\mathbb{V}} \\ 
				& \leq  C_0 (\norm{Z(z,f)}_{\mathbb{V}}+\norm{\varepsilon (\varepsilon\mathtt{I}+R(a)\mathcal{J})^{-1}Z(z,f)}_{\mathbb{V}}). 
			\end{align*}
			By the help of \cite[Lemma 2.2]{NIM}, we deduce
			\begin{align}\label{bounded}
				\norm{\int_{0}^{t}\mathscr{R}(t-s)Bu^{\varepsilon}(s)ds}_{\mathbb{V}} &\leq 2 C_0 \norm{Z(z,f)}_{\mathbb{V}}\nonumber\\
				& \leq 2 C_0 \norm{\tilde{q}-\mathscr{R}(a)(q_0+z)-\int_{0}^{a}\mathscr{R}(a-s)Hf(s)ds}_{\mathbb{V}}\nonumber\\
				& \leq 2 C_0 \left(\norm{\tilde{q}}_{\mathbb{V}}+M(\norm{q_0}_{\mathbb{V}}+\norm{z}_{\mathbb{V}})+M\norm{H}_{\mathcal{L}(\mathbb{V}^*,\mathbb{V})}\int_0^a \norm{f(s)}_{\mathbb{V}^*}ds\right) \nonumber\\
				& \leq 2 C_0 \left(\norm{\tilde{q}}_{\mathbb{V}}+M(\norm{q_0}_{\mathbb{V}}+N_h)+M\norm{H}_{\mathcal{L}(\mathbb{V}^*,\mathbb{V})}\int_0^a \gamma_{\mathfrak{R}}(s)ds\right).
			\end{align}
			Since $z \in h(B_{\mathfrak{R}})$, thus we use hypothesis \textit{(H5)} in the last line to get the final estimate. Therefore, by combining \eqref{Yec} and \eqref{bounded}, we obtain
			\begin{align}\label{18}
				\norm{y^{\varepsilon}(t)}_{\mathbb{V}} \le M\left(\norm{q_0}_{\mathbb{V}}+N_h+\norm{H}_{\mathcal{L}(\mathbb{V}^*,\mathbb{V})}\int_0^a \gamma_{\mathfrak{R}}(s)ds\right)(1+2 C_0)+2 C_0 \norm{\tilde{q}}_{\mathbb{V}}.
			\end{align}
			We deduce from the expression \eqref{18} that 
			\begin{equation*}
				\begin{aligned}
					\|\Gamma_z^\varepsilon(Q_{\mathfrak{R}})\|\leq M\left(\norm{q_0}_{\mathbb{V}}+N_h+\norm{H}_{\mathcal{L}(\mathbb{V}^*,\mathbb{V})}\int_0^a \gamma_{\mathfrak{R}}(s)ds\right)(1+2 C_0)+2 C_0 \norm{\tilde{q}}_{\mathbb{V}}.
				\end{aligned}
			\end{equation*}
			We now assert that there exists $\mathfrak{R}_0\in\mathbb{N}$ such that
			$\Gamma_z^\varepsilon(Q_{\mathfrak{R}_0})\subset B_{\mathfrak{R}_0}$.
			Suppose that this assertion is not true. Then for every
			$\mathfrak{R} \in\mathbb{N}$ there exists $y_{\mathfrak{R}}\in\Gamma_z^\varepsilon(Q_{\mathfrak{R}})$ such that
			$y_{\mathfrak{R}} \notin B_{\mathfrak{R}}$. By the fact that $y_{\mathfrak{R}} \notin B_{\mathfrak{R}}$, we have
			\[
			\begin{aligned}
				\mathfrak{R}
				&< \|\Gamma_z^\varepsilon(Q_{\mathfrak{R}})\| \\
				&\leq M\left(\norm{q_0}_{\mathbb{V}}+N_h+\norm{H}_{\mathcal{L}(\mathbb{V}^*,\mathbb{V})}\int_0^a \gamma_{\mathfrak{R}}(s)ds\right)(1+2 C_0)+2 C_0 \norm{\tilde{q}}_{\mathbb{V}},
			\end{aligned}
			\]
			which implies
			\begin{equation}\label{20}
				\begin{aligned}
					1< \frac{1}{\mathfrak{R}}\|\Gamma_z^\varepsilon(Q_{\mathfrak{R}})\| \leq \frac{1}{\mathfrak{R}}
					\left[
					M\left(\norm{q_0}_{\mathbb{V}}+N_h+\norm{H}_{\mathcal{L}(\mathbb{V}^*,\mathbb{V})}\int_0^a \gamma_{\mathfrak{R}}(s)ds\right)(1+2 C_0)+2 C_0 \norm{\tilde{q}}_{\mathbb{V}}\right].
				\end{aligned}
			\end{equation}
			Owing to Hypotheses \textit{(F2)}$^\prime$, we observe that
			\begin{equation*}
				\liminf_{\mathfrak{R} \to\infty}\frac{1}{\mathfrak{R}}\int_{0}^{a}\gamma_{\mathfrak{R}}(s)ds=0.
			\end{equation*}
			Based on this observation, passing to the limit as
			$\mathfrak{R} \to\infty$ we deduce from \eqref{20} that
			\[
			\begin{aligned}
				1
				&\leq
				\liminf_{\mathfrak{R} \to\infty}
				\frac{1}{\mathfrak{R} }
				\left[M\left(\norm{q_0}_{\mathbb{V}}+N_h+\norm{H}_{\mathcal{L}(\mathbb{V}^*,\mathbb{V})}\int_0^a \gamma_{\mathfrak{R}}(s)ds\right)(1+2 C_0)+2 C_0 \norm{\tilde{q}}_{\mathbb{V}}\right] \\
				&=0,
			\end{aligned}
			\]
			which is a contradiction. Hence, the assertion is true. Thus, we ensure that the operator $\Gamma_z^\varepsilon$ is bounded.
			\vskip 0.1in 
			\noindent  
			\textbf{STEP-II:-}  \emph{The map $\Gamma_z^\varepsilon$ is continuous from ${Q_{\mathfrak{R}}}_w$ into $\mathrm{C}(I, \mathbb{V})$.} Since, the set $Q_{\mathfrak{R}}$ is a compact convex metrizable set in $\mathrm{L}^1(I, \mathbb{V}^*)$ endowed with weak topology. Therefore, it suffices to prove the sequential continuity of the operator $\Gamma_z^\varepsilon$. For proving this, consider a sequence $\{f_n\}_{n\in \mathbb{N}}\subset Q_{\mathfrak{R}}$ such that 
			$$ f_n \xrightarrow{w} f \ \text{ in } \ \mathrm{L}^1(I, \mathbb{V}^*).$$
			Let $y^{\varepsilon}_n=\Gamma_z^\varepsilon(f_n)$, so that the following equation holds:
			\begin{equation}\label{4.71}
				y^{\varepsilon}_n(t)=\mathscr{R}(t)(q_0+z)+\int_{0}^{t}\mathscr{R}(t-s)Hf_n(s)ds+\int_{0}^{t}\mathscr{R}(t-s)Bu^{\varepsilon}_n(s)ds \ \text{ for all } \ t \in I, 
			\end{equation} 
			where   
			\begin{equation*}
				u^{\varepsilon}_n(t)=B^*\mathscr{R}(a-t)^*\mathcal{J}\left((\varepsilon \mathtt{I}+R(a)\mathcal{J})^{-1}\left(\tilde{q}-\mathscr{R}(a)(q_0+z)-\int_{0}^{a}\mathscr{R}(a-s)Hf_n(s)ds\right)\right).
			\end{equation*}
			By virtue of \eqref{compact-con}, we deduce, in particular, that
			\begin{equation*}
				\int_{0}^{a}\mathscr{R}(a-s)Hf_n(s)ds\to \int_{0}^{a}\mathscr{R}(a-s)Hf(s)ds \ \text{ as }\  n\to \infty.
			\end{equation*}
			By the continuity of the operator $(\varepsilon \mathtt{I}+R(a)\mathcal{J})^{-1}$ (see \cite[Lemma 2.24]{arora2025controllability}), we arrive at
			\begin{align*}
				(\varepsilon& \mathtt{I}+R(a)\mathcal{J})^{-1} \bigg(\tilde{q}-\mathscr{R}(a)(q_0+z)-\int_{0}^{a}\mathscr{R}(a-s)Hf_n(s)ds\bigg) \\
				&\to (\varepsilon \mathtt{I}+R(a)\mathcal{J})^{-1} \bigg(\tilde{q}-\mathscr{R}(a)(q_0+z)-\int_{0}^{a}\mathscr{R}(a-s)Hf(s)ds\bigg) \ \text{ in } \ \mathbb{V}\ \text{ as }\ n \to \infty. 
			\end{align*}
			Since the mapping $\mathcal{J} : \mathbb{V} \to \mathbb{V}^*$ is demicontinuous, it is immediate that
			\begin{align*}
				\mathcal{J}\bigg(&(\varepsilon \mathtt{I}+R(a)\mathcal{J})^{-1} \bigg(\tilde{q}-\mathscr{R}(a)(q_0+z)-\int_{0}^{a}\mathscr{R}(a-s)Hf_n(s)ds\bigg) \bigg) \nonumber\\
				&\xrightarrow{w}  
				\mathcal{J}\bigg( (\varepsilon \mathtt{I}+R(a)\mathcal{J})^{-1}  \bigg(\tilde{q}-\mathscr{R}(a)(q_0+z)-\int_{0}^{a}\mathscr{R}(a-s)Hf(s)ds\bigg) \bigg)\ \text{ in } \ \mathbb{V}^* \ \text{ as }\ n \to \infty.
			\end{align*}
			Using the compactness of the semigroup $(\mathscr{R}(t))_{t>0}$ together with the boundedness of the linear operator $B:\mathbb{U}\to\mathbb{V}$, we obtain
			\begin{align*}
				&u^{\varepsilon}_n(t)=B^*\mathscr{R}(a-t)^*\mathcal{J}\left((\varepsilon \mathtt{I}+R(a)\mathcal{J})^{-1}\left(\tilde{q}-\mathscr{R}(a)(q_0+z)-\int_{0}^{a}\mathscr{R}(a-s)Hf_n(s)ds\right)\right)\\
				&\to B^*\mathscr{R}(a-t)^*\mathcal{J}\left((\varepsilon \mathtt{I}+R(a)\mathcal{J})^{-1}\left(\tilde{q}-\mathscr{R}(a)(q_0+z)-\int_{0}^{a}\mathscr{R}(a-s)Hf(s)ds\right)\right)\ \text{ as }\ n\to \infty.
			\end{align*}
			Further, using the Definition \ref{defn-duality} together with \cite[Lemma 2.2]{NIM}, we compute
			\begin{align*}
				\norm{u^{\varepsilon}_n(t)}_{\mathbb{U}} \le& \norm{B^* \mathscr{R}(a-t)^* \mathcal{J}\left((\varepsilon \mathtt{I}+R(a)\mathcal{J})^{-1}\left(\tilde{q}-\mathscr{R}(a)(q_0+z)-\int_{0}^{a}\mathscr{R}(a-s)Hf_n(s)ds\right)\right)}\\
				\le& \frac{1}{\varepsilon} M\|B^*\|_{\mathcal{L}(\mathbb{V}^*, \mathbb{U})} \left\|\left(\varepsilon(\varepsilon \mathtt{I}+R(a)\mathcal{J})^{-1}\left(\tilde{q}-\mathscr{R}(a)(q_0+z)-\int_{0}^{a}\mathscr{R}(a-s)Hf_n(s)ds\right)\right)\right\|_{\mathbb{V}^*}\\
				\leq & \frac{1}{\varepsilon} M \|B^*\|_{\mathcal{L}(\mathbb{V}^*, \mathbb{U})} \bigg\|\tilde{q}-\mathscr{R}(a)(q_0+z)-\int_{0}^{a}\mathscr{R}(a-s)Hf_n(s)ds \bigg\|_{\mathbb{V}} \\       
				\leq & \frac{1}{\varepsilon} M \|B^*\|_{\mathcal{L}(\mathbb{V}^*, \mathbb{U})} \left(\norm{\tilde{q}}_{\mathbb{V}}+M(\norm{q_0}_{\mathbb{V}}+\norm{z}_{\mathbb{V}})+M\norm{H}_{\mathcal{L}(\mathbb{V}^*,\mathbb{V})}\int_{0}^{a}\gamma_{\mathfrak{R}}(s)ds\right),
			\end{align*}
			and by employing \cite[Theorem 2.38]{SMAOMS}, we conclude that 
			$$u^{\varepsilon}_n\to u^{\varepsilon} \ \text{ in } \ \mathrm{L}^2(I,\mathbb{U}).$$
			We now estimate
			\begin{align*}
				&\norm{\int_{0}^{t}\mathscr{R}(t-s)Bu^{\varepsilon}_n(s)ds-\int_{0}^{t}\mathscr{R}(t-s)Bu^{\varepsilon}(s)ds}\\
				&	=\norm{\int_{0}^{t}\mathscr{R}(t-s)(Bu^{\varepsilon}_n(s)-Bu^{\varepsilon}(s))ds}\\
				&	\le M\norm{B}_{\mathcal{L}(\mathbb{U},\mathbb{V})} \int_{0}^{t}\norm{u^{\varepsilon}_n(s)-u^{\varepsilon}(s)}_{\mathbb{U}} ds\\
				&	\le M\norm{B}_{\mathcal{L}(\mathbb{U},\mathbb{V})} a^{\frac{1}{2}}\norm{u^{\varepsilon}_n-u^{\varepsilon}}_{\mathrm{L}^2(I,\mathbb{U})},
			\end{align*}
			where in the last inequality, we use H$\ddot{\text{o}}$lder's inequality.
			From this, we conclude that 
			\begin{equation}\label{4.13}
				\int_{0}^{t}\mathscr{R}(t-s)Bu^{\varepsilon}_n(s)ds\to \int_{0}^{t}\mathscr{R}(t-s)Bu^{\varepsilon}(s)ds \ \text{ in } \ \mathrm{C}(I, \mathbb{V}).
			\end{equation}
			Therefore, passing limit in \eqref{4.71} as $n\to \infty$ and using \eqref{compact-con} and \eqref{4.13} we obtain
			$$ y^{\varepsilon}_n(t) \to y^{\varepsilon}(t)=\mathscr{R}(t)(q_0+z)+\int_{0}^{t}\mathscr{R}(t-s)Hf(s)ds+\int_{0}^{t}\mathscr{R}(t-s)Bu^{\varepsilon}(s)ds \ \text{ for all } \ t \in I.$$
			This shows that the map $\Gamma_z^\varepsilon$ is continuous from ${Q_{\mathfrak{R}}}_w$ into $\mathrm{C}(I, \mathbb{V})$.
			\vskip 0.1in 
			\noindent  
			\textbf{STEP-III:-} Following the same arguments as those in Steps III--VI of Theorem \ref{Theorem-2}, it follows that $ q^{\varepsilon} \in \mathrm{C}(I, \mathbb{V})$ is a mild solution to the problem \eqref{PC1}.
			\vskip 0.1in 
			\noindent 
			\textbf{STEP-IV:-} \emph{It remains to verify that $\norm{q^{\varepsilon}(a)- \tilde{q}}_{\mathbb{V}} \to 0 \text{ as } \varepsilon\to 0^+.$} Since, $q^{\varepsilon}(\cdot)$ is a mild solution of the system  \eqref{PC1}, thus $q^{\varepsilon}(\cdot)$ satisfies 
			\begin{align*}
				q^{\varepsilon}(t)=\mathscr{R}(t)(q_0+z_{\varepsilon})+\int_{0}^{t}\mathscr{R}(t-s)Hf_{\varepsilon}(s)ds+ \int_{0}^{t}\mathscr{R}(t-s)Bu^{\varepsilon}(s)ds,
			\end{align*}
			where $ f_{\varepsilon} \in S_{\partial \mathbb{F}}(q^{\varepsilon})$ and $z_{\varepsilon} \in h(q^{\varepsilon}).$ Recalling the definition of the map $R(a)$, we estimate
			\begin{align}\label{4.58}
				q^{\varepsilon}(a)&=\mathscr{R}(a)(q_0+z_{\varepsilon})+\int_{0}^{a}\mathscr{R}(a-s)Hf_{\varepsilon}(s)ds+ \int_{0}^{a}\mathscr{R}(a-s)Bu^{\varepsilon}(s)ds \nonumber\\
				& = \mathscr{R}(a)(q_0+z_{\varepsilon})+\int_{0}^{a}\mathscr{R}(a-s)Hf_{\varepsilon}(s)ds \nonumber \\
				&\quad+\int_{0}^{a}\mathscr{R}(a-s)BB^*\mathscr{R}(a-s)^*\mathcal{J}\left((\varepsilon \mathtt{I}+R(a)\mathcal{J})^{-1}\bigg(\tilde{q}-\mathscr{R}(a)(q_0+z_{\varepsilon})-\int_{0}^{a}\mathscr{R}(a-s)Hf_{\varepsilon}(s)ds\bigg)\right)ds\nonumber \\
				&= \mathscr{R}(a)(q_0+z_{\varepsilon})+\int_{0}^{a}\mathscr{R}(a-s)Hf_{\varepsilon}(s)ds\nonumber \\
				&\quad+(\varepsilon \mathtt{I}+R(a)\mathcal{J}-\varepsilon \mathtt{I})(\varepsilon \mathtt{I}+R(a)\mathcal{J})^{-1}\bigg(\tilde{q}-\mathscr{R}(a)(q_0+z_{\varepsilon})-\int_{0}^{a}\mathscr{R}(a-s)Hf_{\varepsilon}(s)ds\bigg)\nonumber \\
				&= \mathscr{R}(a)(q_0+z_{\varepsilon})+\int_{0}^{a}\mathscr{R}(a-s)Hf_{\varepsilon}(s)ds+\tilde{q}-\mathscr{R}(a)(q_0+z_{\varepsilon})\nonumber \\
				&-\int_{0}^{a}\mathscr{R}(a-s)Hf_{\varepsilon}(s)ds-\varepsilon(\varepsilon \mathtt{I}+R(a)\mathcal{J})^{-1}\bigg(\tilde{q}-\mathscr{R}(a)(q_0+z_{\varepsilon})-\int_{0}^{a}\mathscr{R}(a-s)Hf_{\varepsilon}(s)ds\bigg)\nonumber \\
				& = \tilde{q}-\varepsilon(\varepsilon \mathtt{I}+R(a)\mathcal{J})^{-1}\bigg(\tilde{q}-\mathscr{R}(a)(q_0+z_{\varepsilon})-\int_{0}^{a}\mathscr{R}(a-s)Hf_{\varepsilon}(s)ds\bigg).
			\end{align}
			By Hypothesis \textit{(H6)}, $z_{\varepsilon} \in h(q^{\varepsilon})$ is bounded in $\mathbb{V}$. Since $\mathbb{V}$ is reflexive, there exists a subsequence (still denoted by the same symbol) such that 
			\begin{align*}
				z_{\varepsilon}  \xrightarrow{w} z\  \text{ in } \ \mathbb{V}.
			\end{align*}
			Moreover, $f_{\varepsilon} \in S_{\partial \mathbb{F}}(q^{\varepsilon})$. By observing \eqref{18}, $q^{\varepsilon}$ is uniformly bounded. Thus, in view of Hypothesis \textit{(F2)}$^\prime$, the sequence $f_{\varepsilon}$ is uniformly integrable in $\mathrm{L}^1(I,\mathbb{V}^*)$ and arguing similar as in STEP-VI of Theorem \ref{Theorem-2}, we deduce that	
			\begin{align*}
				f_{\varepsilon} \xrightarrow{w} f \ \text{ in } \ \mathrm{L}^1(I,\mathbb{V}^*).
			\end{align*}
			Therefore, using the compactness of $(\mathscr{R}(t))_{t > 0}$, we infer \begin{equation}\label{conver-g1}
				\tilde{q}-\mathscr{R}(a)(q_0+z_{\varepsilon} )-\int_{0}^{a}\mathscr{R}(a-s)Hf_{\varepsilon} (s)ds \to \tilde{q}-\mathscr{R}(a)(q_0+z)-\int_{0}^{a}\mathscr{R}(a-s)Hf(s)ds\ \text{ in }\ \mathbb{V}.
			\end{equation}
			Hence, by using the convergence \eqref{conver-g1} together with \cite[Lemma 4.4]{MFPH} in \eqref{4.58}, we obtain $q^{\varepsilon}(a)\to \tilde{q}$ as desired.
		\end{proof}
		
		\begin{theorem}
			Assume that $(\mathscr{R}(t))_{t > 0}$ is a compact operator that satisfies the (R) condition. Let $h : \mathrm{C}(I, \mathbb{V}) \to \mathcal{P}(\mathbb{V})$ satisfies condition (H4), (H6) and \eqref{haus-h}. Additionally, suppose that the conditions (F2)$^\prime$, (F) and 
			\begin{align}\label{condition-con}
				M(2 C_0+1)\bigg(L_{h}+ \norm{H}_{\mathcal{L}(\mathbb{V}^*,\mathbb{V})} \int_{0}^{a} \eta(s)ds\bigg) <1,
			\end{align}
			are satisfied. Then, for every $\varepsilon >0$, the problem \eqref{PC1} is approximately controllable in $I$.
		\end{theorem}
		\begin{proof} 
			Let us now define the multivalued operator $\Gamma^{\varepsilon}: \mathrm{C}(I, \mathbb{V}) \to \mathcal{P}(\mathrm{C}(I, \mathbb{V}))$ by 
			\begin{equation*}
				\Gamma^{\varepsilon}(q)= \mathcal{G}(q)+\Lambda^{\varepsilon}(q),
			\end{equation*}
			where $\mathcal{G}$ is defined in Lemma \ref{lemma-1} and $\Lambda^{\varepsilon} $ is defined as similar in Theorem \ref{Theorem-1} by using value of control given in \eqref{Y2}. We divide the proof into several steps written as follows:
			\vskip 0.1in 
			\noindent 
			\textbf{STEP-I:-} \emph{Show that $\Gamma^{\varepsilon}$ maps bounded sets into bounded sets.} For this, let $q \in \mathrm{C}(I, \mathbb{V})$ with $\norm{q}_{\infty} \leq \mathfrak{R}, \text{ and } y^{\varepsilon} \in \Gamma^{\varepsilon}(q)$. Thus, from \eqref{Def-G} and \eqref{def-y}, we have
			$$	\norm{y^{\varepsilon}(t)}_{\mathbb{V}} \leq M\norm{q_0+z}_{\mathbb{V}}+M\norm{H}_{\mathcal{L}(\mathbb{V}^*,\mathbb{V})}\int_{0}^{t}\norm{f(s)}_{\mathbb{V}^*}ds+\norm{\int_{0}^{t}\mathscr{R}(t-s)Bu^{\varepsilon}(s)ds}_{\mathbb{V}}.$$
			By applying the same argument as in STEP-I of Theorem \ref{MR}, we conclude that $y^{\varepsilon} \in B_{\mathfrak{R}}$. 
			\vskip 0.1in 
			\noindent 
			\textbf{STEP-II:-} \emph{The multivalued operator $\Gamma^{\varepsilon}$ has closed and bounded values.} Firstly, we prove that the multivalued map $\Lambda^{\varepsilon}$ has compact values. For this, it suffices to show that the operator $\Lambda^{\varepsilon}$ is completely continuous and has a closed graph. We begin by proving that it's graph is closed. For this, let $\{q_n\}_{n\in\mathbb{N}}\subset B_{\mathfrak{R}}$ and $\{y^{\varepsilon}_n\}_{n\in\mathbb{N}}\subset \mathrm{C}(I,\mathbb{V})$ be sequences such that
			$$ q_n\to q, \quad y^{\varepsilon}_n\to y^{\varepsilon} \ \text{ with } \  y^{\varepsilon}_n\in\Lambda^{\varepsilon}(q_n).$$
			Since $ y^{\varepsilon}_n\in\Lambda^{\varepsilon}(q_n)$, there exists $f_n\in S_{\partial\mathbb{F}}(q_n)$ such that
			\begin{align}\label{Conver-11}
				y^{\varepsilon}_n(t)=\mathscr{R}(t)q_0+\int_{0}^{t}\mathscr{R}(t-s)Hf_n(s)ds+\int_{0}^{t}\mathscr{R}(t-s)Bu^{\varepsilon}_n(s)ds, ~t\in I. 
			\end{align}
			By assumption \textit{(F2)$^\prime$}, the sequence $\{f_n\}_{n \in \mathbb{N}}$ is uniformly integrable in $\mathrm{L}^1(I,\mathbb{V}^*).$ Thus, by using Dunford–Pettis Theorem \cite[Theorem 4.75]{gasinski2016exercises}, we can find a subsequence of $\{f_n\}_{n \in \mathbb{N}}$ (labeled with the same symbol) such that
			$$ f_n \xrightarrow{w} f \ \text{ in } \ \mathrm{L}^1(I,\mathbb{V}^*).$$
			By following the same argument as in STEP-II of Theorem \ref{MR} and passing limit in \eqref{Conver-11} as $n \to \infty$, we obtain
			$$ y^{\varepsilon}_n(t) \to y^{\varepsilon}(t)=\mathscr{R}(t)q_0+\int_{0}^{t}\mathscr{R}(t-s)Hf(s)ds+\int_{0}^{t}\mathscr{R}(t-s)Bu^{\varepsilon}(s)ds, ~t \in I.$$
			Since $f_n \in S_{\partial \mathbb{F}}(q_n) \text{ for all } n \in \mathbb{N}$ with $f_n \xrightarrow{w} f \text{ in } \mathrm{L}^1(I, \mathbb{V}^*)$ and $q_n \to q \text{ in } \mathrm{C}(I, \mathbb{V}),$ therefore, by \cite[Lemma 5.3]{SMAOMS}, we conclude that $f \in S_{\partial \mathbb{F}}(q).$ Hence, the graph of $\Lambda^{\varepsilon}$ is closed. Next, by proceeding as in the proof of \cite[Theorem~4.1]{arora2025controllability}, with the necessary modifications, we conclude that the multivalued map $\Lambda^{\varepsilon} $ is completely continuous. Consequently, these two facts imply that $\Lambda^{\varepsilon}$ has compact values. Thus, by following the proof of STEP-II of Theorem \ref{Theorem-1}, the multivalued operator $\Gamma^{\varepsilon}$ is closed and bounded valued.
			\vskip 0.1in 
			\noindent 
			\textbf{STEP-III:-} \emph{The multivalued map $\Gamma^{\varepsilon}$  is a set contraction.} For this, consider $p, q, \Tilde{p},\Tilde{q} \in \mathrm{C}(I, \mathbb{V})$ such that $\Tilde{p} \in  \Gamma^{\varepsilon}(p) \text{ and } \Tilde{q} \in  \Gamma^{\varepsilon}(q)$ such that 
			\begin{align*}
				\Tilde{p}(t) &=\mathscr{R}(t)z_p+ \mathscr{R}(t)q_0+\int_{0}^{t}\mathscr{R}(t-s)Hf_1(s)ds+\int_{0}^{t}\mathscr{R}(t-s)Bu^{\varepsilon}_p(s)ds,\\
				\Tilde{q}(t) &=\mathscr{R}(t)z_q+ \mathscr{R}(t)q_0+\int_{0}^{t}\mathscr{R}(t-s)Hf_2(s)ds+\int_{0}^{t}\mathscr{R}(t-s)Bu^{\varepsilon}_q(s)ds,
			\end{align*}
			for some $z_p\in h(p) \text{ and } z_q \in h(q)$. Hence, we can estimate
			\begin{align}\label{estimate-221}
				&\norm{ \Tilde{q}(t)- \Tilde{p}(t)}_{\mathbb{V}}   \nonumber \\
				&\leq M\norm{z_q-z_p}_{\mathbb{V}}+\norm{\int_{0}^{t}\mathscr{R}(t-s)H(f_2(s)-f_1(s))ds}_{\mathbb{V}}+\norm{\int_{0}^{t}\mathscr{R}(t-s)B(u^{\varepsilon}_q(s)-u^{\varepsilon}_p(s))ds}_{\mathbb{V}}.
			\end{align}
			Now, by using the control given in \eqref{Y2}, condition \textit{(R)} and \cite[Lemma 2.2]{NIM}, we calculate 
			\begin{align}\label{4.63}
				&\left\|\int_{0}^{t}\mathscr{R}(t-s)B u^{\varepsilon}_q(s) ds-
				\int_{0}^{t}\mathscr{R}(t-s)B u^{\varepsilon}_p(s) ds\right\|_{\mathbb{V}} \nonumber\\
				&\le \Biggl\|
				\int_{0}^{t}\mathscr{R}(t-s)BB^{*}\mathscr{R}(a-s)^{*}
				\mathcal{J}(\varepsilon\mathtt{I}+R(a)\mathcal{J})^{-1}
				\Bigl(Z(z_q,f_2)-Z(z_p,f_1)\Bigr) ds\Biggr\|_{\mathbb{V}} \nonumber\\
				&\le
				\left\|R(t)R(a)^{-1}R(a)\mathcal{J}
				(\varepsilon\mathtt{I}+R(a)\mathcal{J})^{-1}
				\Bigl(Z(z_q,f_2)-Z(z_p,f_1)\Bigr)
				\right\|_{\mathbb{V}} \nonumber\\
				& \leq C_0\left\|(\varepsilon\mathtt{I}+R(a)\mathcal{J}-\varepsilon\mathtt{I})
				(\varepsilon\mathtt{I}+R(a)\mathcal{J})^{-1}
				\Bigl(Z(z_q,f_2)-Z(z_p,f_1)\Bigr)
				\right\|_{\mathbb{V}} \nonumber\\
				& \leq C_0\left\|Z(z_q,f_2)-Z(z_p,f_1)
				-\varepsilon(\varepsilon\mathtt{I}+R(a)\mathcal{J})^{-1}
				\Bigl(Z(z_q,f_2)-Z(z_p,f_1)\Bigr)
				\right\|_{\mathbb{V}} \nonumber\\
				& \leq 2 C_0\left\|Z(z_q,f_2)-Z(z_p,f_1)\right\|_{\mathbb{V}}\nonumber\\
				& \leq 2 C_0\left\|\mathscr{R}(a)(z_q-z_p)-\int_{0}^{a}\mathscr{R}(a-s)H(f_2(s)-f_1(s))ds \right\|_{\mathbb{V}} \nonumber\\
				& \leq 2 C_0 M\left(\|z_q-z_p\|_{\mathbb{V}}+\norm{H}_{\mathcal{L}(\mathbb{V}^*,\mathbb{V})}\int_{0}^{a}\|f_2(s)-f_1(s)\|_{\mathbb{V}^*}ds \right)
				\nonumber\\
				& \leq 2 C_0 M\left(\left(L_{h}+ \norm{H}_{\mathcal{L}(\mathbb{V}^*,\mathbb{V})}\int_{0}^{a}\eta(s)ds\right)\norm{q-p}_{\infty} + (a+1) \epsilon  \right)
			\end{align}
			where we are using \eqref{estimate-20} and \eqref{estimate-21} in the last inequality. Using \eqref{4.63} in \eqref{estimate-221}, we obtain 
			\begin{align*}
				&\norm{ \Tilde{q}(t)- \Tilde{p}(t)}_{\mathbb{V}}   \nonumber \\
				&\leq (2 C_0+1) M\left(\left(L_{h}+ \norm{H}_{\mathcal{L}(\mathbb{V}^*,\mathbb{V})}\int_{0}^{a}\eta(s)ds\right)\norm{q-p}_{\infty} + (a+1) \epsilon  \right).
			\end{align*}
			Since $\epsilon>0$ is arbitrarily chosen, we can assert that
			\begin{align*}
				\sup_{\Tilde{q} \in \Gamma^{\varepsilon}(q)}d(\Tilde{q},\Gamma^{\varepsilon}(p))\leq M (2 C_0+1)\bigg(L_{h}+ \norm{H}_{\mathcal{L}(\mathbb{V}^*,\mathbb{V})} \int_{0}^{a} \eta(s)ds\bigg) \norm{q-p}_{\infty}.
			\end{align*}
			Similarly, we get 
			\begin{align*}
				\sup_{\Tilde{p} \in \Gamma^{\varepsilon}(p)}d(\Tilde{p},\Gamma^{\varepsilon}(q))\leq M(2 C_0+1)\bigg(L_{h}+ \norm{H}_{\mathcal{L}(\mathbb{V}^*,\mathbb{V})} \int_{0}^{a} \eta(s)ds\bigg) \norm{q-p}_{\infty}.
			\end{align*}
			Repeating this argument, we can affirm that
			\begin{align*}
				d_H{(\Gamma^{\varepsilon}(q),\Gamma^{\varepsilon}(p))}\leq M(2 C_0+1)\bigg(L_{h}+ \norm{H}_{\mathcal{L}(\mathbb{V}^*,\mathbb{V})} \int_{0}^{a} \eta(s)ds\bigg) \norm{q-p}_{\infty},
			\end{align*}
			and combining with \eqref{condition-con}, it follows that $\Gamma^{\varepsilon}$ is a set contraction. Applying Theorem \ref{Fixed-1}, we infer that $\Gamma^{\varepsilon}$ has a fixed point $q^{\varepsilon}(\cdot)\in \mathrm{C}(I, \mathbb{V}).$
			\vskip 0.1in 
			\noindent 
			\textbf{STEP-IV:-} \emph{It remains to verify that $\norm{q^{\varepsilon}(a)- \tilde{q}}_{\mathbb{V}} \to 0 \text{ as } \varepsilon\to 0^+.$} By repeating the arguments presented in STEP-IV of Theorem \ref{MR}, we conclude that the system \eqref{PC1} is approximately controllable on $I$.
		\end{proof}	
		\begin{remark}
			Suppose that assumption \textit{(F2)$^\prime$} in the above theorem is replaced by
			$$\|\partial  \mathbb{F}(t, 0)\|_{\mathbb{V}^*}\in \mathrm{L}^1(I,\mathbb{R}^+),$$
			then the assertion of the above theorem is still valid.
		\end{remark}

		\begin{remark}
			The conclusion of Theorem \ref{Theorem-3} remains valid under the weaker assumption \textit{F(2)}$^{\prime}$ in place of Hypothesis \textit{(F)}$^{\prime}$ by using compactness of the operator $(\mathscr{R}(t))_{t > 0}$. Indeed, by \textit{F(2)}$^{\prime}$, we have
			$$
			f_n \xrightarrow{w}  f \ \text{in } \ \mathrm{L}^1(I,\mathbb{V}^*) \ \text{ as } \ n\to\infty.$$
			Moreover, since $(\mathscr{R}(t))_{t > 0}$ is a compact operator, the multivalued map $\Lambda^{\varepsilon}$ defined similar as in Theorem \ref{Theorem-3} by using value of control given in \eqref{Y2} has compact values. Consequently, by using the fact $\Lambda^{\varepsilon}$ is completely continous, we have
			$$ \beta(\Lambda^{\varepsilon}(W))=0.$$
			Hence, $\Gamma^{\varepsilon}$ is a $\beta$-condensing multivalued operator on $B_{\mathfrak{R}}$ whenever $Mk_h<1. $
			
			Furthermore, since $h:\mathrm{C}(I,\mathbb{V})\to\mathcal{K}v(\mathbb{V})$ is u.s.c. and satisfies condition $\textit{(H6)}$, the arguments employed in the proof of Theorem \ref{Theorem-3} applies to show that $\Gamma^{\varepsilon}$  admits a fixed point $q^{\varepsilon}(\cdot)\in \mathrm{C}(I, \mathbb{V}).$ Finally, repeating the arguments presented in STEP-IV of Theorem \ref{MR}, we conclude that the system \eqref{PC1} is approximately controllable on $I$.
		\end{remark}
		\section{Applications}\label{sec5}\setcounter{equation}{0}
		In this section, we first investigate the existence of a mild solution for a transport equation that generates a non-compact resolvent family and demonstrate the applicability of our developed results. 
		\begin{ex}\label{ex 3.11}
			Consider the following class of equations:
			\begin{align}\label{5.10}
				\dfrac{\partial w(t,\xi)}{\partial t}= \dfrac{\partial w(t,\xi)}{\partial \xi}	+ \displaystyle\int_{0}^{t} e^{-\beta (t-s)}w(s,\xi) ds+\widetilde{\eta}(t,\xi)+ \int_0^{\pi} \mathscr{L}(\xi,\omega) \widetilde{f}(t,\omega)d \omega \ \text{ in }\ (0,a] \times [0,\pi],
			\end{align}
			where
			\begin{align}\label{5.20}
				\widetilde{f}(t,\xi) \in \partial \widetilde{\mathbb{F}}(t,w(t,\xi)) \ \text{ for a.e. }\ t \in (0,a],~ \xi \in [0,\pi].
			\end{align}
			The system is supplemented with the initial and non-local initial conditions
			\begin{align*}
				&w(t,0)=0, ~ t \in [0,a], \\
				&w(0,\xi)\in h((w(\cdot,\xi)),~ \xi \in [0,\pi].
			\end{align*}
			In the system \eqref{5.10}, $\mathscr{L} \in \mathrm{C}([0,\pi] \times [0,\pi],\mathbb{R})$ and define the functional $\widetilde{\mathbb{F}}:[0,a] \times \mathrm{L}^p([0,\pi],\mathbb{R})\to \mathbb{R}$ as follows:
			$$\widetilde{\mathbb{F}}(t,w):=\int_0^{\pi} \widetilde{j}\big(t,w(\tau)\big)d\tau, $$
			where $ \widetilde{j}: [0,a] \times \mathbb{R} \to \mathbb{R}$ is defined as 
			$$ \widetilde{j}\big(t,w(\tau)\big)=\delta(t) \widetilde{g}\big(w(\tau)\big).$$ 
			In the above expression, $\delta : [0,a] \to \mathbb{R}^+$ is a non-negative measurable function with 
			$ \delta \in \mathrm{L}^1([0,a],\mathbb{R}^+),$ and $ \widetilde{g} : \mathbb{R} \to \mathbb{R}$ is defined as 
			$$ \widetilde{g}(r) = \int_{0}^{|r|} \widetilde{\mu}(t)dt,$$
			which is locally Lipschitz, where $\widetilde{\mu}$ is Lipschitz and bounded, that is, $|\widetilde{\mu}(t)| \leq \widetilde{c}$, and $\widetilde{\mu}(0) \leq \widetilde{\mu}_0$ for some positive constant $\widetilde{\mu}_0$. Let us now verify the conditions \textit{(F1)}-\textit{(F2)} and \textit{(F)$^{\prime}$} are satisfied by the functional $\widetilde{\mathbb{F}}$ written as follows:
			\begin{itemize}
				\item[\textit{(F1)}]  Since $\delta : [0,a] \to \mathbb{R}^+ $ is a measurable function, it follows that $$\widetilde{\mathbb{F}}(\cdot,w) \text{ is measurable on } [0,a] \text{ for all } w \in \mathrm{L}^p([0,\pi],\mathbb{R}).$$ Moreover, $\widetilde{\mathbb{F}}(t,\cdot)$ is locally Lipschitz on bounded subsets of $ \mathrm{L}^p([0,\pi],\mathbb{R})$ for \mbox{a.e.} $t \in [0,a]$ which follows from the locally Lipschitz property of the function $\widetilde{g}$, as demonstrated below:
				\begin{align*}
					|\widetilde{\mathbb{F}}(t,w_1) - \widetilde{\mathbb{F}}(t,w_2)| &\leq |\delta(t)| |\widetilde{g}(w_1(\tau))- \widetilde{g}(w_2(\tau))| \\
					&\leq \delta(t) \widetilde{M}\int_{0}^{\pi} |w_1(\tau) - w_2(\tau)| d\tau  \\
					&\leq
					\delta(t) \widetilde{M} (\pi)^{\frac{1}{p'}}\|w_1 - w_2\|_{ \mathrm{L}^p([0,\pi],\mathbb{R})},
				\end{align*}
				where $\widetilde{M}  > 0$ is a locally Lipschitz constant for $\widetilde{g}$ and $\frac{1}{p}+\frac{1}{p'}=1$. 
				\item[\textit{(F2)}] By applying the definition of subdifferential given in \eqref{subdiff-1}, $\partial \widetilde{g}(r)$ is given by
				\begin{equation*}
					\partial \widetilde{g}(r)=
					\begin{cases}
						\begin{aligned}
							[-\widetilde{\mu}_0,\widetilde{\mu}_0],~ & r = 0,\\[2pt]
							\widetilde{\mu}(|r|)\mathrm{sgn}(r),~ & r \neq 0. 
						\end{aligned}
					\end{cases}
				\end{equation*}
				Next, using the definition of the generalized gradient of $\widetilde{j}(t,\cdot)$ (see \cite[Proposition 3.44]{SMAOMS}), we obtain 
				\begin{align}\label{bound-j-11}
					|\partial \widetilde{j}(t,w(\tau))|_{\mathbb{R}}=|\delta(t)\partial \widetilde{g}(w(\tau))|\leq \delta(t)\max\{2\widetilde{\mu}_0,\widetilde{c}\} \ \mbox{ a.e.} \ t \in [0,a].
				\end{align}
				Furthermore, by employing an argument analogous to the proof of part (iv) of \cite[Theorem 3.47]{SMAOMS}, we obtain
				\begin{align*}
					\mathbb{\widetilde{\mathbb{F}}}^0(t,w;v) \leq \int_0^{\pi} \widetilde{j}^0\big(t,w(\tau);v(\tau)\big)d\tau.
				\end{align*} 
				Now, consider $z^* \in \partial \widetilde{\mathbb{F}}(t,w)$ and from the definition of $\partial \widetilde{\mathbb{F}}(t,w),$ it follows that
				$$ \langle z^* , v\rangle_{\mathrm{L}^{p'}([0,\pi],\mathbb{R}),\mathrm{L}^{p}([0,\pi], \mathbb{R})} \leq \widetilde{\mathbb{F}}^{0}(t,w;v)  \leq \int_{0}^{\pi} \widetilde{j}^{0}\big(t,w(\tau);v(\tau)\big)d\tau \  \text{ for all }\ v \in\mathrm{L}^{p}([0,\pi],\mathbb{R}).$$ 
				By applying an argument similar to part (v) of \cite[Theorem 3.47]{SMAOMS}, we arrive at the inclusion
				\begin{align*}
					\partial \widetilde{\mathbb{F}}(t,w) \subseteq \int_0^{\pi} \partial \widetilde{j}\big(t,w(\tau)\big)d\tau \ \text{ for a.e. }\ \tau \in [0,\pi],
				\end{align*}
				which means there exists 
				$\zeta \in \mathrm{L}^{p'}([0,\pi], \mathbb{R})$ such that 
				$\zeta(\tau) \in \partial \widetilde{j}(t,w(\tau))$ for a.e. $ \tau \in [0,\pi]$ and
				\begin{align}\label{equality-1}
					\langle z^*, v\rangle_{\mathrm{L}^{p'}([0,\pi],\mathbb{R}),\mathrm{L}^{p}([0,\pi], \mathbb{R})}
					= \int_0^{\pi} \zeta(\tau)v(\tau)d \tau \ \text{ for all }\ v \in \mathrm{L}^{p}([0,\pi],\mathbb{R}).
				\end{align}

				Since the function $\partial \widetilde{j}$ satisfies the growth condition \eqref{bound-j-11}, we have
				\begin{align}\label{equality-21}
					\|\widetilde{\zeta} \|_{\mathrm{L}^{p'}([0,\pi],\mathbb{R})}^{p'}
					= \int_{0}^{\pi} |\widetilde{\zeta} (\tau)|^{p'}d \tau \leq\int_{0}^{\pi}  (\max\{2\widetilde{\mu}_0,\widetilde{c}\}\delta(t))^{p'} d \tau
					\leq (\max\{2\widetilde{\mu}_0,\widetilde{c}\} \delta(t))^{p'}\pi.
				\end{align}
				By combining \eqref{equality-1} and \eqref{equality-21}, we get $\|z^*\|_{\mathrm{L}^{p'}([0,\pi],\mathbb{R})} \leq \max\{2\widetilde{\mu}_0,\widetilde{c}\} \delta(t)(\pi)^{\frac{1}{p'}},$  which shows that condition \textit{(F2)} is satisfied.
				
				\item[\textit{(F)$^{\prime}$}]  Firstly, we show that 
				\begin{align}\label{5.830}
					d_{H}(\partial \widetilde{g}(r_1),\partial \widetilde{g}(r_2))
					\leq C |r_1-r_2| \ \text{ for all } \ r_1,r_2\in \mathbb{R}.   
				\end{align}
				To verify \eqref{5.830}, we consider the following cases.
				\vskip 0.2cm
				\noindent
				\textbf{Case 1:} \emph{$r_1,r_2 \neq 0$ and both have the same sign.} In this case, 
				$$\partial \widetilde{g}(r_i)=\widetilde{\mu}(|r_i|)\mathrm{sgn}(r_i)\ \text{ for }\ i=1,2 .$$ 
				Let $a_i\in \partial \widetilde{g}(r_i)$ for $i=1,2$. Then
				\[
				a_i=\widetilde{\mu}(|r_i|)\mathrm{sgn}(r_i).
				\]
				Hence,
				$$
				|a_1-a_2| =|\widetilde{\mu}(|r_1|)-\widetilde{\mu}(|r_2|)| \leq L ||r_1|-|r_2||  \leq L |r_1-r_2|,$$
				since $\widetilde{\mu}$ is  Lipschitz.
				\vskip 0.2cm
				\noindent
				\textbf{Case 2:} \emph{$r_1\neq 0 \text{ and } r_2 = 0$.} In this case,
				$$
				\partial \widetilde{g}(r_1)=\{a_1\},~
				\partial \widetilde{g}(0)=[-\widetilde{\mu}_0,\widetilde{\mu}_0].$$
				Thus,
				\begin{equation}\label{dh}
					d_{H}(\partial \widetilde{g}(r_1),\partial \widetilde{g}(0))=
					\begin{cases}
						\begin{aligned}
							0,~ & |a_1| \leq \widetilde{\mu}_0,\\[2pt]
							|a_1|-\widetilde{\mu}_0 ,~ & |a_1| > \widetilde{\mu}_0. 
						\end{aligned}
					\end{cases}
				\end{equation}
				Using \eqref{dh} and the condition $\widetilde{\mu}(0)\le \widetilde{\mu}_0$, we obtain
				$$ 
				d_{H}(\partial \widetilde{g}(r_1),\partial \widetilde{g}(0))=|a_1|-\widetilde{\mu}_0 \leq |\widetilde{\mu}(|r_1|)|-\widetilde{\mu}(0) \leq |\widetilde{\mu}(|r_1|)-\widetilde{\mu}(0)| \leq L |r_1|.$$
				\vskip 0.2cm
				\noindent
				\textbf{Case 3:} \emph{$r_1 >0$ and $r_2 <0$.} Since, $r_1 >0$ and $r_2 <0$, then
				$$ a_1=\widetilde{\mu}(r_1),~ a_2=-\widetilde{\mu}(-r_2).$$ Note that, there exist an $\epsilon>0$ such that $|r_1|,|r_2|\geq\epsilon$ and $|r_1-r_2|\geq r_1-r_2=r_1+|r_2|\geq 2\epsilon$. This implies that $\frac{\mu_0}{\epsilon}|r_1-r_2|\geq 2\mu_0$.
				We now calculate
				\begin{align*}
					|a_1-a_2|& \leq |\widetilde{\mu}(r_1)-\widetilde{\mu}(0)|+|\widetilde{\mu}(|r_2|)-\widetilde{\mu}(0)|+2|\widetilde{\mu(0)}|
					\nonumber\\&	\leq L(r_1+|r_2|)+2\mu_0 = \left(L+\frac{\mu_0}{\epsilon}\right)|r_1-r_2|.
				\end{align*}
				Finally, from \eqref{5.830}, we have 
				\begin{align}\label{5.870}
					d_{H} (\partial  \widetilde{j}\big(t,r_1),\partial \widetilde{j}\big(t,r_2))\leq C\delta(t)|r_1-r_2|.
				\end{align}
				The Clarke subdifferential of 
				$\widetilde{\mathbb{F}}$ satisfies
				$$
				\partial \widetilde{\mathbb{F}}(t,w) \subseteq
				\left\{
				\widetilde{\zeta} \in \mathrm{L}^{p'}([0,\pi], \mathbb{R}) :
				\widetilde{\zeta}(\tau)\in \partial \widetilde{j}\big(t,w(\tau)\big)
				\text{ a.e. } \tau\in[0,\pi]
				\right\}.$$
				Consider, $\widetilde{\zeta} _1\in \partial \widetilde{\mathbb{F}}(t,w_1),
				~\widetilde{\zeta} _2\in \partial \widetilde{\mathbb{F}}(t,w_2).
				$ Then $ \widetilde{\zeta}_1(\tau)\in \partial \widetilde{j}(t,w_1(\tau)),~
				\widetilde{\zeta}_2(\tau)\in \partial \widetilde{j}(t,w_2(\tau)), $ for a.e. $\tau\in[0,\pi],$ respectively. Thus, by using \eqref{5.870}, we obtain
				$$
				|\widetilde{\zeta}_1(\tau)-\widetilde{\zeta}_2(\tau)| \leq C\delta(t)|w_1(\tau)-w_2(\tau)| \text{ a.e. } \tau\in[0,\pi].$$
				Now, squaring and integrating over $[0,\pi]$, we obtain
				$$\|\widetilde{\zeta}_1-\widetilde{\zeta}_2\|_{\mathrm{L}^{p'}}\leq C\delta(t)\|w_1-w_2\|_{\mathrm{L}^{p'}}.$$
				Hence,
				\begin{align}\label{hd-1}
					d_H\big(\partial \widetilde{\mathbb{F}}(t,w_1),
					\partial \widetilde{\mathbb{F}}(t,w_2)\big)
					\leq C\delta(t) \|w_1-w_2\|_{\mathrm{L}^{p'}}.
				\end{align}
				Now, let $K \in \mathrm{L}^{p'}([0,\pi],\mathbb{R}) $ be bounded set. Then, by the definition of $\chi(K),$ there exist finitely many points $w_1, \ldots, w_n \in K$ and $r>0$ such that 
				$$ K \subset \bigcup_{i=1}^{n} B(w_i,r),~r < \chi(K) + \varepsilon.$$
				Thus, for any $w \in K$, pick $i$ with $\|w - w_i\|_{\mathrm{L}^{p'}} < r$. Then, from \eqref{hd-1}, we get 
				$$d_H(\partial \widetilde{\mathbb{F}}(t,w), \partial \widetilde{\mathbb{F}}(t,w_i)) \leq C\delta(t)r .$$
				Thus,
				$$
				\partial \widetilde{\mathbb{F}}(t,K) = \bigcup_{w \in K} \partial \widetilde{\mathbb{F}}(t,w) \subset \bigcup_{i=1}^{n}  N_{C\delta(t)r}\big(\partial \widetilde{\mathbb{F}}(t,w_i)\big),$$ 
				where the notation $$N_{C\delta(t)r}\big(\partial \widetilde{\mathbb{F}}(t,w)\big):=\left\{x\in \mathrm{L}^{p'}([0,\pi], \mathbb{R}):
				d\big(x,\partial\widetilde{\mathbb F}(t,w)\big)
				<C\delta(t)r\right\}$$ with $ d(x,\partial\widetilde{\mathbb F}(t,w)) = \inf\limits_{y \in \partial\widetilde{\mathbb F}(t,w)} d(x,y).$
				By the covering definition,
				$$\chi(\partial \widetilde{\mathbb{F}}(t,K)) \leq C\delta(t)r < C\delta(t)(\chi(K)+\varepsilon).$$
				Letting $\varepsilon \to 0$ gives
				$$\chi(\partial \widetilde{\mathbb{F}}(t,K)) \leq C\delta(t)\chi(K) \text{ a.e. } t \ge 0 .$$
			\end{itemize}
			Hence, assumption \textit{(F)$^{\prime}$} holds.
			
			\textbf{Abstract Formulation:} Now, we reformulate the problem \eqref{5.10}-\eqref{5.20} by rewriting  as an abstract problem in the state space $\mathbb{V} = \mathrm{L}^p([0,\pi],\mathbb{R})$ with $ p \in [2, \infty) $, which is a reflexive Banach space and the control space $ \mathbb{U} = \mathrm{L}^2([0,\pi],\mathbb{R})$. Note that
			$\mathbb{V}^* = \mathrm{L}^{p'}([0,\pi],\mathbb{R}),$
			which is uniformly convex. To this aim, let us define 
			\begin{align*}
				q(t)(\xi) = w(t,\xi),~ t \in [0,a],\  \xi \in [0,\pi].
			\end{align*}
			Thus, the abstract reformulation of equation \eqref{5.10}-\eqref{5.20} may be given as the follows:
			\begin{equation}\label{5.141}
				\left\{
				\begin{aligned}
					q'(t) &\in \mathscr{A}q(t)+ \int_0^t G(t-s)q(s)ds + Bu(t)+H \partial \mathbb{F}(t, q(t)),~ t\in (0,a],\\
					q(0) &\in h(q(\cdot)),
				\end{aligned}
				\right.
			\end{equation}
			where, $h : \mathrm{C}([0,a],\mathbb{V}) \to \mathcal{P}(\mathbb{V})$ is given as follows:
			\begin{align}\label{def-h}
				h(q):=  \left\{ z \in W^{1,p}([0,\pi],\mathbb{R}) : |z(\xi)| \le \|q\|_{\infty}^{1/2},~ |z'(\xi)| \le \|q\|_{\infty}^{1/2} \text{ a.e. } \xi \in [0,\pi]\right\},
			\end{align}
			and the operators $\mathscr{A}: \mathscr{D}(\mathscr{A}) \subset \mathbb{V} \to \mathbb{V}$ defined as
			$$
			\mathscr{A} g= \frac{\partial}{\partial \xi} g \ \text{ with } \ \mathscr{D}(\mathscr{A})=\{g \in \mathrm{H}^1([0,\pi],\mathbb{R}): g(\pi)=0\},
			$$
			with
			$$ N(t) g =0 \ \text{ for }\  g \in \mathbb{V} \ \text{ and } \
			G(t) g =  e^{-\beta t} g\ \text{ for }\ g \in  \mathbb{V}.$$
			In addition, using the notation $u(t)(\xi)=u(t,\xi)$, and define the operator $B : \mathbb{U} \to \mathbb{V}$ by
			\begin{align*}
				Bu(t)(\xi) := \widetilde{\eta}(t,\xi) = \int_0^\pi \mathscr{K}(\xi,\omega)\,u(t)(\omega)\,d\omega,
				\quad t \in [0,a],~ \xi \in [0,\pi],
			\end{align*}
			where the kernel $\mathscr{K} \in \mathrm{C}([0,\pi]\times[0,\pi],\mathbb{R})$ is symmetric, that is, $\mathscr{K}(\theta,\xi)=\mathscr{K}(\xi,\theta)$ for all $\theta,\xi \in [0,\pi]$. For instance,
			\begin{align*}
				\mathscr{K}(\theta,\xi)=
				\begin{cases}
					\xi(\pi-\theta), & 0 \le \xi \le \theta \le \pi, \\[4pt]
					(\pi-\xi)\theta, & 0 \le \theta \le \xi \le \pi,
				\end{cases}
				\quad \text{or} \quad
				\mathscr{K}(\theta,\xi)=\min\{\xi,\theta\}.
			\end{align*}
			Then $B$ is a bounded operator (see application section of \cite{SAMTDJ}) and satisfies
			\begin{align}\label{bound-B}
				\|B\|_{\mathcal{L}(\mathbb{U},\mathbb{V})} \le C_k \pi^{\frac{1}{2}+\frac{1}{p}},
			\end{align}
			where $C_k=\sup\limits_{(\theta,\xi)\in[0,\pi]^2}\mathscr{K}(\theta,\xi)$. Similarly, define $H : \mathbb{V}^* \to \mathbb{V}$ by
			\[
			H(x^*)(\xi)=\int_0^\pi \mathscr{L}(\xi,\omega)\,x^*(\omega)\,d\omega,
			\quad x^* \in \mathbb{V}^*,~\xi \in [0,\pi],
			\]
			for which
			\begin{align}\label{bound-H}
				\|H\|_{\mathcal{L}(\mathbb{V}^*,\mathbb{V})} \le C_l\,\pi^{\frac{1}{p}},
			\end{align}
			where $C_l=\sup\limits_{(\theta,\xi)\in[0,\pi]^2}\mathscr{L}(\theta,\xi)$. 
		
		One can easily verify that $\mathscr{D}(\mathscr{A})$ is dense in $\mathbb{V}$ and the operator $\mathscr{A}$ with this domain is closed. Moreover, the operator $\mathscr{A}$ is the infinitesimal generator of a contraction semigroup $\mathcal{S}(t)$ on $\mathbb{V}$.
		The semigroup $(\mathcal{S}(t))_{t \geq 0}$ (see \cite[Example 2.3.7]{MTGW}) is given by
		\begin{align*}
			(\mathcal{S}(t)g)(\xi)=g(\xi+t), \quad \xi \in [0,\pi],\ t\in [0,a].
		\end{align*}
		Thus, 
		$$\mathcal{F}(\lambda)=\left(\lambda \mathtt{I}-\mathscr{A}-\frac{\mathtt{I}}{\lambda+\beta}\right)^{-1}.$$
		By using the definition of resolvent family defined as \eqref{2.31}, the Laplace  transform of $\mathscr{R}(t)$ is
		\begin{align}\label{4.53}
			\widehat{\mathscr{R}}(\lambda) = \left(\lambda \mathtt{I}-\mathscr{A}-\frac{\mathtt{I}}{\lambda+\beta} \right)^{-1}=\left(\lambda \mathtt{I}-\mathscr{A}\right)^{-1}\left(\mathtt{I}-\frac{\left(\lambda \mathtt{I}-\mathscr{A}\right)^{-1}}{\lambda+\beta} \right)^{-1}.
		\end{align}
		It is easy to see that there exists a large number $r$ such that $|\lambda|>r$ such that
		$$\bigg\|\frac{\left(\lambda \mathtt{I}-\mathscr{A}\right)^{-1}}{\lambda+\beta}\bigg\|_{\mathcal{L}(\mathbb{V})}\leq \frac{1}{|\lambda||\lambda+\beta|}<1.$$
		Thus, we can write \eqref{4.53} as 
		\begin{align*}
			\widehat{\mathscr{R}}(\lambda) = \left(\lambda \mathtt{I}-\mathscr{A}\right)^{-1} \sum_{n=0}^{\infty}\left(\frac{\left(\lambda \mathtt{I}-\mathscr{A}\right)^{-1}}{\lambda+\beta} \right)^{n}.
		\end{align*}
		Hence, we obtain
		\begin{align}\label{4.55}
			\widehat{\mathscr{R}}(\lambda) = \sum_{n=0}^{\infty}\left(\lambda \mathtt{I}-\mathscr{A}\right)^{-(n+1)}(\lambda+\beta)^{-n}.
		\end{align}
		Since, by using the relation
		$$
		(\lambda \mathtt{I} - \mathscr{A})^{-1}g
		= \int_{0}^{\infty} e^{-\lambda s} \mathcal{S}(s)gds,
		$$
		which implies 
		\begin{align*}
			\mathfrak{L}(\mathcal{S}(t))= (\lambda \mathtt{I} - \mathscr{A})^{-1},
		\end{align*} 
		where $\mathfrak{L}$ represents the Laplace transform. By using the fundamental properties of the Laplace transform, we get 
		$$\mathfrak{L}(t^n \mathcal{S}(t))= n! (\lambda \mathtt{I} - \mathscr{A})^{-(n+1)},$$ 
		and $$ \mathfrak{L}(t^{n-1} e^{-\beta t})= (n-1)! (\lambda +\beta)^{-n}.$$ 
		Thus, by taking the inverse Laplace transform to \eqref{4.55}, we obtain
		\begin{align*}
			\mathscr{R}(t)g = \mathcal{S}(t)g+\sum_{n=1}^{\infty} \int_{0}^{t}\frac{(t-s)^{n-1}}{(n-1)!} e^{-\beta(t-s)} \frac{s^n}{n !}\mathcal{S}(s)gds.
		\end{align*}
		Since, the semigroup $(\mathcal{S}(t))_{t>0}$ is not compact (see \cite[Example 4.31]{KJENR}), which in turn implies that the resolvent family $(\mathscr{R}(t))_{t>0}$ is not compact.
		
		Now, we verify that the multivalued mapping $h(\cdot)$ satisfies all the hypotheses given in Theorem \ref{Theorem-3}.
		\begin{itemize}
			\item[\textit{(H1)}] We now show that the multimap $h$ maps bounded set of $\mathrm{C}([0,a],\mathbb{V})$ into bounded sets in $\mathbb{V}.$ For this, let 
			$$Q_{\mathfrak{R}} = \{q \in \mathrm{C}([0,a],\mathbb{V}) : \|q\|_{\mathrm{C}([0,a], \mathbb{V})}\leq \mathfrak{R}\},~\mathfrak{R} > 0,$$ be a bounded set in $\mathrm{C}([0,a],\mathbb{V})$ and $z \in h(q)$. Then, we have
			$$|z(\xi)| \leq \|q\|_{\infty}^{1/2},~|z'(\xi)| \leq \|q\|_{\infty}^{1/2} \text{ a.e. } \xi \in [0,\pi].$$
			Hence,
			$$ \|z\|_{\mathbb{V}} \leq \pi^{1/p} \|q\|_{\infty}^{1/2} \leq \pi^{1/p} \mathfrak{R}^{1/2}.$$
			Thus,
			$$ \|h(Q_{\mathfrak{R}})\|_{\mathbb{V}} = \sup\{\|z\|_{\mathbb{V}} : z \in h(q),~ q \in Q_{\mathfrak{R}}\}
			\leq \pi^{1/p} \mathfrak{R}^{1/2},$$ 
			showing that $h$ indeed maps bounded sets into bounded sets.
			
			\item[\textit{(H3)}] Next, we prove that the multimap $h$ has closed graph. Assume that $q_n \in \mathrm{C}([0,a],\mathbb{V})$, with $z_n \in h(q_n)$ satisfy
			$$q_n \to q \text{ in } \mathrm{C}([0,a],\mathbb{V}), ~ z_n \to z \text{ in } \mathbb{V}.$$ 
			From the definition of the multimap $h$, each $z_n \in W^{1,p}([0,\pi],\mathbb{R})$ and 
			\begin{align}\label{5.171}
				|z_n(\xi)| \le \|q_n\|_{\infty}^{1/2},~ |z^{\prime}_n(\xi)| \le \|q_n\|_{\infty}^{1/2} \ \text{ a.e. }\ \xi \in [0,\pi]. 
			\end{align}
			Since, $q_n \to q \text{ in } \mathrm{C}([0,a],\mathbb{V}),$ thus the sequence $\|q_n\|_{\infty}$ is bounded; hence
			\begin{align}\label{5.181}
				|z_n(\xi)| \le C,~ |z^{\prime}_n(\xi)| \le C \ \text{ a.e. }\ \xi \in [0,\pi].
			\end{align}
			By using \eqref{5.181}, we get 
			$$\|z_n\|_{W^{1,p}([0,\pi],\mathbb{R})}\leq 2C \pi^{1/p},$$
			thus by Banach-Alaoglu theorem, there exists a subsequence $\{z_n\}_{n \in \mathbb{N}}$ (denoted by same symbol) such that 
			\begin{align}\label{5.191}
				z_n \xrightarrow{w} z \ \text{ in }\ W^{1,p}([0,\pi],\mathbb{R}).
			\end{align}
			Since, by the Rellich–Kondrachov theorem, 
			$$ W^{1,p}([0,\pi],\mathbb{R}) \underbrace{\hookrightarrow}_{\text{compact}} \mathrm{C}([0,\pi], \mathbb{R}),$$ 
			we obtain $$z_n(\xi) \to z(\xi) \ \text{ for all }\ \xi \in [0,\pi].$$
			Passing to the limit $n \to \infty$ in \eqref{5.171} gives
			$$|z(\xi)| \le \|q\|_{\infty}^{1/2} \ \text{ for all }\ \xi \in [0,\pi].$$
			Furthermore, from \eqref{5.191}, we have
			$z'_n \xrightarrow{w} z'$ in $\mathrm{L}^{p}([0,\pi],\mathbb{R})$ up to a subsequence. Therefore, an application of Mazur's Lemma \cite{IERT} yields a sequence
			\begin{equation*}
				\tilde{\phi}_{n} = \sum_{k=n}^{N} \Lambda_{k} z'_k,\text{ where } \sum_{k=n}^{N} \Lambda_{k} = 1 \text{ and }\Lambda_{k} \geq 0,~ n \leq k \leq N,
			\end{equation*}
			such that $\tilde{\phi}_{n} \to z'$ as $n \to \infty$ and 
			$\tilde{\phi}_{n}(\xi) \to z'(\xi)$ a.e. $\xi \in [0,\pi]$ up to a subsequence. 
			Hence, we obtain
			\begin{equation*}
				|\tilde{\phi}_{n}(\xi)| \leq \sum_{k=n}^{N} 
				\Lambda_{k}\|q_{k}\|_{\infty}^{1/2} \ \text{ for a.e. } \ \xi \in [0,\pi].
			\end{equation*}
			Letting $n \to \infty$, we obtain $$|z'(\xi)| \leq \|q\|_{\infty}^{1/2}\ \text{ for a.e. } \ \xi \in [0,\pi].$$
			Therefore, $z \in h(q)$ proves that the graph of $h$ is closed, which in turn implies that $h$ is closed. 
			
			
			Now, consider a bounded sequence 
			$\{q_n\}_{n \in \mathbb{N}}$ in $ \mathrm{C}([0,a],\mathbb{V})$. We prove that the sequence 
			$\{h(q_n)\}_{n \in \mathbb{N}}$ admits a convergent subsequence in $\mathbb{V}$. 
			Let $z_n \in h(q_n)$ for each $n \in \mathbb{N}$. From the definition of the multimap $h$, we have
			\begin{equation*}
				|z_n(\xi)| \le \|q_n\|_{\infty}^{1/2},~ 
				|z_n'(\xi)| \le \|q_n\|_{\infty}^{1/2}\ \text{ for all } n \in \mathbb{N}, \text{ a.e. } \xi \in [0,\pi].
			\end{equation*}
			Hence, $\{z_n\}_{n \in \mathbb{N}}$ 
			is a bounded sequence in $W^{1,p}([0,\pi],\mathbb{R})$, and by the Rellich–Kondrachov compact embedding theorem, we have 
			$$ W^{1,p}([0,\pi],\mathbb{R}) \underbrace{\hookrightarrow}_{\text{compact}} \mathrm{C}([0,\pi], \mathbb{R})\underbrace{\hookrightarrow}_{\text{continuous}} \mathrm{L}^{p}([0,\pi],\mathbb{R}).$$
			Therefore, $\{z_n\}_{n \in \mathbb{N}}$ is relatively 
			compact in $\mathrm{L}^{p}([0,\pi],\mathbb{R})$, and thus (upto a subsequence) we obtain
			$$z_n \to z \text{ in } \mathrm{L}^{p}([0,\pi],\mathbb{R}).$$
			This shows that the multivalued map $h$ is compact. Thus, combining the closedness of $h$ with the compact multivalued map, it implies $h$ has compact values, and invoking \cite[Theorem 1.1.5]{kamenskii}, we conclude the multimap $h$ is u.s.c.
			\item[\textit{(H5)}] Since, $h(W)$ is relatively compact for all bounded sets $W \subseteq \mathrm{C}([0,a], \mathbb{V})$, thus it satisfies the assumption \textit{(H5)}. 
		\end{itemize}		
		Since, by using \eqref{bound-B}, \eqref{bound-H}, along with
		$$\Phi_{\mathfrak{R}}=\pi^{1/p} \mathfrak{R}^{1/2}, ~\gamma(t) =\max\{2\widetilde{\mu}_0,\widetilde{c}\} \delta(t)(\pi)^{\frac{1}{p'}},~ q_0=0, ~ k_h=0 \ \text{ and }\ k(t)=C\delta(t),$$ we deduce that, for sufficiently large $\mathfrak{R}>0$, one can always choose a function $\delta\in \mathrm{L}^1([0,a],\mathbb{R}^+)$ such that
		$$ \int_{0}^{a} \delta(s) ds <\frac{1}{\max\{2\widetilde{\mu}_0,\widetilde{c}\} C_l\pi}\left(\frac{\mathfrak{R}}{M}-\pi^{1/p} \mathfrak{R}^{1/2}-C_k \pi^{\frac{1}{2}+\frac{1}{p}} a^{1/2}\norm{u}_{\mathrm{L}^2(I, \mathbb{U})}\right)$$
		and 
		$$ \int_{0}^{a} \delta(s)ds < \frac{1}{2M \widetilde{C}} ,$$ 
		where $\widetilde{C}=C C_l \pi^{1/p}$. Hence, the conditions \eqref{assum-1} and \eqref{3.38} are satisfied. Therefore, all hypotheses of Theorem \ref{Theorem-3} are satisfied, and the system \eqref{5.141} has a mild solution.
	\end{ex}
	Now, we investigate the approximate controllability of a heat equation that demonstrates the applicability of our developed results. 
	\begin{ex}\label{5.1111}
		Consider the following class of equations:
		\begin{align}\label{5.11}
			\dfrac{\partial w(t,\xi)}{\partial t}= \dfrac{\partial^{2}w(t,\xi)}{\partial \xi^{2}}+\eta(t,\xi)+\int_0^{\pi} \mathscr{L}(\xi,\omega) f(t,\omega)d \omega \ \text{ in }\ (0,a] \times [0,\pi],
		\end{align}
		where
		\begin{align}\label{5.21}
			f(t,\xi) \in \partial \mathbb{F}(t,w(t,\xi))\ \text{ for a.e. }\ t \in (0,a],~ \xi \in [0,\pi].
		\end{align}
		The system is supplemented with the boundary and non-local initial conditions
		\begin{align*}
			&w(t,0)=w(t,\pi)=0, ~ t \in [0,a], \\
			&w(0,\xi)\in h(w(\cdot,\xi)),~ \xi \in [0,\pi].
		\end{align*}
		Thus, combining \eqref{5.11} and \eqref{5.21},  we obtain the following differential inclusion:
		\begin{align}\label{5.51}
			\dfrac{\partial w(t,\xi)}{\partial t} \in \dfrac{\partial^{2}w(t,\xi)}{\partial \xi^{2}}+\eta(t,\xi)+ \int_0^{\pi} \mathscr{L}(\xi,\omega) \partial \mathbb{F}(t,w(t,\omega))d \omega\ \text{ in }\ (0,a] \times [0,\pi].
		\end{align}      
		Now, we define the functional $\mathbb{F}:[0,a] \times \mathrm{L}^p([0,\pi],\mathbb{R})\to \mathbb{R}$ as follows:
		$$\mathbb{F}(t,w):= \int_0^{\pi} j\big(t,w(\tau)\big)d\tau, $$
		where $j: [0,a] \times \mathbb{R} \to \mathbb{R}$ is defined as 
		$$j\big(t,w(\tau)\big)=\delta(t)g\big(w(\tau)\big).$$ 
		In the above expression, $\delta : [0,a] \to \mathbb{R} $ is a non-negative measurable function with 
		$ \delta \in \mathrm{L}^1([0,a],\mathbb{R}^+),$ and $g : \mathbb{R} \to \mathbb{R}$ is defined as 
		$$g(r) = |r| + \sin r,$$
		which is locally Lipschitz. Let us now verify the conditions \textit{(F1)}-\textit{(F2)} satisfied by the functional $\mathbb{F}$ written as follows:
		\begin{itemize}
			\item[\textit{(F1)}]  Since $\delta : [0,a] \to \mathbb{R} $ is a measurable function, it follows that $$\mathbb{F}(\cdot,w) \text{ is measurable on } [0,a] \text{ for all } w \in \mathrm{L}^p([0,\pi],\mathbb{R}).$$ Moreover, $\mathbb{F}(t,\cdot)$ is locally Lipschitz on bounded subsets of $ \mathrm{L}^p([0,\pi],\mathbb{R})$ for \mbox{a.e.} $t \in [0,a]$ which follows from the locally Lipschitz property of the function $g$
			together with H{\"o}lder's inequality, as demonstrated below:
			\begin{align*}
				|\mathbb{F}(t,w_1) - \mathbb{F}(t,w_2)|&= \left|\int_0^{\pi} \delta(t)    \big(g(w_1(\tau)) - g(w_2(\tau))\big)d\tau \right|  \\
				&\leq |\delta(t)| \int_0^{\pi} |g(w_1(\tau)) - g(w_2(\tau))|d\tau   \\
				&\leq \delta(t) M \int_0^{\pi} |w_1(\tau) - w_2(\tau)| d \tau,  \\
				&\leq \delta(t) M (\pi)^{\frac{1}{p'}}\|w_1 - w_2\|_{ \mathrm{L}^p([0,\pi],\mathbb{R})},
			\end{align*}
			where $M > 0$ is a local Lipschitz constant for $g$. 
			\item[\textit{(F2)}] By applying \cite[Proposition 3.32]{SMAOMS}, we obtain the subdifferential of $g(r)$ which is given by
			\begin{equation*}
				\partial g(r)=
				\begin{cases}
					\begin{aligned}
						\cos r - 1,~ & r < 0,\\[2pt]
						[0,2],~ & r = 0,\\[2pt]
						\cos r + 1,~ & r > 0. 
					\end{aligned}
				\end{cases}
			\end{equation*}
			Next, using the definition of the generalized gradient of $j(t,\cdot)$ (see \cite[Proposition 3.44]{SMAOMS}), we obtain 
			\begin{align*}
				|\partial j(t,w(\tau))| =|\delta(t)\partial g(w(\tau))| \leq 2\delta(t) \ \mbox{ a.e.}\ t \in [0,a].
			\end{align*}
			Furthermore, by employing an argument analogous to the proof of Example \ref{ex 3.11}, we have
			\begin{align}\label{equality-211}
				\|\zeta\|_{\mathrm{L}^{p'}([0,\pi],\mathbb{R})}^{p'}
				= \int_{0}^{\pi} |\zeta(\tau)|^{p'}d \tau \leq\int_{0}^{\pi}  (2\delta(t))^{p'} d \tau
				\leq (2\delta(t))^{p'}\pi,
			\end{align}
			and for $z^* \in \partial \mathbb{F}(t,w)$, we have  equality \eqref{equality-1}. Thus, by comparing \eqref{equality-1} and \eqref{equality-211},
			we get $\|z^*\|_{\mathrm{L}^{p'}([0,\pi],\mathbb{R})} \leq 2\delta(t)(\pi)^{\frac{1}{p'}},$  which shows that condition \textit{(F2)} is satisfied and further \textit{(F2)}$^\prime$ holds.
		\end{itemize}
		\textbf{Abstract Formulation:}	Now, we reformulate the problem \eqref{5.51} by rewriting as an abstract problem in the state space $\mathbb{V} = \mathrm{L}^p([0,\pi],\mathbb{R})$ with $ p \in [2, \infty) $ with  
		$\mathbb{V}^* = \mathrm{L}^{p'}([0,\pi],\mathbb{R})$ and the control space $ \mathbb{U} = \mathrm{L}^2([0,\pi],\mathbb{R})$. To this aim, let us define 
		\begin{align*}
			q(t)(\xi) = w(t,\xi),~ t \in [0,a], \xi \in [0,\pi].
		\end{align*}		
		Now, we define the operators 
		$\mathscr{A}: \mathscr{D}(\mathscr{A}) \subset \mathbb{V} \to \mathbb{V}$
		as
		$$ \mathscr{A} g = g'', ~ \mathscr{D}(\mathscr{A}) = W^{2,p}([0,\pi],\mathbb{R}) \cap W^{1,p}_0([0,\pi],\mathbb{R}).$$
		Moreover, the operator $ \mathscr{A} $ can be written as
		$$ \mathscr{A} g = \sum_{k=1}^{\infty} -k^2 \widehat{g}_k v_k, ~ g \in  \mathscr{D}(\mathscr{A}),$$
		and $ (\mathscr{R}(t))_{t \geq 0}$ generated by the operator $\mathscr{A} $ is given by
		$$ \mathscr{R}(t)=\sum_{k=1}^{\infty} e^{-k^2t}\widehat{g}_k v_k,~ g \in \mathbb{V},$$
		where 
		\begin{align}\label{fourier-s}
			v_k(\xi) = \sqrt{\frac{2}{\pi}} \sin(k\xi) \ \text{ and }\
			\widehat{g}_k := \int_{0}^{\pi} g(\xi) v_k(\xi) d\xi.
		\end{align}
		Also, the operator $B : \mathbb{U} \to \mathbb{V}$ is defined as
		\begin{align*}
			Bu(t)(\xi) := \eta(t, \xi) =\sum_{k=1}^{\infty} \frac{1}{e^{k^l}} \langle u(t), v_k \rangle v_k ,~ t \in [0,a], \xi \in [0,\pi],
		\end{align*}
		which is a bounded operator. In order to show this, let us estimate
		\begin{align*}
			\|Bu(t)\|_{\mathbb{V}}
			= \left\| \sum_{k=1}^{\infty} \frac{1}{e^{k^l}} \langle u(t), v_k \rangle v_k \right\|_{\mathbb{V}} \leq \sum_{k=1}^{\infty} \frac{|\langle u(t), v_k \rangle|}{e^{k^l}} \|v_k\|_{\mathbb{V}} \leq \sqrt{2}\pi^{\frac{1}{p}-\frac{1}{2}}\sum_{k=1}^{\infty} \frac{|\langle u(t), v_k \rangle|}{e^{k^l}}.
		\end{align*}
		By using the Cauchy-Schwarz inequality and $\|v_k\|_{\mathrm{L}^{2}} = 1$ for all $k \in \mathbb{N}$, we obtain
		\begin{align*}
			\|Bu(t)\|_{\mathbb{V}}
			\leq \sqrt{2}\pi^{\frac{1}{p}-\frac{1}{2}} \|u(t)\|_{\mathbb{U} } \sum_{k=1}^{\infty} \frac{1}{e^{k^l}}.
		\end{align*}
		Since the series $\sum_{k=1}^{\infty} \frac{1}{e^{k^l}}$ is convergent for any $l>0$, thus we get $\|B\|_{\mathcal{L}(\mathbb{U} ,\mathbb{V})} \leq C_p$.
		Define $H : \mathbb{V}^* \to \mathbb{V}$ via
		$$H(x^*)(\xi) = \int_0^\pi \mathscr{L}(\xi,\omega)x^*(\omega)d\omega,~ x^* \in \mathbb{V}^*,~\xi \in [0,\pi].$$
		Now, the abstract reformulation of equation \eqref{5.51} can be given as the follows:
		\begin{equation}\label{5.1411}
			\left\{
			\begin{aligned}
				q'(t) &\in \mathscr{A} q(t)+ Bu(t) + H \partial \mathbb{F}(t, q(t)),~ t\in (0,a],\\
				q(0) &\in h(q(\cdot)),
			\end{aligned}
			\right.
		\end{equation}
		where, $h : \mathrm{C}([0,a],\mathbb{V}) \to \mathcal{P}(\mathbb{V})$ is a multivaled map defined in \eqref{def-h} that satisfies all hypotheses of Theorem \ref{MR}, as verified in Example~\ref{ex 3.11}. 
		Also, the adjoint operator 
		$$B^*: \mathbb{V}^*\to \mathbb{U} $$ is given by
		\begin{align*}
			B^* w^*(\xi) := \sum_{k=1}^{\infty} \frac{\widehat{w^*_k}}{e^{k^l}} v_k(\xi),
		\end{align*}
		where $\widehat{w^*_k}:= \int_{0}^{\pi} w^*(\xi) v_k(\xi) d\xi .$ 
		
		Now, we prove that Hypothesis (R) is satisfied. To prove this, let us estimate
		\begin{align*}
			R(t)g^*
			= \int_0^t \mathscr{R}(t-s) B B^* \mathscr{R}(a-s)^*g^* ds
			&= \sum_{k=1}^{\infty} \frac{1}{e^{2k^l}} \left( \int_0^t e^{-k^2(a+t-2s)} ds \right) \widehat{g^*_k} v_k \\
			&=\sum_{k=1}^{\infty} \frac{e^{-k^2(a-t)}-e^{-k^2(a+t)}}{2k^2 e^{2k^l}}  \widehat{g^*_k} v_k .
		\end{align*}
		We now first verify that the corresponding linear system associated with problem \eqref{5.1411} is approximately controllable. To accomplish this, consider
		$$B^* \mathscr{R}(a-t)^*y^* =  0 \ \text{ for any }\ y^* \in \mathbb{V}^* ,$$
		which further implies $$\sum_{k=1}^{\infty} \frac{e^{-k^2(a-t)}\widehat{y^*_k}}{e^{2k^l}} v_k= 0\ \text{ for any } \ y^* \in \mathbb{V}^* .$$
		Taking the duality pairing with $v_m $, it follows that
		$$ \frac{e^{-k^2(a-t)}\widehat{y^*_m}}{e^{2m^l}} = 0\ \text{ for all }\ m \in \mathbb{N},$$
		which further implies 
		$$\widehat{y^*_m}=0 \ \text{ for all }\ m \in \mathbb{N}.$$ 
		Hence $y^*=0 \text{ in }  \mathbb{V}^* .$ Therefore, the linear system corresponding to \eqref{5.1411} is approximately controllable (see Lemma \ref{lemma 5.1}). Consequently, by \cite[Remark 4.2]{MFPH}, the operator $R(a)$ is injective. Hence, its inverse 
		$$ R(a)^{-1} :\operatorname{Im}R(a) \to \mathbb{V}^* $$
		is well-defined as a linear operator. 
		
		Let $y \in \operatorname{Im}R(a)$. Then there exists $x \in \mathbb{V}$ such that $y = R(a)x$. Thus,
		\begin{equation*}
			y = R(a)x =  \sum_{k=1}^{\infty} \frac{1-e^{-2k^2a}}{2k^2 e^{2k^l}} \widehat{x}_k v_k= \sum_{k=1}^{\infty} \widehat{y}_k v_k,
		\end{equation*}
		where
		\begin{equation*}
			\widehat{y}_k  = \frac{1-e^{-2k^2a}}{2k^2 e^{2k^l}} \widehat{x}_k.
		\end{equation*}
		Therefore, the inverse operator $R(a)^{-1}$ is given by
		\begin{equation*}
			R(a)^{-1} y = \sum_{k=1}^{\infty} \frac{2k^2 e^{2k^l}}{1-e^{-2k^2a}} \widehat{y}_k v_k.
		\end{equation*}
		Now, we show that the operator
		$$
		R(t)R(a)^{-1} : \operatorname{Im}R(a) \to \mathbb{V} $$
		is bounded. Let $y \in \operatorname{Im}R(a)$. Then, we have 
		\begin{align}\label{R-bound}
			R(t)R(a)^{-1} y 
			&= R(t)\sum_{k=1}^{\infty} \frac{2k^2 e^{2k^l}}{1-e^{-2k^2a}} \widehat{y}_k v_k = \sum_{k=1}^{\infty}  
			\frac{e^{-k^2(a-t)}-e^{-k^2(a+t)}}
			{1-e^{-2k^2a}} \, \widehat{y}_k v_k = \sum_{k=1}^{\infty} \rho_k(t,a)\widehat{y}_k v_k,
		\end{align}
		where
		\begin{align}\label{5.961}
			\rho_k(t,a) := e^{-k^2(a-t)}\cdot \frac{1-e^{-2k^2t}}{1-e^{-2k^2a}},
		\end{align}
		or equivalently, we may write 
		$$\rho_k= e^{-k^2(a-t)}
		\cdot
		\frac{2e^{-k^2t}\sinh(k^2 t)}{2e^{-k^2 a}\sinh(k^2 a)}=
		\frac{\sinh(k^2 t)}{\sinh(k^2 a)}.$$
		Since \(0\le t\le a\) and \(\sinh\) is positive and increasing on \((0,\infty)\), therefore, we have 
		$$ 0\le \sinh(k^2 t)\le \sinh(k^2 a),$$
		which implies
		$$ 0\le \rho_k = \frac{\sinh(k^2 t)}{\sinh(k^2 a)}\leq 1.$$
		Hence, we get 
		\begin{align}\label{M1}
			\sup_{k\ge1} |\rho_k|\le 1.
		\end{align}
		Now, it is easy to verify that $ \rho_k $ is decreasing. Moreover, since $ 0\le  \rho_k \le 1$, we obtain
		\begin{align}\label{M2}
			\sup_{j\geq 0}\sum_{k=2^j}^{2^{j+1}}|\rho_k-\rho_{k+1}|
			=\sup_{j\geq 0}\sum_{k=2^j}^{2^{j+1}}(\rho_k-\rho_{k+1})=\sup_{j\geq 0}(\rho_{2^j}-\rho_{{2^{j+1}+1}})
			\le 1.
		\end{align}
		Define $T:\mathrm{L}^p([0,2\pi],\mathbb{R})\to\mathrm{L}^p([0,2\pi],\mathbb{R})$ as
		$$(Tg)(\xi)=\sum_{k=1}^{\infty}\rho_k(t,a)
		\bigl[a_k\cos(k\xi)+b_k\sin(k\xi)\bigr],$$ where $a_k$ and $b_k$ are the Fourier coefficients of $g$. By \eqref{M1}-\eqref{M2} and Marcinkiewicz's multiplier theorem \cite[Theorem 1]{JMAR}, the operator $T$  is bounded on $\mathrm{L}^p([0,2\pi],\mathbb{R})$ and 
		$$
		\|T g\|_{\mathrm{L}^p([0,2\pi],\mathbb{R})}
		\leq M_p\|g\|_{\mathrm{L}^p([0,2\pi],\mathbb{R})} , \ g\in\mathrm{L}^p([0,2\pi],\mathbb{R}) 
		$$ where $M_p$ is a positive constant.
		
		Further, for any  $g=\mathrm{L}^p([0,\pi],\mathbb{R})$, we define its odd $2\pi$-periodic extension $\tilde{g}$ on $[0,2\pi]$ by
		$$
		\tilde g(\xi)=
		\begin{cases}
			g(\xi), & \xi\in [0,\pi),\\[4pt]
			-g(2\pi-\xi), & \xi\in(\pi,2\pi].
		\end{cases}
		$$
		As this extension is odd, its constant and cosine coefficients vanish. Also, it is immediate that
		\begin{align}\label{iso}
			\|\widetilde{g}\|_{\mathrm{L}^p([0,2\pi],\mathbb{R})}= 2^{1/p}\|g\|_{\mathrm{L}^p([0,\pi],\mathbb{R})}.
		\end{align}
		Its Fourier series therefore consists only of sine terms, that is, 
		$$\tilde{g}(\xi)=\sum_{k=1}^{\infty}b_k\sin(k\xi),$$ where $b_k=\frac{2}{\pi}\int_0^\pi g(x)sin(k\xi)$. Consequently, we find 
		\begin{align*}\tilde{g}(\xi)=\sum_{n=1}^{\infty} \hat{g}_kv_k(\xi), \ \text{ where }\ 
			v_k(\xi) = \sqrt{\frac{2}{\pi}} \sin(k\xi)\  \text{ and }\ 
			\widehat{g}_k = \int_{0}^{\pi} g(\xi) v_k(\xi) d\xi.
		\end{align*}
		Subsequently, we have 
		$$T\tilde{g}=\sum_{k=1}^{\infty}\rho(k,a)\hat{g}_kv_k.$$
		Further, define $S:\mathrm{L}^p([0,\pi],\mathbb{R})\to\mathrm{L}^p([0,\pi],\mathbb{R})$ as
		$$
		S g:=\left. T\tilde{g}\right|_{[0,\pi]}, \ g\in\mathrm{L}^p([0,\pi],\mathbb{R}), $$
		where $\tilde{g}$ is the odd $2\pi$-periodic extension on $[0,2\pi]$ of $g$. Since $T$ preserves oddness, $T\widetilde{g}$ is the odd $2\pi$-periodic extension of $Sg$. Thus,
		\begin{align*}
			(T\widetilde{g})(\xi)=
			\begin{cases}
				(Sg)(\xi), & \xi\in [0,\pi),\\
				-(Sg)(2\pi-\xi), & \xi\in(\pi,2\pi].
			\end{cases}
		\end{align*}
		Consequently,
		\begin{align}\label{iso2}
			\|T\widetilde{g}\|_{\mathrm{L}^p([0,2\pi],\mathbb{R})}^p
			&=
			\int_0^\pi |(Sg)(\xi)|^p\,d\xi
			+
			\int_\pi^{2\pi}|(Sg)(2\pi-\xi)|^p\,d\xi\nonumber\\
			&=2\|Sg\|_{\mathrm{L}^p([0,\pi],\mathbb{R})}^p.
		\end{align}
		Using the estimates \eqref{iso}-\eqref{iso2} along with the boundedness of $T$, we obtain
		\[
		\begin{aligned}
			\|Sg\|_{\mathrm{L}^p([0,\pi],\mathbb{R})}
			&=
			2^{-1/p}
			\|T\widetilde{g}\|_{\mathrm{L}^p([0,2\pi],\mathbb{R})}\\
			&\leq
			2^{-1/p}M_p
			\|\widetilde{g}\|_{\mathrm{L}^p([0,2\pi],\mathbb{R})}\\
			&=
			M_p\|g\|_{\mathrm{L}^p([0,\pi],\mathbb{R})}.
		\end{aligned}
		\]
		Hence, $S$ is a bounded linear operator on
		$\mathrm{L}^p([0,\pi],\mathbb{R})$ with
		$\|S\|_{\mathcal{L}(\mathrm{L}^p([0,\pi],\mathbb{R})}\leq M_p$. Observe that 
		$$
		Sy=R(t)R(a)^{-1}y,
		\qquad y\in\operatorname{Im}R(a)\subset \mathbb{V}=\mathrm{L}^p([0,\pi],\mathbb{R}).
		$$
		Therefore,
		$$
		\|R(t)R(a)^{-1}y\|_{\mathbb{V}}
		\leq M_p\|y\|_{\mathbb{V}},
		\qquad y\in\operatorname{Im}R(a),
		$$
		which implies
		$$
		\|R(t)R(a)^{-1}\|
		_{\mathcal{L}(\operatorname{Im}R(a),\mathbb{V})}
		\leq M_p.
		$$
		which proves Hypothesis (R). Finally, by invoking Theorem \ref{MR}, we conclude that the control system \eqref{5.1411} is approximately controllable.
	\end{ex}
	
	Now, we investigate the approximate controllability of a heat equation involving a weakly singular kernel that demonstrates the applicability of our developed results. 
	\begin{ex}\label{5.111}
		Consider the following class of equations:
		\begin{align}\label{5.1}
			&\dfrac{\partial}{\partial t}\Bigg[w(t,\xi)+ \displaystyle\int_{0}^{t} (t-s)^{\gamma} e^{-\kappa(t-s)}w(s,\xi) ds\Bigg] \nonumber\\
			&= \dfrac{\partial^{2}w(t,\xi)}{\partial \xi^{2}}
			+ \dfrac{\alpha}{\Gamma(\nu)} \displaystyle\int_{0}^{t} e^{-\beta (t-s)} (t-s)^{\nu-1}
			w(s,\xi) ds+
			\eta(t,\xi)+ f(t,\omega)\ \text{ in }\ (0,a] \times [0,\pi],
		\end{align}
		where
		\begin{align}\label{5.2}
			f(t,\xi) \in \partial \mathbb{F}(t,w(t,\xi)) \text{ for a.e. } t \in (0,a],~ \xi \in [0,\pi].
		\end{align}
		The system is supplemented with the boundary and non-local initial conditions
		\begin{align*}
			&w(t,0)=w(t,\pi)=0, ~ t \in [0,a], \\
			&w(0,\xi)\in h((w(\cdot,\xi)),~ \xi \in [0,\pi].
		\end{align*}
		In the equation \eqref{5.1}, $\alpha$, $\kappa,\beta,\gamma$ are positive constants and $0<\nu<1$. Thus, combining \eqref{5.1} and \eqref{5.2},  we obtain the following differential inclusion:
		\begin{align}\label{5.5}
			&\dfrac{\partial}{\partial t}\Bigg[w(t,\xi)+ \displaystyle\int_{0}^{t} (t-s)^{\gamma} e^{-\kappa(t-s)}w(s,\xi) ds\Bigg] \nonumber\\
			&\in \dfrac{\partial^{2}w(t,\xi)}{\partial \xi^{2}}
			+ \dfrac{\alpha}{\Gamma(\nu)} \displaystyle\int_{0}^{t} e^{-\beta (t-s)} (t-s)^{\nu-1}
			w(s,\xi)ds+
			\eta(t,\xi)+ \partial \mathbb{F}(t,w(t,\omega)) \ \text{ in }\ (0,a] \times [0,\pi].
		\end{align}      
		Now, we define the functional $\mathbb{F}:[0,a] \times \mathrm{L}^2([0,\pi],\mathbb{R})\to \mathbb{R}$ as follows:
		$$\mathbb{F}(t,w):= \int_0^{\pi} j\big(t,w(\tau)\big)d\tau, $$
		where $j: [0,a] \times \mathbb{R} \to \mathbb{R}$ is defined as 
		$$j\big(t,w(\tau)\big)=\delta(t)g\big(w(\tau)\big).$$ 
		In the above expression, $\delta : [0,a] \to \mathbb{R} $ is a non-negative measurable function with 
		$ \delta \in \mathrm{L}^1([0,a],\mathbb{R}^+),$ and $g : \mathbb{R} \to \mathbb{R}$ is defined as 
		$$g(r) = |r| + \sin r,$$
		which is locally Lipschitz. From the Example \ref{5.1111}, \textit{(F1)}-\textit{(F2)} are satisfied by the functional $\mathbb{F}$ for $p=2$. 
		
		\textbf{Abstract Formulation:}	Now, we reformulate the problem \eqref{5.5}
		in $\mathbb{V} = \mathrm{L}^2([0,\pi],\mathbb{R})$ with
		$\mathbb{V}^* = \mathrm{L}^{2}([0,\pi],\mathbb{R})$ and the control space $ \mathbb{U} = \mathrm{L}^2([0,\pi],\mathbb{R})$. To this aim, let us define 
		\begin{align*}
			q(t)(\xi) = w(t,\xi),~ t \in [0,a], \xi \in [0,\pi].
		\end{align*}		
		Now, we define the operators $\mathscr{A} :  \mathscr{D}(\mathscr{A}) \subset \mathbb{V} \to \mathbb{V}$
		as
		$$ \mathscr{A} g = g'', ~  \mathscr{D}(\mathscr{A}) = H^2([0,\pi],\mathbb{R}) \cap H^1_0([0,\pi],\mathbb{R}),$$ and 
		\begin{align*}
			N(t) g = t^{\gamma} e^{-\kappa t} g \ \text{ for } \ g \in \mathbb{V} \text{ and }
			G(t) g = t^{\nu-1} e^{-\beta t} g \ \text{ for }\ g \in  \mathbb{V}.
		\end{align*}
		Also, the operator $B : \mathbb{U} \to \mathbb{V}$ is defined as
		\begin{align*}
			Bu(t)(\xi) := \eta(t, \xi) =u(t)(\xi) ,~ t \in [0,a], \xi \in [0,\pi],
		\end{align*}
		which is a bounded.
		Define $H : \mathbb{V}^* \to \mathbb{V}$ via
		$$H(x^*)(\xi) = x^*(\xi),~ x^* \in \mathbb{V}^*,~\xi \in [0,\pi].$$
		By \cite[Example 5.1]{arora2025controllability}, the assumptionns \textit{(A1)}-\textit{(A2)} are satsified. Moreover, it is easy to see that the conditions \textit{(A3)} is satisfied with $b(t)=t^{\nu-1} e^{-\beta t} $ such that 
		$$\hat{b}(\lambda)=\frac{\Gamma(\nu)}{(\lambda+\beta)^{\nu}} \text{ and } \norm{G(t) g}_{\mathbb{V}}\leq t^{\nu-1} e^{-\beta t} \norm{g}_1 \ \text{ for all }\ g \in  \mathscr{D}(\mathscr{A}),$$
		and the Laplace transformation of the
		operator $G(\cdot)$ is given as
		$$ \widehat{G}(\lambda)= \frac{\Gamma(\nu)}{(\lambda+\beta)^{\nu}} \mathtt{I}. $$
		In view of above expression, $\widehat{G}$ has an analytical extension to $\Lambda_\vartheta$ such that $$\|\widehat{G}(\lambda)\|_{\mathcal{L}(\mathbb{V})} \to 0 \ \text{ as }\ |\lambda| \to \infty.$$
		Finally, to verify condition \textit{(A4)}, we choose  $\widetilde{D} = \mathrm{C}_{0}^{\infty}([0,\pi],\mathbb{R})$ (space of infinitely differentiable functions with compact support $[0,\pi]$). Clearly, $\widetilde{D}$ is dense in $H^2([0,\pi],\mathbb{R}) \cap H_{0}^1([0,\pi],\mathbb{R})$ under the graph norm, and
		$$\|\mathscr{A}\widehat{G}(\lambda)g\|_{\mathbb{V}} \leq \frac{\Gamma(\nu)}{(\lambda+\beta)^{\nu}} \|g\|_1. $$ 
		Hence, by using the Theorem \ref{Thm-2.2}, the linear system corresponding to \eqref{5.5} admits a resolvent operator $\mathscr{R}(\cdot)$ on $\mathbb{V}$.\\
		Now, the abstract reformulation of equation \eqref{5.5} may be given as follows:
		\begin{equation}\label{5.14}
			\left\{
			\begin{aligned}
				\dfrac{d}{dt}\left(q(t)+\int_0^t N(t-s)q(s)ds\right)&\in \mathscr{A}q(t)+\int_0^t G(t-s)q(s)ds + Bu(t) + H \partial \mathbb{F}(t, q(t)),~ t\in (0,a],\\
				q(0) &\in h(q(\cdot)),
			\end{aligned}
			\right.
		\end{equation}
		where, $h : \mathrm{C}([0,a],\mathbb{V}) \to \mathcal{P}(\mathbb{V})$ is a multivaled map defined in \eqref{def-h} for $p=2$  satisfies all hypotheses of Theorem \ref{MR}, as verified in Example~\ref{ex 3.11}. 
		We now prove that Hypothesis (R) is satisfied. To prove this, let us first calculate the resolvent family corresponding to the linear system of equation \eqref{5.111}.
		Applying the Laplace transform to the resolvent equation (see \eqref{2.3} in Definition \ref{def-res}) and using the convolution theorem:
		\begin{align*}
			\lambda\bigl(\mathtt{I}+\widehat{N}(\lambda)\bigr)\widehat{\mathscr{R}}(\lambda) - \mathtt{I}
			= \mathscr{A}\widehat{\mathscr{R}}(\lambda) + \widehat{G}(\lambda)\widehat{\mathscr{R}}(\lambda).
		\end{align*}
		After rearranging, we get 
		\begin{align}\label{lap-res}
			\widehat{\mathscr{R}}(\lambda) = \big(m(\lambda) \mathtt{I}-\mathscr{A}\big)^{-1},
		\end{align}
		where $m(\lambda)$ is defined by
		\begin{align*}
			m(\lambda)
			= \lambda\bigl(\mathtt{I}+\widehat{N}(\lambda)\bigr) - \widehat{G}(\lambda)
			= \lambda\left(1+\frac{\Gamma(\gamma+1)}{(\lambda+\kappa)^{\gamma+1}}\right)
			- \frac{\alpha}{(\lambda+\beta)^{\nu}}.
		\end{align*}
		By projecting \eqref{lap-res} onto the $k$-th eigenspace and using $\mathscr{A}v_k=-\mu_k v_k$, where $\mu_k=k^2$, we define
		\begin{align}\label{5.83}
			\widehat{r}_k(\lambda):=\frac{1}{m(\lambda)+\mu_k}.
		\end{align}
		As we know, the poles $\{\lambda_{k,j}\}$ of $\widehat{r}_k$ are the roots of the equation  $m(\lambda)+\mu_k=0$. By taking the inverse Laplace transform, the general scalar resolvent formula is
		\begin{equation}\label{5.86}
			r_k(t) 
			= \sum_j \left( \sum_{\ell=0}^{n_j-1} \frac{d^{\ell}_{k,j}}{\ell!} t^\ell \right) e^{\lambda_{k,j} t}, \quad t \geq 0.
		\end{equation}
		where $n_j$ is the multiplicity of the pole $\lambda_{k,j}$ and 
		$$d^{\ell}_{k,j} = \frac{1}{\ell!} \lim_{\lambda\to \lambda_{k,j}}
		\frac{d^\ell}{d\lambda^\ell}
		\left((\lambda - \lambda_{k,j})^{n_j} \widehat{r}_k(\lambda)\right), \quad \ell = 0,1,\ldots,n_j-1 .$$
		From \eqref{5.83}, we have 
		$$\lambda \widehat{r}_k(\lambda) = \frac{\lambda}{m(\lambda) + \mu_k}.$$
		As $\lambda \to \infty$, the dominant term in $m(\lambda)$ is $\lambda$, that is, $m(\lambda) \sim \lambda$ or $\lim\limits_{\lambda\to\infty}\frac{m(\lambda)}{\lambda}=1$.
		Therefore,
		\begin{align}\label{5.87}
			\frac{\lambda}{m(\lambda) + \mu_k} \sim \frac{\lambda}{\lambda + \mu_k} \to 1 \ \text{ as }\ \lambda\to\infty. 
		\end{align}
		Using the well-known property of the Laplace transform 
		$$\lim\limits_{\lambda \to \infty} \lambda \widehat{r}_k(\lambda) =r_k(0^+),$$ 
		where $r_k(0^+)$ is the right-hand limit of $r_k$ at $0$. Thus, combining it with \eqref{5.87}, we deduce that $r_k(0) = 1 .$  Hence, the resolvent family corresponding to the linear system of equation \eqref{5.111} is given by 
		\begin{equation}\label{express-R}
			\mathscr{R}(t)g = \sum_{k=1}^{\infty} r_k(t) \widehat{g}_k v_k,
		\end{equation}
		where $\widehat{g}_k$ and $r_k(t)$ are defined in \eqref{fourier-s} and \eqref{5.86}, respectively. Consequently, this yields $\mathscr{R}(0)=\mathtt{I}$. Observing the expression \eqref{express-R} for $  g^* \in \mathbb{V}^*$ and $t \in [0, a]$, we have
		\begin{equation*}
			\mathscr{R}(t)^*g^* = \sum_{k=1}^{\infty} \overline{r_k(t)} \widehat{g^*_k} v_k.
		\end{equation*}
		Also, the adjoint operator 
		$B^*: \mathbb{V}^*\to \mathbb{U} $ is given by
		$ B^* w^*(\xi) := w^*(\xi).$
		Now, let us estimate 
		$$
		\mathscr{R}(t-s) B B^* \mathscr{R}(a-s)^*g^* =\sum_{k=1}^{\infty} r_k(t-s) \overline{r_k(a-s)} \widehat{g^*_k} v_k.
		$$
		Integrating over $s \in [0,t] $, we obtain
		$$
		R(t)g^*
		= \int_0^t \mathscr{R}(t-s) B B^* \mathscr{R}(a-s)^*g^* ds
		= \sum_{k=1}^{\infty} \int_0^t r_k(t-s) \overline{r_k(a-s)} ds \widehat{g^*_k} v_k.
		$$
		Now, firstly verify that the corresponding linear system \eqref{ineq-14} associated with problem \eqref{5.14} is approximately controllable. To accomplish this, consider	
		$$B^* \mathscr{R}(a-t)^*y^* =  0 \ \text{ for any }\ y^* \in \mathbb{V}^* ,$$
		which further implies $$\sum_{k=1}^{\infty} \overline{r_k(a-t)}\widehat{y^*_k} v_k= 0 \ \text{ for any }\ y^* \in \mathbb{V}^* .$$
		Taking the duality pairing with $v_m $, it follows that
		$$ \overline{r_m(a-t)}\widehat{y^*_m} = 0 \ \text{ for all }\ m \in \mathbb{N}.$$
		Since $\overline{r_m(t)} \not\equiv  0 \text{ on } [0,a],$ we deduce that $$y^*_m =0 \ \text{ for all }\ m \in \mathbb{N}.$$ 
		Hence $y^*=0 \text{ in }  \mathbb{V}^* .$ Therefore, the linear system corresponding to \eqref{5.14} is approximately controllable (see Lemma \ref{lemma 5.1}). Consequently, by \cite[Remark 4.2]{MFPH}, the operator $R(a)$ is injective. Hence, its inverse 
		$$ R(a)^{-1} :\operatorname{Im}R(a) \to \mathbb{V}^* $$
		is well-defined as a linear operator. 
		
		Let $y \in \operatorname{Im}R(a)$. Then there exists $x \in \mathbb{V}$ such that $y = R(a)x$. Thus,
		\begin{equation*}
			y = R(a)x =  \sum_{k=1}^{\infty} \left( \int_0^a |r_k(a-s)|^2 ds \right) \widehat{x}_k v_k= \sum_{k=1}^{\infty} \widehat{y}_k v_k,
		\end{equation*}
		where
		\begin{equation*}
			\widehat{y}_k  = \left( \int_0^a |r_k(a-s)|^2 ds \right) \widehat{x}_k.
		\end{equation*}
		Therefore, the inverse operator $R(a)^{-1}$ is given by
		\begin{equation*}
			R(a)^{-1} y = \sum_{k=1}^{\infty} \frac{1}{\left( \int_0^a |r_k(a-s)|^2 ds \right)} \widehat{y}_k v_k.
		\end{equation*}
		Now, we show that the operator
		$$
		R(t)R(a)^{-1} : \operatorname{Im}R(a) \to \mathbb{V} $$
		is bounded. Let $y \in \operatorname{Im}R(a)$. Then,
		\begin{align}\label{R-bound-1}
			R(t)R(a)^{-1} y 
			= \sum_{k=1}^{\infty}  
			\frac{\displaystyle \int_0^t r_k(t-s)\,\overline{r_k(a-s)} ds}
			{\displaystyle \int_0^a |r_k(a-s)|^2  ds} \, \widehat{y}_k v_k = \sum_{k=1}^{\infty} \rho_k(t,a)\widehat{y}_k v_k,
		\end{align}
		where
		\begin{align}\label{5.961}
			\rho_k(t,a) :=  
			\frac{\displaystyle \int_0^t r_k(t-s)\,\overline{r_k(a-s)} ds}
			{\displaystyle \int_0^a |r_k(a-s)|^2  ds}.
		\end{align}
		Using the change of variable $\varsigma_1= a-s$, we obtain
		\begin{align*}
			\int_0^a |r_k(a-s)|^2 ds
			&= \int_0^a |r_k(\varsigma_1)|^2 d\varsigma_1.
		\end{align*}
		Similarly, letting $\varsigma_2 = t-s$, we have
		\begin{align*}
			\int_0^t r_k(t-s)\,\overline{r_k(a-s)} ds
			&= \int_0^t r_k(\varsigma_2)\,\overline{r_k(a-t+\varsigma_2)} d\varsigma_2.
		\end{align*}
		Therefore, by using Cauchy-Schwarz inequality, \eqref{5.961} becomes
		\begin{align*}
			|\rho_k(t,a)| \leq \frac{\left(\displaystyle \int_0^t |r_k(\varsigma_2)|^2 d\varsigma_2\right)^{1/2}
				\left(\displaystyle  \int_0^t |r_k(a-t+\varsigma_2)|^2 d \varsigma_2 \right)^{1/2}}{\displaystyle \int_0^a |r_k(\varsigma_1)|^2 d\varsigma_1}
		\end{align*}
		Next, using the substitution $\varsigma = a-t+\varsigma_2$, we get
		\begin{align*}
			\int_0^t |r_k(a-t+\varsigma_2)|^2 d \varsigma_2
			&=
			\int_{a-t}^{a} |r_k(\varsigma)|^2 d\varsigma.
		\end{align*}
		Hence,
		\begin{align*}
			|\rho_k(t,a)|
			&\leq
			\frac{
				\left(\displaystyle  \int_0^t |r_k(\varsigma_2)|^2 d\varsigma_2\right)^{1/2}
				\left(\displaystyle  \int_{a-t}^{a} |r_k(\varsigma)|^2 d\varsigma \right)^{1/2}
			}{
				\displaystyle \int_0^a |r_k(\varsigma_1)|^2 d\varsigma_1
			}.
		\end{align*}
		Since
		$$
		\int_0^t |r_k(\varsigma_2)|^2 d \varsigma_2 \leq
		\int_0^a |r_k(\varsigma_1)|^2 d \varsigma_1 \text{ and } \int_{a-t}^{a} |r_k(\varsigma)|^2 d \varsigma
		\le
		\int_0^a |r_k(\varsigma_1)|^2 d \varsigma_1,
		$$
		it follows that
		\begin{align}\label{5.96}
			|\rho_k(t,a)|
			\le 1 \ \text{ for all } k \geq 1, 0 \leq t \leq a.
		\end{align}
		Since $\{v_k\}$ forms an orthonormal basis of $\mathbb{V}$, it follows from \cite[Theorem 7.4.2]{SK} together with \eqref{R-bound-1} and \eqref{5.96} that
		\begin{align*}
			\|R(t)R(a)^{-1} y \|^2_{\mathbb{V}}
			= \sum_{k=1}^{\infty} |\rho_k(t,a)|^2  |\widehat{y}_k|^2 \leq \|y\|^2_{\mathbb{V}}.
		\end{align*}
		Taking the supremum over $\|y\|_{\mathbb{V}} = 1$ completes the proof. Finally, by invoking Theorem \ref{MR}, we conclude that the control system \eqref{5.14} is approximately controllable.
	\end{ex}
	
	\begin{theorem}\label{Theorem-B}
		Let $\mathbb{V}$ be a Banach space with an unconditional Schauder basis 
		$\{v_n\}_{n \geq 1}$. Assume
		$$
		\mathscr{A} v_n=\overline{\lambda}_n v_n \text{ with } \overline{\lambda}_n<0 \text{ and } \mathscr{P}(t) v_n=e^{\overline{\lambda}_n t}v_n $$
		be the associated diagonal operator generated by $\mathscr{A}$. Let $Q \in \mathcal{L}(\mathbb{V})$ be diagonal with respect to the basis, that is $Q v_n= q_n v_n \text{ for } q_n>0$, then the operator
		$$ R(t)R(a)^{-1}:\operatorname{Im}R(a) \to \mathbb{V} $$
		is bounded.
	\end{theorem}
	\begin{proof} Since all the operators are diagonal with respect to the basis 
		$\{v_n\}_{n\geq 1}$, and proceeding similarly as in Example~\ref{5.1111}, we obtain
		\begin{align*}
			R(t)v_n
			= \int_0^t \mathscr{P}(t-s)Q\mathscr{P}(a-s)^*v_n ds
			= q_n \left( \int_0^t e^{\overline{\lambda}_n(a+t-2s)} ds \right)  v_n= q_n \frac{e^{\overline{\lambda}_n(a-t)}-e^{\overline{\lambda}_n(a+t)}}{2|\overline{\lambda}_n|}  v_n.
		\end{align*}
		Proceeding similarly as in Example~\ref{5.1111}, for 
		$$
		y=\sum_{n=1}^{\infty} y_n v_n \in \operatorname{Im}R(a),
		$$
		we estimate
		\begin{align}\label{R-bound1}
			R(t)R(a)^{-1} y 
			= \sum_{n=1}^{\infty}  
			e^{\overline{\lambda}_n(a-t)}\cdot \frac{1-e^{2\overline{\lambda}_nt}}{1-e^{2\overline{\lambda}_n a}} y_n v_n = \sum_{k=1}^{\infty} \overline{\mu}_n(t,a)y_n v_n,
		\end{align}
		where
		\begin{align}\label{5.961}
			\overline{\mu}_n(t,a) := e^{\overline{\lambda}_n(a-t)}\cdot \frac{1-e^{2\overline{\lambda}_nt}}{1-e^{2\overline{\lambda}_n a}}.
		\end{align}
		Because $\overline{\lambda}_n<0$ and $0 \leq t \leq a$, thus we get 
		$$ 0 \leq e^{\overline{\lambda}_n(a-t)} \leq 1 \ \text{ and } \ 0 \leq 
		\frac{1-e^{2\overline{\lambda}_n t}}
		{1-e^{2\overline{\lambda}_n a}} \leq 1, $$
		which further implies that 
		$$ 0\leq \overline{\mu}_n(t,a)\leq 1. $$
		Consequently, by applying \cite[Proposition 1.c.7]{JLLT} together with \eqref{R-bound1}, there exists a constant $C>0$ such that
		$$ \left\|R(t)R(a)^{-1} y \right\|_{\mathbb{V}}=\left\|
		\sum_{n=1}^{\infty}\overline{\mu}_n(t,a)y_n v_n \right\|
		\leq C \sup_n |\overline{\mu}_n(t,a)| \|
		\sum_{n=1}^{\infty} y_n v_n \|.
		$$
		Since
		$$ \sup_n |\overline{\mu}_n(t,a)|\leq 1, $$
		we obtain $$
		\|R(t)R(a)^{-1}\|_{\mathcal{L}(\operatorname{Im}R(a),\mathbb{V})} \leq M_p, $$
		which completes the proof.
	\end{proof}
	
	\begin{ex}\label{Example 5.6}
		Consider the following system:
		\begin{align}\label{Ex 5.6}
			\frac{\partial w(t,\xi)}{\partial t}
			=
			-\mathcal{D}^{2}w(t,\xi)
			+\eta(t,\xi)
			+
			\int_{0}^{1} \mathscr{L}(\xi,\omega) f(t,\omega)d \omega \ \text{ in }\ (0,a] \times [0,1],
		\end{align}
		with the boundary and non-local initial conditions
		\begin{align*}
			&w(t,0)=w(t,\pi)=0, ~ t \in [0,a], \\
			&w(0,\xi)\in h((w(\cdot,\xi)),~ \xi \in [0,1].
		\end{align*}
		where $\eta(\cdot,\cdot), \mathscr{L}(\cdot,\cdot)$, $f$ and non local conditions are defined in the Example \ref{5.1111}. Now, let us denote by $
		\mathfrak{D} $
		the collection of all dyadic intervals in $[0, 1)$, that is, 
		$$
		\mathfrak{D} = \left\{ \left[ \frac{k}{2^j}, \frac{k+1}{2^j} \right) : j \geq 0, \, 0 \leq k \leq 2^j - 1 \right\}.$$
		For each dyadic interval
		$$
		\mathfrak{I} =\left[ \frac{k}{2^j}, \frac{k+1}{2^j} \right),$$
		its length is given by $|\mathfrak{I}|=2^{-j}$ and the integer
		$j = j(\mathfrak{I})$ is called the dyadic scale of $\mathfrak{I}.$
		Moreover, each interval $\mathfrak{I}$ can be decomposed into its left and right halves as
		$$
		\mathfrak{I} = \mathfrak{I}_L \cup \mathfrak{I}_R,$$
		where
		$$
		\mathfrak{I}_L = \left[ \frac{k}{2^j}, \frac{k+1/2}{2^j} \right) \ \text{ and } \ \mathfrak{I}_R = \left[ \frac{k+1/2}{2^j}, \frac{k+1}{2^j} \right).$$
		Next, let
		$$
		\{h_{\mathfrak{I}}\}_{\mathfrak{I} \in \mathfrak{D}}
		$$
		denote the Haar basis (for more information, see \cite[Subsection 4.2.b]{THJML}) associated with the dyadic intervals $\mathfrak{I}\subset[0,1)$, defined by
		$$
		h_{\mathfrak{I}}(x) = |\mathfrak{I}|^{-1/2} (\mathbf{1}_{\mathfrak{I}_L}(x) - \mathbf{1}_{\mathfrak{I}_R}(x)).
		$$
		By \cite[Theorem 2.c.5]{JLLT}, one can see that the Haar basis is an unconditional Schauder basis of $\mathrm{L}^p(0,1)$ for $1 <p <\infty$.
		
		Now, the abstract reformulation of the problem \eqref{Ex 5.6} in $\mathbb{V} = \mathrm{L}^p([0,1],\mathbb{R})$ may be given as follows: 
		\begin{equation}
			q'(t)\in \mathscr{A} q(t)+ Bu(t) + H \partial \mathbb{F}(t, q(t)),~ t\in (0,a],
		\end{equation}
		where the Walsh-Laplacian (dyadic Laplacian) is defined as
		$$ \mathscr{A}= -\mathcal{D}^2 \ \text{ with } \ \mathscr{A} h_\mathfrak{I} = -2^{2 j(\mathfrak{I})} h_\mathfrak{I}.$$
		Equivalently, if
		$g=\sum_{\mathfrak{I}\in \mathfrak{D} }  c_\mathfrak{I} h_\mathfrak{I} $ where $c_\mathfrak{I} := \langle g, h_\mathfrak{I} \rangle,$
		then
		$$
	\mathscr{A}g=-\sum_{\mathfrak{I}\in \mathfrak{D} } 2^{2j(\mathfrak{I})} c_\mathfrak{I} h_\mathfrak{I}$$
		with domain of $\mathscr{A}$ is
		$$
		\mathscr{D}(\mathscr{A})=\left\{
		g=\sum_{\mathfrak{I}\in \mathfrak{D}} c_\mathfrak{I} h_\mathfrak{I} :
		\sum_{\mathfrak{I}\in \mathfrak{D}} 2^{2j(\mathfrak{I})} c_\mathfrak{I} h_\mathfrak{I}
		\text{ converges in } \mathrm{L}^p(0,1)\right\}$$
		and generates the analytic operator
		$$
		\mathscr{R}(t) = e^{t\mathscr{A}},$$
		given by
		$$
		\mathscr{R}(t) h_\mathfrak{I}  = e^{-2^{2j(\mathfrak{I} )}t} h_\mathfrak{I}.$$
		Since all the operators are diagonal with respect to the Haar basis, and proceeding similarly as in Example~\ref{5.1111}, we obtain
		\begin{align*}
			R(t)h_\mathfrak{I} 
			= \int_0^t \mathscr{R}(t-s) B B^* \mathscr{R}(a-s)^*h_\mathfrak{I}  ds
			&= \frac{1}{e^{2 {j(\mathfrak{I})}^l}}\left( \int_0^t e^{-2^{2j(\mathfrak{I})}(a+t-2s)} ds \right)  h_\mathfrak{I}  \\
			&= \frac{e^{-2^{2j(\mathfrak{I})}(a-t)}-e^{-2^{2j(\mathfrak{I})}(a+t)}}{2^{2j(\mathfrak{I})+1} e^{2 {j(\mathfrak{I})}^l}} h_\mathfrak{I} .
		\end{align*}
		Thus, for 
		$$
		y=\sum_{\mathfrak{I}\in \mathfrak{D} } c_\mathfrak{I} h_\mathfrak{I} \in \operatorname{Im}R(a),
		$$
		we estimate
		\begin{align*}
			R(t)R(a)^{-1} y
			= \sum_{\mathfrak{I}\in \mathfrak{D} }  
			e^{-2^{2j(\mathfrak{I})}(a-t)}\cdot \frac{1-e^{-2^{2j(\mathfrak{I})+1}t}}{1-e^{-2^{2j(\mathfrak{I})+1}a}} c_\mathfrak{I} h_\mathfrak{I} =\sum_{\mathfrak{I}\in \mathfrak{D} } \overline{\mu}_{\mathfrak{I}}(t,a)c_\mathfrak{I} h_\mathfrak{I},
		\end{align*}
		where
		\begin{align}\label{5.961}
			\overline{\mu}_{\mathfrak{I}}(t,a) := e^{-2^{2j(\mathfrak{I})}(a-t)}\cdot \frac{1-e^{-2^{2j(\mathfrak{I})+1}t}}{1-e^{-2^{2j(\mathfrak{I})+1}a}}.
		\end{align}
		Since $ 0\leq \overline{\mu}_{\mathfrak{I}}(t,a)\leq 1, $ 
		thus by applying Theorem \ref{Theorem-B}, there exists a unconditional constant $M_p>0$ such that
		$$\|R(t)R(a)^{-1}\|_{\mathcal{L}(\operatorname{Im}R(a),\mathbb{V})} \leq M_p, $$
		where $M_p=\max\{p-1,\frac{1}{p-1}\}$ (see \cite[Theorem 6.1.6]{FANJK}), which completes the proof. Proceeding analogously to Example \ref{5.1111} and applying Theorem \ref{MR}, we deduce that the control system \eqref{Example 5.6} is approximately controllable.
	\end{ex}
	
	\section{Conclusion and Future Work}
	In this paper, we developed the existence of a mild solution of the neutral integro-differential hemivariational inequalities via differential inclusion with non-local inclusion conditions. For this, we employed the theory of generalized Clarke subdifferential, generalized Clarke directional derivative, resolvent operator theory, and fixed-point theorems for multivalued maps. Subsequently, we investigated the approximate controllability of the inclusion system \eqref{PC1} with non-local conditions in a Banach space having a uniformly convex dual. As an application, we have provided an example to determine the approximate controllability of the
	neutral integro-differential inclusions relevant to the theory of heat conduction of material. Further, one can extend the above derived results to fractional order integro-differential and stochastic systems.
	
	A crucial distinction between the control constructions in Mahmudov’s work \cite{NIM} and the recent paper \cite{SV-2025} lies in the issue of adaptedness. In Mahmudov’s formulation, the control is explicitly constructed in a time-local manner, with all stochastic integrals taken over the interval $[0,t]$, ensuring that the control at time $t$ depends only on past information. Consequently, the control process is progressively measurable and belongs to the admissible space $\mathrm{L}^2_{\mathcal{F}}(0, a;\mathbb{U})$, making it well-defined within the stochastic framework. In contrast, the control defined in equation (4.2) of the paper \cite{SV-2025} involves integrals over the entire interval $[0,a]$, including stochastic integrals up to the terminal time $a$. As a result, the control at time $t$ depends on future values of the Wiener process and the state, rendering it anticipative rather than adapted. This violates the standard requirement that admissible stochastic controls be $\mathcal{F}_t$-measurable. Similar issues have also been observed in the literature, see for instance \cite{LLZLMB,YPYKC}. We plan to address and resolve this issue in our future work by developing an adapted control framework consistent with the stochastic setting.

	\medskip\noindent
	\textbf{Declaration of competing interest:}
	The authors have no competing interests to declare.
	
	\medskip\noindent
	\textbf{Acknowledgments:}
	The first author gratefully acknowledges (Prime Minister Research Fellowship, PMRF ID: 2803609) for financial support to carry out her research work. S. Arora also gratefully acknowledges financial support from ANID–Fondecyt (Project No. 3240359). Support for M. T. Mohan's research received from the National Board of Higher Mathematics (NBHM), Department of Atomic Energy, Government of India (Project No. 02011/13/2025/NBHM(R.P)/R\&D II/1137).

\end{document}